%
%
%
%
%
%
\magnification=\magstephalf      
%
%
\vsize=7.5truein                 
\hsize=5.2truein                 
\newskip\stdskip                 
\stdskip=6pt plus3pt minus3pt    
\medskipamount=\stdskip          
\parindent=0pt                   
\parskip=\stdskip                
\abovedisplayskip=\stdskip       
\belowdisplayskip=\stdskip       
\mathsurround=0.75pt             
\overfullrule=0pt                
%
%
\def\ppar{\par\goodbreak\vskip 8pt plus 4pt minus 4pt}     
%
%
\def\stdspace{\hskip 0.75em plus 0.15em\ignorespaces}
\let\qua\stdspace 
%
%
%
%
%
%
%
\def\hexnumber#1{\ifcase#1 0\or 1\or 2\or 3\or 4\or 5\or 6\or 7\or 8\or
 9\or A\or B\or C\or D\or E\or F\fi}
%
%
\font\thirtnmsa=msam10 scaled 1315    
\font\tenmsa=msam10          \font\ninemsa=msam9
\font\sevenmsa=msam7         \font\sixmsa=msam6
\font\fivemsa=msam5
%
%
\newfam\msafam                  \textfont\msafam=\tenmsa
\scriptfont\msafam=\sevenmsa    \scriptscriptfont\msafam=\fivemsa
\edef\hexa{\hexnumber\msafam}        
\def\msa{\fam\msafam\tenmsa}         
%
%
\font\thirtnmsb=msbm10 scaled 1315   
\font\tenmsb=msbm10      \font\ninemsb=msbm9
\font\sevenmsb=msbm7     \font\sixmsb=msbm6
\font\fivemsb=msbm5
%
\newfam\msbfam                   \textfont\msbfam=\tenmsb       
\scriptfont\msbfam=\sevenmsb     \scriptscriptfont\msbfam=\fivemsb
\edef\hexb{\hexnumber\msbfam}    
\def\msb{\fam\msbfam\tenmsb}     
%
%
\font\thirtneufm=eufm10 scaled 1315   
\font\teneufm=eufm10                 \font\nineeufm=eufm9
\font\seveneufm=eufm7                \font\sixeufm=eufm6
\font\fiveeufm=eufm5
%
\newfam\eufmfam                    \textfont\eufmfam=\teneufm
\scriptfont\eufmfam=\seveneufm     \scriptscriptfont\eufmfam=\fiveeufm
\edef\hexf{\hexnumber\eufmfam}      
\def\frak{\fam\eufmfam\teneufm}     
%
%
%
\font\thirtnrm=cmr10 scaled 1315    
\font\ninerm=cmr9                   \font\sixrm=cmr6   
%
\font\thirtni=cmmi10 scaled 1315    
\font\ninei=cmmi9                   \font\sixi=cmmi6  
%
\font\thirtnsy=cmsy10 scaled 1315   
\font\ninesy=cmsy9                  \font\sixsy=cmsy6  
%
\font\thirtnbf=cmbx10 scaled 1315   
\font\ninebf=cmbx9                  \font\sixbf=cmbx6  
%
%
\font\thirtnex=cmex10 scaled 1315   
\font\nineex=cmex9                  
%
%
\font\thirtnit=cmti10 scaled 1315  
\font\nineit=cmti9                  
%
\font\thirtnsl=cmsl10 scaled 1315  
\font\ninesl=cmsl9                  
%
\font\thirtntt=cmtt10 scaled 1315  
\font\ninett=cmtt9                  
%
%
%
%
\def\small{%
%
%
\textfont0=\ninerm \scriptfont0=\sixrm \scriptscriptfont0=\fiverm
\def\rm{\fam0\ninerm}
%
%
\textfont1=\ninei \scriptfont1=\sixi \scriptscriptfont1=\fivei
%
%
\textfont2=\ninesy \scriptfont2=\sixsy \scriptscriptfont2=\fivesy
%
%
\textfont3=\nineex \scriptfont3=\nineex \scriptscriptfont3=\nineex
%
%
\textfont\bffam=\ninebf \scriptfont\bffam=\sixbf
\scriptscriptfont\bffam=\fivebf \def\bf{\fam\bffam\ninebf}%
%
%
\textfont\itfam=\nineit \def\it{\fam\itfam\nineit}%
\textfont\slfam=\ninesl \def\sl{\fam\slfam\ninesl}%
\textfont\ttfam=\ninett \def\tt{\fam\ttfam\ninett}%
%
%
%
\textfont\msafam=\ninemsa \scriptfont\msafam=\sixmsa
\scriptscriptfont\msafam=\fivemsa \def\msa{\fam\msafam\ninemsa}%
%
%
\textfont\msbfam=\ninemsb \scriptfont\msbfam=\sixmsb
\scriptscriptfont\msbfam=\fivemsb \def\msb{\fam\msbfam\ninemsb}%
%
%
\textfont\eufmfam=\nineeufm  \scriptfont\eufmfam=\sixeufm
\scriptscriptfont\eufmfam=\fiveeufm \def\frak{\fam\eufmfam\nineeufm}%
%
%
%
\normalbaselineskip=11pt%
\setbox\strutbox=\hbox{\vrule height8pt depth3pt width0pt}%
%
%
\normalbaselines\rm
%
%
\stdskip=4pt plus2pt minus2pt    
\medskipamount=\stdskip          
\parskip=\stdskip                
\abovedisplayskip=\stdskip       
\belowdisplayskip=\stdskip       
\def\ppar{\par\goodbreak\vskip 6pt plus 3pt minus 3pt}%
%
%
\def\section##1{\global\advance\sectionnumber by 1
\vskip-\lastskip\penalty-800\vskip 20pt plus10pt minus5pt 
\egroup{\bf\number\sectionnumber\quad##1}\bgroup\small         
\vskip 6pt plus3pt minus3pt
\nobreak\resultnumber=1}
}    
%
\def\beginsmall{\bgroup\small}
\let\endsmall\egroup
%
%
%
%
\def\large{%
\textfont0=\thirtnrm \scriptfont0=\ninerm \scriptscriptfont0=\sevenrm
\def\rm{\fam0\thirtnrm}%
\textfont1=\thirtni \scriptfont1=\ninei \scriptscriptfont1=\seveni
\textfont2=\thirtnsy \scriptfont2=\ninesy \scriptscriptfont2=\sevensy
\textfont3=\thirtnex \scriptfont3=\thirtnex \scriptscriptfont3=\thirtnex
\textfont\bffam=\thirtnbf \scriptfont\bffam=\ninebf
\scriptscriptfont\bffam=\sevenbf \def\bf{\fam\bffam\thirtnbf}%
\textfont\itfam=\thirtnit \def\it{\fam\itfam\thirtnit}%
\textfont\slfam=\thirtnsl \def\sl{\fam\slfam\thirtnsl}%
\textfont\ttfam=\thirtntt \def\tt{\fam\ttfam\thirtntt}%
\textfont\msafam=\thirtnmsa \scriptfont\msafam=\ninemsa
\scriptscriptfont\msafam=\sevenmsa \def\msa{\fam\msafam\thirtnmsa}%
\textfont\msbfam=\thirtnmsb \scriptfont\msbfam=\ninemsb
\scriptscriptfont\msbfam=\sevenmsb \def\msb{\fam\msbfam\thirtnmsb}%
\textfont\eufmfam=\thirtneufm  \scriptfont\eufmfam=\nineeufm
\scriptscriptfont\eufmfam=\seveneufm \def\frak{\fam\eufmfam\teneufm}%
\normalbaselineskip=16pt%
\setbox\strutbox=\hbox{\vrule height11.5pt depth4.5pt width0pt}%
\normalbaselines\rm}%
\let\Large\large   
%
\def\Bbb#1{{\msb#1}}

%

%
\mathchardef\plussquare="0\hexa01
\mathchardef\nge="3\hexb0B
\mathchardef\maltesecross="0\hexa7A
\mathchardef\del="0\hexf01
%
%
%
%
\font\sc=cmcsc10
%
%
%
%
\def\sqr#1#2{{\vcenter{\vbox{\hrule  height.#2truept
	\hbox{\vrule width.#2truept height#1truept 
	\kern#1truept \vrule width.#2truept}
	\hrule height.#2truept}}}}
\def\sq{\sqr55}    
%
%
%
%
\newcount\sectionnumber            
\newcount\resultnumber             
\sectionnumber=0\resultnumber=1    
%
%
%
\def\section#1{\global\advance\sectionnumber by 1
\xdef\nextkey{\number\sectionnumber}
\vskip-\lastskip\penalty-800\vskip 20pt plus10pt minus5pt 
{\large\bf\number\sectionnumber\quad#1}         
\vskip 8pt plus4pt minus4pt
\nobreak\resultnumber=1}      
%
%
%
%
%
\def\sh#1{\vskip-\lastskip\ppar{\bf #1}\par\nobreak\medskip}         
%
%
%
%

%
\def\proc#1{\xdef\nextkey{\number\sectionnumber.\number\resultnumber}%
\vskip-\lastskip\ppar\bf%
\noindent#1\ \number\sectionnumber.\number\resultnumber
\stdspace\sl\global\advance\resultnumber by 1\ignorespaces}
 
%
%
\def\qed{\hfill$\sq$\par\goodbreak\rm}   
%
%
%
%
%
%
%
%
\def\proclaim#1{\vskip-\lastskip\ppar\bf%
\noindent#1\stdspace\sl\ignorespaces} 
\let\endproclaim\endproc
%
%
%
%
\def\rk#1{\vskip-\lastskip\ppar{\bf #1}\stdspace\ignorespaces}                

%
%
%
%
%
%
\def\label{\xdef\nextkey{\number\sectionnumber.\number\resultnumber}%
\number\sectionnumber.\number\resultnumber
\global\advance\resultnumber by 1}
%
%
%
%
%
%
%
%
%
%
%
%
%
%
%
%
\newcount\refnumber              
\refnumber=1                     
\long\def\reflist#1\endreflist{%
\long\def\thereflist{#1}{\def\refkey##1##2\par{\xdef##1{\number\refnumber}%
\global\advance\refnumber by 1}%
\def\key##1##2\par{\expandafter\xdef%
\csname##1\endcsname{\number\refnumber}%
\global\advance\refnumber by 1}#1\par}}
\long\def\references{%
\penalty-800\vskip-\lastskip\vskip 15pt plus10pt minus5pt 
{\large\bf References}\ppar 
{\leftskip=25pt\frenchspacing    
\small\parskip=3pt plus2pt       
\def\refkey##1##2\par{\noindent  
\llap{[##1]\stdspace}\ignorespaces##2\par}         
\def\key##1##2\par{\noindent  
\llap{[\ref{##1}]\stdspace}\ignorespaces##2\par}  
\def\,{\thinspace}\thereflist\par}}
%
%
%
\newcount\footnotenumber         
\footnotenumber=1                
\def\fnote#1{\xdef\nextkey{\number\footnotenumber}%
{\small\ifnum\footnotenumber>9\parindent=14pt%
\else\parindent=10pt\fi\footnote{$^{\number\footnotenumber}$}%
{\hglue-5pt#1}\global\advance\footnotenumber by 1}}
%
%
%
%
%
%
%
\newcount\figurenumber          
\figurenumber=1                 
\def\caption#1{\xdef\nextkey{\number\figurenumber}%
\cl{\small Figure \number\figurenumber: #1}%
\global\advance\figurenumber by 1}
\def\figurelabel{\xdef\nextkey{\number\figurenumber}%
\cl{\small Figure \number\figurenumber}%
\global\advance\figurenumber by 1}
\long\def\figure#1\endfigure{{\xdef\nextkey{\number\figurenumber}%
\let\captiontext\relax\def\caption##1{\xdef\captiontext{##1}}%
\midinsert\cl{\ignorespaces#1\unskip\unskip\unskip\unskip}\vglue6pt\cl{\small 
Figure \number\figurenumber\ifx\captiontext\relax\else: \captiontext
\fi}\endinsert\global\advance\figurenumber by 1}}
%
%
%
%
%
%
%
\def\nextkey{??}   
%
\def\key#1{\expandafter\xdef\csname #1\endcsname{\nextkey}}
\def\ref#1{\expandafter\ifx\csname #1\endcsname\relax
\immediate\write16{Reference {#1} undefined}??\else
\csname #1\endcsname\fi}
%
%
%
%
%
%
%
\newread\gtinfile
\newwrite\gtreffile
\def\useforwardrefs{
\openin\gtinfile\jobname.ref
\ifeof\gtinfile
\closein\gtinfile
\immediate\write16{No file \jobname.ref}
\else
\closein\gtinfile
\input \jobname.ref
\fi
\immediate\openout\gtreffile \jobname.ref
%
%
\def\key##1{{\def\\{\noexpand}%
\expandafter\xdef\csname ##1\endcsname{\nextkey}%
\immediate\write\gtreffile{\\\expandafter\\\def\\\csname ##1\\\endcsname%
{\nextkey}}}}
%
%
\long\def\reflist##1\endreflist{%
\long\def\thereflist{##1}{\def\refkey####1####2\par{\xdef####1{%
\number\refnumber}{\def\\{\noexpand}\immediate\write\gtreffile
{\\\def\\####1{\number\refnumber}}}\global\advance\refnumber by 1}%
\def\key####1####2\par{\expandafter\xdef%
\csname####1\endcsname{\number\refnumber}%
{\def\\{\noexpand}\immediate\write\gtreffile
{\\\expandafter\\\def\\\csname ####1\\\endcsname{\number\refnumber}}}
\global\advance\refnumber by 1}##1\par}}
\long\def\biblio##1\endbiblio{\reflist##1\endreflist\references}%
%
%
\def\numkey##1{{\def\\{\noexpand}%
\xdef##1{\number\sectionnumber.\number\resultnumber}
\immediate\write\gtreffile{\\\def\\##1%
{\number\sectionnumber.\number\resultnumber}}}}
\def\seckey##1{{\def\\{\noexpand}\xdef##1{\number\sectionnumber}
\immediate\write\gtreffile{\\\def\\##1{\number\sectionnumber}}}}
\def\figkey##1{\xdef##1{\number\figurenumber}%
{\def\\{\noexpand}\immediate\write\gtreffile%
{\\\def\\##1{\number\figurenumber}}}
\number\figurenumber\global\advance\figurenumber by 1}
}   
%
%
%
%
\def\figkey#1{\xdef#1{\number\figurenumber}%
\number\figurenumber\global\advance\figurenumber by 1}
\def\fig#1#2\endfig{%
\midinsert\cl{#2}\vglue6pt\cl{\small Figure #1}\endinsert}
\def\newfig{\number\figurenumber\global\advance\figurenumber by 1}
\def\numkey#1{\xdef#1{\number\sectionnumber.\number\resultnumber}}
\def\seckey#1{\xdef#1{\number\sectionnumber}}
%
%
%
%
%
%
%
%
%
\def\verb{\catcode`\"=\active}       
\def\brev{\catcode`\"=12}            
\brev                                
\verb                                
{\obeyspaces\gdef {\ }}              
{\catcode`\`=\active\gdef`{\relax\lq}}
\def"{%
\begingroup\baselineskip=12pt\def\par{\leavevmode\endgraf}%
\tt\obeylines\obeyspaces\parskip=0pt\parindent=0pt%
\catcode`\$=12\catcode`\&=12\catcode`\^=12\catcode`\#=12%
\catcode`\_=12\catcode`\~=12%
\catcode`\{=12\catcode`\}=12\catcode`\%=12\catcode`\\=12%
\catcode`\`=\active\let"\endgroup}
\brev      
%
%
%
%
%
%
\def\items{\par\leftskip = 25pt}           
\def\enditems{\par\leftskip = 0pt}         
\def\item#1{\par\leavevmode\llap{#1\stdspace}%
\ignorespaces}                             
%
%

%
%
\def\np{\vfil\eject}         
\def\nl{\hfil\break}         
\def\cl{\centerline}         
\def\gt{{\mathsurround=0pt\it $\cal G\mskip-2mu$eometry \&\ 
$\cal T\!\!$opology}}        
\def\agt{{\mathsurround=0pt\it$\cal A\mskip-.7mu$lgebraic \&\ 
$\cal G\mskip-2mu$eometric $\cal T\!\!$opology}}  
%
%
%

%
%
%
%
%
\def\title#1{\def\thetitle{#1}}

\def\author#1{\edef\previousauthors{\theauthors}
 \ifx\theauthors\relax\def\theauthors{#1}\else
 \def\theauthors{\previousauthors\par#1}\fi}

%
\def\address#1{\edef\previousaddresses{\theaddress}
 \ifx\theaddress\relax\def\theaddress{#1}\else
 \def\theaddress{\previousaddresses\par\vskip 2pt\par#1}\fi}
\def\secondaddress#1{\edef\previousaddresses{\theaddress}
 \ifx\theaddress\relax\def\theaddress{#1}\else
 \def\theaddress{\previousaddresses\par{\rm and}\par#1}\fi}   

\def\email#1{\edef\previousemails{\theemail}
 \ifx\theemail\relax\def\theemail{#1}\else
 \def\theemail{\previousemails\hskip 0.75em\relax#1}\fi}
\def\secondemail#1{\edef\previousemails{\theemail}
 \ifx\theemail\relax\def\theemail{#1}\else
 \def\theemail{\previousemails\hskip 0.75em{\rm and}\hskip 0.75em
 \relax#1}\fi}
\def\url#1{\edef\previousurls{\theurl}
 \ifx\theurl\relax\def\theurl{#1}\else
 \def\theurl{\previousurls\hskip 0.75em\relax#1}\fi}
\def\secondurl#1{\edef\previousurls{\theurl}
 \ifx\theurl\relax\def\theurl{#1}\else
 \def\theurl{\previousurls\hskip 0.75em{\rm and}\hskip 0.75em
 \relax#1}\fi}
\long\def\abstract#1\endabstract{\long\def\theabstract{#1}}
\def\primaryclass#1{\def\theprimaryclass{#1}}
\let\subjclass\primaryclass                        
\def\secondaryclass#1{\def\thesecondaryclass{#1}}
\def\keywords#1{\def\thekeywords{#1}}
%
%
\let\\\par\let\thetitle\relax\let\theshorttitle\relax
\let\theauthors\relax\let\theshortauthors\relax
\let\theaddress\relax\let\theshortaddress\relax
\let\theemail\relax\let\theurl\relax
\let\theabstract\relax\let\theprimaryclass\relax
\let\thesecondaryclass\relax\let\thekeywords\relax
%
%
%
%
\long\def\maketitlepage{    

\vglue 0.2truein   

%
{\parskip=0pt\leftskip 0pt plus 1fil\def\\{\par\smallskip}{\large
\bf\thetitle}\par\medskip}   

\vglue 0.15truein 

%
{\parskip=0pt\leftskip 0pt plus 1fil\def\\{\par}{\sc\theauthors}
\par\medskip}%
 
\vglue 0.1truein 

%
{\small\parskip=0pt
{\leftskip 0pt plus 1fil\def\\{\par}{\sl\theaddress}\par}
\ifx\theemail\relax\else  
\vglue 5pt \def\\{\stdspace{\rm and}\stdspace} 
\cl{Email:\stdspace\tt\theemail}\fi
\ifx\theurl\relax\else    
\vglue 5pt \def\\{\stdspace{\rm and}\stdspace} 
\cl{URL:\stdspace\tt\theurl}\fi\par}

\vglue 7pt 

{\bf Abstract}

\vglue 5pt

\theabstract

\vglue 7pt 

{\bf AMS Classification numbers}\quad Primary:\quad \theprimaryclass\par

Secondary:\quad \thesecondaryclass

\vglue 5pt 

{\bf Keywords:}\quad \thekeywords

\np  

}    
%
%
\long\def\makeshorttitle{    


%
{\parskip=0pt\leftskip 0pt plus 1fil\def\\{\par\smallskip}{\large
\bf\thetitle}\par\medskip}   

\vglue 0.05truein 

%
{\parskip=0pt\leftskip 0pt plus 1fil\def\\{\par}{\sc\theauthors}
\par\medskip}%
 
\vglue 0.03truein 

%
{\small\parskip=0pt
{\leftskip 0pt plus 1fil\def\\{\par}{\sl\ifx\theshortaddress\relax
\theaddress\else\theshortaddress\fi}\par}
\ifx\theemail\relax\else  
\vglue 5pt \def\\{\stdspace{\rm and}\stdspace} 
\cl{Email:\stdspace\tt\theemail}\fi
\ifx\theurl\relax\else    
\vglue 5pt \def\\{\stdspace{\rm and}\stdspace} 
\cl{URL:\stdspace\tt\theurl}\fi\par}

\vglue 10pt 


{\small\leftskip 25pt\rightskip 25pt{\bf Abstract}\stdspace\theabstract

{\bf AMS Classification}\stdspace\theprimaryclass
\ifx\thesecondaryclass\relax\else; \thesecondaryclass\fi\par
{\bf Keywords}\stdspace \thekeywords\par}
\vglue 7pt
}    
\let\maketitle\makeshorttitle        
%
%

\def\volumenumber#1{\def\thevolumenumber{#1}}
\def\volumeyear#1{\def\thevolumeyear{#1}}
\def\pagenumbers#1#2{\def\startpage{#1}\def\finishpage{#2}}
\def\published#1{\def\publishdate{#1}}
\def\received#1{\def\receiveddate{#1}}
\def\revised#1{\def\reviseddate{#1}}
\let\reviseddate\relax
\volumenumber{X}
\volumeyear{20XX}
\pagenumbers{1}{XXX}
\published{XX Xxxember 20XX}

\long\def\makeagttitle{   
\agt\hfill      
\hbox to 60truept{\vbox to 0pt{\vglue -14truept{\bf [Logo here]}\vss}\hss}
\break
{\small Volume \thevolumenumber\ (\thevolumeyear)
\startpage--\finishpage\nl
Published: \publishdate}

\vglue .2truein

{\parskip=0pt\leftskip 0pt plus 1fil\def\\{\par\smallskip}{\large
\bf\thetitle}\par\medskip}   
\vglue 0.05truein 

%
{\parskip=0pt\leftskip 0pt plus 1fil\def\\{\par}{\sc\theauthors}
\par\medskip}%
 
\vglue 0.03truein 


{\small\leftskip 25truept\rightskip 25truept{\bf Abstract}\stdspace\theabstract

{\bf AMS Classification}\stdspace\theprimaryclass
\ifx\thesecondaryclass\relax\else; \thesecondaryclass\fi\par
{\bf Keywords}\stdspace \thekeywords\par}\vglue 7truept

}   


\def\Addresses{\bigskip
{\small \parskip 0pt \leftskip 0pt \rightskip 0pt plus 1fil \def\\{\par}
\sl\theaddress\par\medskip \rm Email:\stdspace\tt\theemail\par
\ifx\theurl\relax\else\smallskip \rm URL:\stdspace\tt\theurl\par\fi}}

\def\agtart{
\hoffset 14truemm
\voffset 31truemm
\font\phead=cmsl9 scaled 950
\font\pnum=cmbx10 scaled 913
\font\pfoot=cmsl9 scaled 950
\headline{\vbox to 0pt{\vskip -4.5mm\line{\small\phead\ifnum
\count0=\startpage ISSN numbers are printed here
\hfill {\pnum\folio}\else\ifodd\count0\def\\{ }%
\ifx\theshorttitle\relax\thetitle\else\theshorttitle\fi\hfill{\pnum\folio}
\else\def\\{ and }{\pnum\folio}\hfill\ifx\theshortauthors\relax\theauthors
\else\theshortauthors\fi\fi\fi}\vss}}
\footline{\vbox to 0pt{\vglue 0mm\line{\small\pfoot\ifnum\count0=\startpage
Copyright declaration is printed here\hfill\else
\agt, Volume \thevolumenumber\ (\thevolumeyear)\hfill\fi}\vss}}
\let\maketitle\makeagttitle\let\makeshorttitle\makeagttitle}


\def\ifplaintex{\expandafter\ifx\csname documentclass\endcsname\relax}


\ifplaintex 
\hoffset 14truemm
\voffset 31truemm
\else
\headsep 23pt
\footskip 35pt
\hoffset -4truemm
\voffset 12.5truemm
\fi

\expandafter\ifx\csname beginpicture\endcsname\relax
\expandafter\ifx\csname documentclass\endcsname\relax
\input pictex \else\font\fiverm=cmr5
\input prepictex \input pictex \input postpictex \fi\fi

\def\gt{{\mathsurround=0pt\it $\cal G\mskip-2mu$eometry \&\ 
$\cal T\!\!$opology}}        

\def\gtp{{\mathsurround=0pt\it $\cal G\mskip-2mu$eometry \&\ 
$\cal T\!\!$opology $\cal P\!$ublications}}  


\def\lognumber#1{\def\thelognumber{#1}}
\def\volumenumber#1{\def\thevolumenumber{#1}}
\def\papernumber#1{\def\thepapernumber{#1}}
\def\volumeyear#1{\def\thevolumeyear{#1}}

\def\pagenumbers#1#2{\def\startpage{#1}\def\finishpage{#2}}
\def\published#1{\def\publishdate{#1}}
\def\proposed#1{\def\theproposer{#1}}
\def\seconded#1{\def\theseconders{#1}}
\def\received#1{\def\receiveddate{#1}}
\def\revised#1{\def\reviseddate{#1}}
\def\accepted#1{\def\accepteddate{#1}}

\long\def\asciiabstract#1{\long\def\theasciiabstract{#1}}


\let\\\par\let\thelognumber\relax
\let\thevolumenumber\relax\let\thepapernumber\relax
\let\thevolumeyear\relax\let\thesamplenumber\relax\let\startpage\relax
\let\finishpage\relax\let\publishdate\relax\let\receiveddate\relax
\let\reviseddate\relax\let\accepteddate\relax\let\theasciititle\relax
\let\theasciiauthors\relax
\let\theasciiabstract\relax
\let\theasciiemail\relax\let\theshortauthors\relax\let\theshorttitle\relax

\long\def\maketitlep{   

\count0=\startpage

\gt\hfill      
\beginpicture
\setcoordinatesystem units <0.33truein, 0.33truein> point at 2.2 0.9
\setplotsymbol ({$\cal G$})
\plotsymbolspacing=9truept
\circulararc 315 degrees from 0 1 center at 0 0
\setplotsymbol ({$\cal T$})
\circulararc 315 degrees from 1 -1 center at 1 0
\endpicture
%
\break
{\small\ifx\thesamplenumber\relax 
Volume \else Sample
\fi\thevolumenumber\ (\thevolumeyear)
\startpage--\finishpage\nl
Published: \publishdate}
\vglue 0.5truein plus 0.4fil minus 0.1truein

{\parskip=0pt\leftskip 0pt plus 1fil\def\\{\par\smallskip}{\ifplaintex\large
\else\Large\fi\bf\thetitle}\par\medskip}   

\vglue 0pt plus 0.1fil 

{\parskip=0pt\leftskip 0pt plus 1fil\def\\{\par}{\sc\theauthors}
\par\medskip}

\vglue 0pt plus 0.1fil 

{\small\parskip=0pt\let\newline\\
{\leftskip 0pt plus 1fil\def\\{\par}{\sl\theaddress}\par}
\expandafter\ifx\theemail\relax    
\relax\else\vglue 5pt plus 0.02fil minus 2pt\def\\{\stdspace{\rm 
and}\stdspace} 
\cl{Email:\stdspace\tt\theemail}\fi
\ifx\theurl\relax                  
\relax\else\vglue 5pt plus 0.02fil minus 2pt\def\\{\stdspace{\rm 
and}\stdspace}
\cl{URL:\stdspace\tt\theurl}\fi\par}

\vglue 7pt plus 0.3fil minus 3pt

{\bf Abstract}
\vglue 5pt plus 0.1fil minus 2pt

\theabstract

\vglue 7pt plus 0.3fil minus 3pt

{\bf AMS Classification numbers}\quad Primary:\quad \theprimaryclass

Secondary:\quad \thesecondaryclass

\vglue 5pt plus 0.3fil minus 2pt

{\bf Keywords}\quad \thekeywords

\vglue 10pt plus 0.5fil minus 5pt

{\small  Proposed: \theproposer\hfill Received: \receiveddate\nl
Seconded: \theseconders\hfill 
\ifx\reviseddate\relax                         
Accepted: \accepteddate                        
\else
Revised: \reviseddate                          
\fi}
\eject
}       

\let\maketitlepage\maketitlep
\let\maketitle\maketitlepage


\font\phead=cmsl9 scaled 950
\font\lhead=cmsl9 scaled 1050
\font\pnum=cmbx10 scaled 913
\font\lnum=cmbx10 
\font\pfoot=cmsl9 scaled 950
\font\lfoot=cmsl9 scaled 1050
\ifplaintex
\headline{\vbox to 0pt{\vskip -4.5mm\line{\small\phead\ifnum
\count0=\startpage ISSN 1364-0380 (on line)
1465-3060 (printed) \hfill {\pnum\folio}\else\ifodd\count0\def\\{ }%
\ifx\theshorttitle\relax\thetitle\else\theshorttitle\fi\hfill{\pnum\folio}
\else\def\\{ and }{\pnum\folio}\hfill\ifx\theshortauthors\relax\theauthors
\else\theshortauthors\fi\fi\fi}\vss}}
\footline{\vbox to 0pt{\vglue 0mm\line{\small\pfoot\ifnum\count0=\startpage
\copyright\ \gtp\hfill\else
\gt, Volume \thevolumenumber\ (\thevolumeyear)\hfill\fi}\vss
}}
\else
\makeatletter
\def\@oddhead{{\small\lhead\ifnum\count0=\startpage ISSN 1364-0380 (on line)
1465-3060 (printed) \hfill {\lnum\number\count0}\else\ifodd\count0
\def\\{ }\ifx\theshorttitle\relax \thetitle \else\theshorttitle\fi\hfill
{\lnum\number\count0}\else\def\\{ and }{\lnum\number\count0}
\hfill\ifx\theshortauthors\relax 
\theauthors\else\theshortauthors\fi\fi\fi}}\def\@evenhead{\@oddhead}
\def\@oddfoot{\small\lfoot\ifnum\count0=\startpage\copyright\ \gtp\hfill\else
\gt, Volume \thevolumenumber\ (\thevolumeyear)\hfill\fi}
\def\@evenfoot{\@oddfoot}
\makeatother
\fi


\newwrite\gtoutfile
\long\gdef\makeheadfile{  
{\def\\{, }\def\s{ }
\immediate\openout\gtoutfile head.xxx
\immediate\write\gtoutfile{To: math@arxiv.org}
\immediate\write\gtoutfile{Subject: put or rep NNNNN:pppp}
\immediate\write\gtoutfile{--text follows this line--}
\immediate\write\gtoutfile{Proxy-for: \ifx\theasciiauthors\relax
\theauthors\else\theasciiauthors\fi\s<\ifx\theasciiemail\relax\theemail\else\theasciiemail\fi>}
\immediate\write\gtoutfile{\noexpand\\}
\immediate\write\gtoutfile{Authors: \ifx\theasciiauthors\relax
\theauthors\else\theasciiauthors\fi}
\immediate\write\gtoutfile{Title: \ifx\theasciititle\relax
\thetitle\else\theasciititle\fi}
\immediate\write\gtoutfile{Subj-class: GT or SG or MG etc}
\immediate\write\gtoutfile{MSC-class: \theprimaryclass\ifx\thesecondaryclass\relax\else, \thesecondaryclass\fi}
\immediate\write\gtoutfile{Journal-ref: Geom. Topol. \thevolumenumber
(\thevolumeyear) \startpage-\finishpage}
\immediate\write\gtoutfile{Comments: Published by Geometry and Topology at}
\immediate\write\gtoutfile{\s\s http://www.maths.warwick.ac.uk/gt/GTVol\thevolumenumber/paper\thepapernumber.abs.html}
\immediate\write\gtoutfile{\noexpand\\}
\immediate\write\gtoutfile{}
\ifx\theasciiabstract\relax
\immediate\write\gtoutfile{\theabstract}\else
\immediate\write\gtoutfile{\theasciiabstract}\fi
\immediate\write\gtoutfile{}
\immediate\write\gtoutfile{\noexpand\\}
\immediate\write\gtoutfile{}
\immediate\closeout\gtoutfile}}  

\def\maketitlepage{\maketitlep\makeheadfile}
\let\maketitle\maketitlepage


\def\ifplaintex{\expandafter\ifx\csname documentclass\endcsname\relax}


\ifplaintex 
\hoffset 14truemm
\voffset 31truemm
\else
\headsep 23pt
\footskip 35pt
\hoffset -4truemm
\voffset 12.5truemm
\fi

\expandafter\ifx\csname beginpicture\endcsname\relax
\expandafter\ifx\csname documentclass\endcsname\relax
\input pictex \else\font\fiverm=cmr5
\input prepictex \input pictex \input postpictex \fi\fi

\def\gt{{\mathsurround=0pt\it $\cal G\mskip-2mu$eometry \&\ 
$\cal T\!\!$opology}}        

\def\gtp{{\mathsurround=0pt\it $\cal G\mskip-2mu$eometry \&\ 
$\cal T\!\!$opology $\cal P\!$ublications}}  


\def\lognumber#1{\def\thelognumber{#1}}
\def\volumenumber#1{\def\thevolumenumber{#1}}
\def\papernumber#1{\def\thepapernumber{#1}}
\def\volumeyear#1{\def\thevolumeyear{#1}}

\def\pagenumbers#1#2{\def\startpage{#1}\def\finishpage{#2}}
\def\published#1{\def\publishdate{#1}}
\def\proposed#1{\def\theproposer{#1}}
\def\seconded#1{\def\theseconders{#1}}
\def\received#1{\def\receiveddate{#1}}
\def\revised#1{\def\reviseddate{#1}}
\def\accepted#1{\def\accepteddate{#1}}

\long\def\asciiabstract#1{\long\def\theasciiabstract{#1}}


\let\\\par\let\thelognumber\relax
\let\thevolumenumber\relax\let\thepapernumber\relax
\let\thevolumeyear\relax\let\thesamplenumber\relax\let\startpage\relax
\let\finishpage\relax\let\publishdate\relax\let\receiveddate\relax
\let\reviseddate\relax\let\accepteddate\relax\let\theasciititle\relax
\let\theasciiauthors\relax
\let\theasciiabstract\relax
\let\theasciiemail\relax\let\theshortauthors\relax\let\theshorttitle\relax

\long\def\maketitlep{   

\count0=\startpage

\gt\hfill      
\beginpicture
\setcoordinatesystem units <0.33truein, 0.33truein> point at 2.2 0.9
\setplotsymbol ({$\cal G$})
\plotsymbolspacing=9truept
\circulararc 315 degrees from 0 1 center at 0 0
\setplotsymbol ({$\cal T$})
\circulararc 315 degrees from 1 -1 center at 1 0
\endpicture
%
\break
{\small\ifx\thesamplenumber\relax 
Volume \else Sample
\fi\thevolumenumber\ (\thevolumeyear)
\startpage--\finishpage\nl
Published: \publishdate}
\vglue 0.5truein plus 0.4fil minus 0.1truein

{\parskip=0pt\leftskip 0pt plus 1fil\def\\{\par\smallskip}{\ifplaintex\large
\else\Large\fi\bf\thetitle}\par\medskip}   

\vglue 0pt plus 0.1fil 

{\parskip=0pt\leftskip 0pt plus 1fil\def\\{\par}{\sc\theauthors}
\par\medskip}

\vglue 0pt plus 0.1fil 

{\small\parskip=0pt\let\newline\\
{\leftskip 0pt plus 1fil\def\\{\par}{\sl\theaddress}\par}
\expandafter\ifx\theemail\relax    
\relax\else\vglue 5pt plus 0.02fil minus 2pt\def\\{\stdspace{\rm 
and}\stdspace} 
\cl{Email:\stdspace\tt\theemail}\fi
\ifx\theurl\relax                  
\relax\else\vglue 5pt plus 0.02fil minus 2pt\def\\{\stdspace{\rm 
and}\stdspace}
\cl{URL:\stdspace\tt\theurl}\fi\par}

\vglue 7pt plus 0.3fil minus 3pt

{\bf Abstract}
\vglue 5pt plus 0.1fil minus 2pt

\theabstract

\vglue 7pt plus 0.3fil minus 3pt

{\bf AMS Classification numbers}\quad Primary:\quad \theprimaryclass

Secondary:\quad \thesecondaryclass

\vglue 5pt plus 0.3fil minus 2pt

{\bf Keywords}\quad \thekeywords

\vglue 10pt plus 0.5fil minus 5pt

{\small  Proposed: \theproposer\hfill Received: \receiveddate\nl
Seconded: \theseconders\hfill 
\ifx\reviseddate\relax                         
Accepted: \accepteddate                        
\else
Revised: \reviseddate                          
\fi}
\eject
}       

\let\maketitlepage\maketitlep
\let\maketitle\maketitlepage


\font\phead=cmsl9 scaled 950
\font\lhead=cmsl9 scaled 1050
\font\pnum=cmbx10 scaled 913
\font\lnum=cmbx10 
\font\pfoot=cmsl9 scaled 950
\font\lfoot=cmsl9 scaled 1050
\ifplaintex
\headline{\vbox to 0pt{\vskip -4.5mm\line{\small\phead\ifnum
\count0=\startpage ISSN 1364-0380 (on line)
1465-3060 (printed) \hfill {\pnum\folio}\else\ifodd\count0\def\\{ }%
\ifx\theshorttitle\relax\thetitle\else\theshorttitle\fi\hfill{\pnum\folio}
\else\def\\{ and }{\pnum\folio}\hfill\ifx\theshortauthors\relax\theauthors
\else\theshortauthors\fi\fi\fi}\vss}}
\footline{\vbox to 0pt{\vglue 0mm\line{\small\pfoot\ifnum\count0=\startpage
\copyright\ \gtp\hfill\else
\gt, Volume \thevolumenumber\ (\thevolumeyear)\hfill\fi}\vss
}}
\else
\makeatletter
\def\@oddhead{{\small\lhead\ifnum\count0=\startpage ISSN 1364-0380 (on line)
1465-3060 (printed) \hfill {\lnum\number\count0}\else\ifodd\count0
\def\\{ }\ifx\theshorttitle\relax \thetitle \else\theshorttitle\fi\hfill
{\lnum\number\count0}\else\def\\{ and }{\lnum\number\count0}
\hfill\ifx\theshortauthors\relax 
\theauthors\else\theshortauthors\fi\fi\fi}}\def\@evenhead{\@oddhead}
\def\@oddfoot{\small\lfoot\ifnum\count0=\startpage\copyright\ \gtp\hfill\else
\gt, Volume \thevolumenumber\ (\thevolumeyear)\hfill\fi}
\def\@evenfoot{\@oddfoot}
\makeatother
\fi


\newwrite\gtoutfile
\long\gdef\makeheadfile{  
{\def\\{, }\def\s{ }
\immediate\openout\gtoutfile head.xxx
\immediate\write\gtoutfile{To: math@arxiv.org}
\immediate\write\gtoutfile{Subject: put or rep NNNNN:pppp}
\immediate\write\gtoutfile{--text follows this line--}
\immediate\write\gtoutfile{Proxy-for: \ifx\theasciiauthors\relax
\theauthors\else\theasciiauthors\fi\s<\ifx\theasciiemail\relax\theemail\else\theasciiemail\fi>}
\immediate\write\gtoutfile{\noexpand\\}
\immediate\write\gtoutfile{Authors: \ifx\theasciiauthors\relax
\theauthors\else\theasciiauthors\fi}
\immediate\write\gtoutfile{Title: \ifx\theasciititle\relax
\thetitle\else\theasciititle\fi}
\immediate\write\gtoutfile{Subj-class: GT or SG or MG etc}
\immediate\write\gtoutfile{MSC-class: \theprimaryclass\ifx\thesecondaryclass\relax\else, \thesecondaryclass\fi}
\immediate\write\gtoutfile{Journal-ref: Geom. Topol. \thevolumenumber
(\thevolumeyear) \startpage-\finishpage}
\immediate\write\gtoutfile{Comments: Published by Geometry and Topology at}
\immediate\write\gtoutfile{\s\s http://www.maths.warwick.ac.uk/gt/GTVol\thevolumenumber/paper\thepapernumber.abs.html}
\immediate\write\gtoutfile{\noexpand\\}
\immediate\write\gtoutfile{}
\ifx\theasciiabstract\relax
\immediate\write\gtoutfile{\theabstract}\else
\immediate\write\gtoutfile{\theasciiabstract}\fi
\immediate\write\gtoutfile{}
\immediate\write\gtoutfile{\noexpand\\}
\immediate\write\gtoutfile{}
\immediate\closeout\gtoutfile}}  

\def\maketitlepage{\maketitlep\makeheadfile}
\let\maketitle\maketitlepage


\def\ifplaintex{\expandafter\ifx\csname documentclass\endcsname\relax}


\ifplaintex 
\hoffset 14truemm
\voffset 31truemm
\else
\headsep 23pt
\footskip 35pt
\hoffset -4truemm
\voffset 12.5truemm
\fi

\expandafter\ifx\csname beginpicture\endcsname\relax
\expandafter\ifx\csname documentclass\endcsname\relax
\input pictex \else\font\fiverm=cmr5
\input prepictex \input pictex \input postpictex \fi\fi

\def\gt{{\mathsurround=0pt\it $\cal G\mskip-2mu$eometry \&\ 
$\cal T\!\!$opology}}        

\def\gtp{{\mathsurround=0pt\it $\cal G\mskip-2mu$eometry \&\ 
$\cal T\!\!$opology $\cal P\!$ublications}}  


\def\lognumber#1{\def\thelognumber{#1}}
\def\volumenumber#1{\def\thevolumenumber{#1}}
\def\papernumber#1{\def\thepapernumber{#1}}
\def\volumeyear#1{\def\thevolumeyear{#1}}

\def\pagenumbers#1#2{\def\startpage{#1}\def\finishpage{#2}}
\def\published#1{\def\publishdate{#1}}
\def\proposed#1{\def\theproposer{#1}}
\def\seconded#1{\def\theseconders{#1}}
\def\received#1{\def\receiveddate{#1}}
\def\revised#1{\def\reviseddate{#1}}
\def\accepted#1{\def\accepteddate{#1}}

\long\def\asciiabstract#1{\long\def\theasciiabstract{#1}}


\let\\\par\let\thelognumber\relax
\let\thevolumenumber\relax\let\thepapernumber\relax
\let\thevolumeyear\relax\let\thesamplenumber\relax\let\startpage\relax
\let\finishpage\relax\let\publishdate\relax\let\receiveddate\relax
\let\reviseddate\relax\let\accepteddate\relax\let\theasciititle\relax
\let\theasciiauthors\relax
\let\theasciiabstract\relax
\let\theasciiemail\relax\let\theshortauthors\relax\let\theshorttitle\relax

\long\def\maketitlep{   

\count0=\startpage

\gt\hfill      
\beginpicture
\setcoordinatesystem units <0.33truein, 0.33truein> point at 2.2 0.9
\setplotsymbol ({$\cal G$})
\plotsymbolspacing=9truept
\circulararc 315 degrees from 0 1 center at 0 0
\setplotsymbol ({$\cal T$})
\circulararc 315 degrees from 1 -1 center at 1 0
\endpicture
%
\break
{\small\ifx\thesamplenumber\relax 
Volume \else Sample
\fi\thevolumenumber\ (\thevolumeyear)
\startpage--\finishpage\nl
Published: \publishdate}
\vglue 0.5truein plus 0.4fil minus 0.1truein

{\parskip=0pt\leftskip 0pt plus 1fil\def\\{\par\smallskip}{\ifplaintex\large
\else\Large\fi\bf\thetitle}\par\medskip}   

\vglue 0pt plus 0.1fil 

{\parskip=0pt\leftskip 0pt plus 1fil\def\\{\par}{\sc\theauthors}
\par\medskip}

\vglue 0pt plus 0.1fil 

{\small\parskip=0pt\let\newline\\
{\leftskip 0pt plus 1fil\def\\{\par}{\sl\theaddress}\par}
\expandafter\ifx\theemail\relax    
\relax\else\vglue 5pt plus 0.02fil minus 2pt\def\\{\stdspace{\rm 
and}\stdspace} 
\cl{Email:\stdspace\tt\theemail}\fi
\ifx\theurl\relax                  
\relax\else\vglue 5pt plus 0.02fil minus 2pt\def\\{\stdspace{\rm 
and}\stdspace}
\cl{URL:\stdspace\tt\theurl}\fi\par}

\vglue 7pt plus 0.3fil minus 3pt

{\bf Abstract}
\vglue 5pt plus 0.1fil minus 2pt

\theabstract

\vglue 7pt plus 0.3fil minus 3pt

{\bf AMS Classification numbers}\quad Primary:\quad \theprimaryclass

Secondary:\quad \thesecondaryclass

\vglue 5pt plus 0.3fil minus 2pt

{\bf Keywords}\quad \thekeywords

\vglue 10pt plus 0.5fil minus 5pt

{\small  Proposed: \theproposer\hfill Received: \receiveddate\nl
Seconded: \theseconders\hfill 
\ifx\reviseddate\relax                         
Accepted: \accepteddate                        
\else
Revised: \reviseddate                          
\fi}
\eject
}       

\let\maketitlepage\maketitlep
\let\maketitle\maketitlepage


\font\phead=cmsl9 scaled 950
\font\lhead=cmsl9 scaled 1050
\font\pnum=cmbx10 scaled 913
\font\lnum=cmbx10 
\font\pfoot=cmsl9 scaled 950
\font\lfoot=cmsl9 scaled 1050
\ifplaintex
\headline{\vbox to 0pt{\vskip -4.5mm\line{\small\phead\ifnum
\count0=\startpage ISSN 1364-0380 (on line)
1465-3060 (printed) \hfill {\pnum\folio}\else\ifodd\count0\def\\{ }%
\ifx\theshorttitle\relax\thetitle\else\theshorttitle\fi\hfill{\pnum\folio}
\else\def\\{ and }{\pnum\folio}\hfill\ifx\theshortauthors\relax\theauthors
\else\theshortauthors\fi\fi\fi}\vss}}
\footline{\vbox to 0pt{\vglue 0mm\line{\small\pfoot\ifnum\count0=\startpage
\copyright\ \gtp\hfill\else
\gt, Volume \thevolumenumber\ (\thevolumeyear)\hfill\fi}\vss
}}
\else
\makeatletter
\def\@oddhead{{\small\lhead\ifnum\count0=\startpage ISSN 1364-0380 (on line)
1465-3060 (printed) \hfill {\lnum\number\count0}\else\ifodd\count0
\def\\{ }\ifx\theshorttitle\relax \thetitle \else\theshorttitle\fi\hfill
{\lnum\number\count0}\else\def\\{ and }{\lnum\number\count0}
\hfill\ifx\theshortauthors\relax 
\theauthors\else\theshortauthors\fi\fi\fi}}\def\@evenhead{\@oddhead}
\def\@oddfoot{\small\lfoot\ifnum\count0=\startpage\copyright\ \gtp\hfill\else
\gt, Volume \thevolumenumber\ (\thevolumeyear)\hfill\fi}
\def\@evenfoot{\@oddfoot}
\makeatother
\fi


\newwrite\gtoutfile
\long\gdef\makeheadfile{  
{\def\\{, }\def\s{ }
\immediate\openout\gtoutfile head.xxx
\immediate\write\gtoutfile{To: math@arxiv.org}
\immediate\write\gtoutfile{Subject: put or rep NNNNN:pppp}
\immediate\write\gtoutfile{--text follows this line--}
\immediate\write\gtoutfile{Proxy-for: \ifx\theasciiauthors\relax
\theauthors\else\theasciiauthors\fi\s<\ifx\theasciiemail\relax\theemail\else\theasciiemail\fi>}
\immediate\write\gtoutfile{\noexpand\\}
\immediate\write\gtoutfile{Authors: \ifx\theasciiauthors\relax
\theauthors\else\theasciiauthors\fi}
\immediate\write\gtoutfile{Title: \ifx\theasciititle\relax
\thetitle\else\theasciititle\fi}
\immediate\write\gtoutfile{Subj-class: GT or SG or MG etc}
\immediate\write\gtoutfile{MSC-class: \theprimaryclass\ifx\thesecondaryclass\relax\else, \thesecondaryclass\fi}
\immediate\write\gtoutfile{Journal-ref: Geom. Topol. \thevolumenumber
(\thevolumeyear) \startpage-\finishpage}
\immediate\write\gtoutfile{Comments: Published by Geometry and Topology at}
\immediate\write\gtoutfile{\s\s http://www.maths.warwick.ac.uk/gt/GTVol\thevolumenumber/paper\thepapernumber.abs.html}
\immediate\write\gtoutfile{\noexpand\\}
\immediate\write\gtoutfile{}
\ifx\theasciiabstract\relax
\immediate\write\gtoutfile{\theabstract}\else
\immediate\write\gtoutfile{\theasciiabstract}\fi
\immediate\write\gtoutfile{}
\immediate\write\gtoutfile{\noexpand\\}
\immediate\write\gtoutfile{}
\immediate\closeout\gtoutfile}}  

\def\maketitlepage{\maketitlep\makeheadfile}
\let\maketitle\maketitlepage


\def\ifplaintex{\expandafter\ifx\csname documentclass\endcsname\relax}


\ifplaintex 
\hoffset 14truemm
\voffset 31truemm
\else
\headsep 23pt
\footskip 35pt
\hoffset -4truemm
\voffset 12.5truemm
\fi

\expandafter\ifx\csname beginpicture\endcsname\relax
\expandafter\ifx\csname documentclass\endcsname\relax
\input pictex \else\font\fiverm=cmr5
\input prepictex \input pictex \input postpictex \fi\fi

\def\gt{{\mathsurround=0pt\it $\cal G\mskip-2mu$eometry \&\ 
$\cal T\!\!$opology}}        

\def\gtp{{\mathsurround=0pt\it $\cal G\mskip-2mu$eometry \&\ 
$\cal T\!\!$opology $\cal P\!$ublications}}  


\def\lognumber#1{\def\thelognumber{#1}}
\def\volumenumber#1{\def\thevolumenumber{#1}}
\def\papernumber#1{\def\thepapernumber{#1}}
\def\volumeyear#1{\def\thevolumeyear{#1}}

\def\pagenumbers#1#2{\def\startpage{#1}\def\finishpage{#2}}
\def\published#1{\def\publishdate{#1}}
\def\proposed#1{\def\theproposer{#1}}
\def\seconded#1{\def\theseconders{#1}}
\def\received#1{\def\receiveddate{#1}}
\def\revised#1{\def\reviseddate{#1}}
\def\accepted#1{\def\accepteddate{#1}}

\long\def\asciiabstract#1{\long\def\theasciiabstract{#1}}


\let\\\par\let\thelognumber\relax
\let\thevolumenumber\relax\let\thepapernumber\relax
\let\thevolumeyear\relax\let\thesamplenumber\relax\let\startpage\relax
\let\finishpage\relax\let\publishdate\relax\let\receiveddate\relax
\let\reviseddate\relax\let\accepteddate\relax\let\theasciititle\relax
\let\theasciiauthors\relax
\let\theasciiabstract\relax
\let\theasciiemail\relax\let\theshortauthors\relax\let\theshorttitle\relax

\long\def\maketitlep{   

\count0=\startpage

\gt\hfill      
\beginpicture
\setcoordinatesystem units <0.33truein, 0.33truein> point at 2.2 0.9
\setplotsymbol ({$\cal G$})
\plotsymbolspacing=9truept
\circulararc 315 degrees from 0 1 center at 0 0
\setplotsymbol ({$\cal T$})
\circulararc 315 degrees from 1 -1 center at 1 0
\endpicture
%
\break
{\small\ifx\thesamplenumber\relax 
Volume \else Sample
\fi\thevolumenumber\ (\thevolumeyear)
\startpage--\finishpage\nl
Published: \publishdate}
\vglue 0.5truein plus 0.4fil minus 0.1truein

{\parskip=0pt\leftskip 0pt plus 1fil\def\\{\par\smallskip}{\ifplaintex\large
\else\Large\fi\bf\thetitle}\par\medskip}   

\vglue 0pt plus 0.1fil 

{\parskip=0pt\leftskip 0pt plus 1fil\def\\{\par}{\sc\theauthors}
\par\medskip}

\vglue 0pt plus 0.1fil 

{\small\parskip=0pt\let\newline\\
{\leftskip 0pt plus 1fil\def\\{\par}{\sl\theaddress}\par}
\expandafter\ifx\theemail\relax    
\relax\else\vglue 5pt plus 0.02fil minus 2pt\def\\{\stdspace{\rm 
and}\stdspace} 
\cl{Email:\stdspace\tt\theemail}\fi
\ifx\theurl\relax                  
\relax\else\vglue 5pt plus 0.02fil minus 2pt\def\\{\stdspace{\rm 
and}\stdspace}
\cl{URL:\stdspace\tt\theurl}\fi\par}

\vglue 7pt plus 0.3fil minus 3pt

{\bf Abstract}
\vglue 5pt plus 0.1fil minus 2pt

\theabstract

\vglue 7pt plus 0.3fil minus 3pt

{\bf AMS Classification numbers}\quad Primary:\quad \theprimaryclass

Secondary:\quad \thesecondaryclass

\vglue 5pt plus 0.3fil minus 2pt

{\bf Keywords}\quad \thekeywords

\vglue 10pt plus 0.5fil minus 5pt

{\small  Proposed: \theproposer\hfill Received: \receiveddate\nl
Seconded: \theseconders\hfill 
\ifx\reviseddate\relax                         
Accepted: \accepteddate                        
\else
Revised: \reviseddate                          
\fi}
\eject
}       

\let\maketitlepage\maketitlep
\let\maketitle\maketitlepage


\font\phead=cmsl9 scaled 950
\font\lhead=cmsl9 scaled 1050
\font\pnum=cmbx10 scaled 913
\font\lnum=cmbx10 
\font\pfoot=cmsl9 scaled 950
\font\lfoot=cmsl9 scaled 1050
\ifplaintex
\headline{\vbox to 0pt{\vskip -4.5mm\line{\small\phead\ifnum
\count0=\startpage ISSN 1364-0380 (on line)
1465-3060 (printed) \hfill {\pnum\folio}\else\ifodd\count0\def\\{ }%
\ifx\theshorttitle\relax\thetitle\else\theshorttitle\fi\hfill{\pnum\folio}
\else\def\\{ and }{\pnum\folio}\hfill\ifx\theshortauthors\relax\theauthors
\else\theshortauthors\fi\fi\fi}\vss}}
\footline{\vbox to 0pt{\vglue 0mm\line{\small\pfoot\ifnum\count0=\startpage
\copyright\ \gtp\hfill\else
\gt, Volume \thevolumenumber\ (\thevolumeyear)\hfill\fi}\vss
}}
\else
\makeatletter
\def\@oddhead{{\small\lhead\ifnum\count0=\startpage ISSN 1364-0380 (on line)
1465-3060 (printed) \hfill {\lnum\number\count0}\else\ifodd\count0
\def\\{ }\ifx\theshorttitle\relax \thetitle \else\theshorttitle\fi\hfill
{\lnum\number\count0}\else\def\\{ and }{\lnum\number\count0}
\hfill\ifx\theshortauthors\relax 
\theauthors\else\theshortauthors\fi\fi\fi}}\def\@evenhead{\@oddhead}
\def\@oddfoot{\small\lfoot\ifnum\count0=\startpage\copyright\ \gtp\hfill\else
\gt, Volume \thevolumenumber\ (\thevolumeyear)\hfill\fi}
\def\@evenfoot{\@oddfoot}
\makeatother
\fi


\newwrite\gtoutfile
\long\gdef\makeheadfile{  
{\def\\{, }\def\s{ }
\immediate\openout\gtoutfile head.xxx
\immediate\write\gtoutfile{To: math@arxiv.org}
\immediate\write\gtoutfile{Subject: put or rep NNNNN:pppp}
\immediate\write\gtoutfile{--text follows this line--}
\immediate\write\gtoutfile{Proxy-for: \ifx\theasciiauthors\relax
\theauthors\else\theasciiauthors\fi\s<\ifx\theasciiemail\relax\theemail\else\theasciiemail\fi>}
\immediate\write\gtoutfile{\noexpand\\}
\immediate\write\gtoutfile{Authors: \ifx\theasciiauthors\relax
\theauthors\else\theasciiauthors\fi}
\immediate\write\gtoutfile{Title: \ifx\theasciititle\relax
\thetitle\else\theasciititle\fi}
\immediate\write\gtoutfile{Subj-class: GT or SG or MG etc}
\immediate\write\gtoutfile{MSC-class: \theprimaryclass\ifx\thesecondaryclass\relax\else, \thesecondaryclass\fi}
\immediate\write\gtoutfile{Journal-ref: Geom. Topol. \thevolumenumber
(\thevolumeyear) \startpage-\finishpage}
\immediate\write\gtoutfile{Comments: Published by Geometry and Topology at}
\immediate\write\gtoutfile{\s\s http://www.maths.warwick.ac.uk/gt/GTVol\thevolumenumber/paper\thepapernumber.abs.html}
\immediate\write\gtoutfile{\noexpand\\}
\immediate\write\gtoutfile{}
\ifx\theasciiabstract\relax
\immediate\write\gtoutfile{\theabstract}\else
\immediate\write\gtoutfile{\theasciiabstract}\fi
\immediate\write\gtoutfile{}
\immediate\write\gtoutfile{\noexpand\\}
\immediate\write\gtoutfile{}
\immediate\closeout\gtoutfile}}  

\def\maketitlepage{\maketitlep\makeheadfile}
\let\maketitle\maketitlepage


\def\ifplaintex{\expandafter\ifx\csname documentclass\endcsname\relax}


\ifplaintex 
\hoffset 14truemm
\voffset 31truemm
\else
\headsep 23pt
\footskip 35pt
\hoffset -4truemm
\voffset 12.5truemm
\fi

\expandafter\ifx\csname beginpicture\endcsname\relax
\expandafter\ifx\csname documentclass\endcsname\relax
\input pictex \else\font\fiverm=cmr5
\input prepictex \input pictex \input postpictex \fi\fi

\def\gt{{\mathsurround=0pt\it $\cal G\mskip-2mu$eometry \&\ 
$\cal T\!\!$opology}}        

\def\gtp{{\mathsurround=0pt\it $\cal G\mskip-2mu$eometry \&\ 
$\cal T\!\!$opology $\cal P\!$ublications}}  


\def\lognumber#1{\def\thelognumber{#1}}
\def\volumenumber#1{\def\thevolumenumber{#1}}
\def\papernumber#1{\def\thepapernumber{#1}}
\def\volumeyear#1{\def\thevolumeyear{#1}}

\def\pagenumbers#1#2{\def\startpage{#1}\def\finishpage{#2}}
\def\published#1{\def\publishdate{#1}}
\def\proposed#1{\def\theproposer{#1}}
\def\seconded#1{\def\theseconders{#1}}
\def\received#1{\def\receiveddate{#1}}
\def\revised#1{\def\reviseddate{#1}}
\def\accepted#1{\def\accepteddate{#1}}

\long\def\asciiabstract#1{\long\def\theasciiabstract{#1}}


\let\\\par\let\thelognumber\relax
\let\thevolumenumber\relax\let\thepapernumber\relax
\let\thevolumeyear\relax\let\thesamplenumber\relax\let\startpage\relax
\let\finishpage\relax\let\publishdate\relax\let\receiveddate\relax
\let\reviseddate\relax\let\accepteddate\relax\let\theasciititle\relax
\let\theasciiauthors\relax
\let\theasciiabstract\relax
\let\theasciiemail\relax\let\theshortauthors\relax\let\theshorttitle\relax

\long\def\maketitlep{   

\count0=\startpage

\gt\hfill      
\beginpicture
\setcoordinatesystem units <0.33truein, 0.33truein> point at 2.2 0.9
\setplotsymbol ({$\cal G$})
\plotsymbolspacing=9truept
\circulararc 315 degrees from 0 1 center at 0 0
\setplotsymbol ({$\cal T$})
\circulararc 315 degrees from 1 -1 center at 1 0
\endpicture
%
\break
{\small\ifx\thesamplenumber\relax 
Volume \else Sample
\fi\thevolumenumber\ (\thevolumeyear)
\startpage--\finishpage\nl
Published: \publishdate}
\vglue 0.5truein plus 0.4fil minus 0.1truein

{\parskip=0pt\leftskip 0pt plus 1fil\def\\{\par\smallskip}{\ifplaintex\large
\else\Large\fi\bf\thetitle}\par\medskip}   

\vglue 0pt plus 0.1fil 

{\parskip=0pt\leftskip 0pt plus 1fil\def\\{\par}{\sc\theauthors}
\par\medskip}

\vglue 0pt plus 0.1fil 

{\small\parskip=0pt\let\newline\\
{\leftskip 0pt plus 1fil\def\\{\par}{\sl\theaddress}\par}
\expandafter\ifx\theemail\relax    
\relax\else\vglue 5pt plus 0.02fil minus 2pt\def\\{\stdspace{\rm 
and}\stdspace} 
\cl{Email:\stdspace\tt\theemail}\fi
\ifx\theurl\relax                  
\relax\else\vglue 5pt plus 0.02fil minus 2pt\def\\{\stdspace{\rm 
and}\stdspace}
\cl{URL:\stdspace\tt\theurl}\fi\par}

\vglue 7pt plus 0.3fil minus 3pt

{\bf Abstract}
\vglue 5pt plus 0.1fil minus 2pt

\theabstract

\vglue 7pt plus 0.3fil minus 3pt

{\bf AMS Classification numbers}\quad Primary:\quad \theprimaryclass

Secondary:\quad \thesecondaryclass

\vglue 5pt plus 0.3fil minus 2pt

{\bf Keywords}\quad \thekeywords

\vglue 10pt plus 0.5fil minus 5pt

{\small  Proposed: \theproposer\hfill Received: \receiveddate\nl
Seconded: \theseconders\hfill 
\ifx\reviseddate\relax                         
Accepted: \accepteddate                        
\else
Revised: \reviseddate                          
\fi}
\eject
}       

\let\maketitlepage\maketitlep
\let\maketitle\maketitlepage


\font\phead=cmsl9 scaled 950
\font\lhead=cmsl9 scaled 1050
\font\pnum=cmbx10 scaled 913
\font\lnum=cmbx10 
\font\pfoot=cmsl9 scaled 950
\font\lfoot=cmsl9 scaled 1050
\ifplaintex
\headline{\vbox to 0pt{\vskip -4.5mm\line{\small\phead\ifnum
\count0=\startpage ISSN 1364-0380 (on line)
1465-3060 (printed) \hfill {\pnum\folio}\else\ifodd\count0\def\\{ }%
\ifx\theshorttitle\relax\thetitle\else\theshorttitle\fi\hfill{\pnum\folio}
\else\def\\{ and }{\pnum\folio}\hfill\ifx\theshortauthors\relax\theauthors
\else\theshortauthors\fi\fi\fi}\vss}}
\footline{\vbox to 0pt{\vglue 0mm\line{\small\pfoot\ifnum\count0=\startpage
\copyright\ \gtp\hfill\else
\gt, Volume \thevolumenumber\ (\thevolumeyear)\hfill\fi}\vss
}}
\else
\makeatletter
\def\@oddhead{{\small\lhead\ifnum\count0=\startpage ISSN 1364-0380 (on line)
1465-3060 (printed) \hfill {\lnum\number\count0}\else\ifodd\count0
\def\\{ }\ifx\theshorttitle\relax \thetitle \else\theshorttitle\fi\hfill
{\lnum\number\count0}\else\def\\{ and }{\lnum\number\count0}
\hfill\ifx\theshortauthors\relax 
\theauthors\else\theshortauthors\fi\fi\fi}}\def\@evenhead{\@oddhead}
\def\@oddfoot{\small\lfoot\ifnum\count0=\startpage\copyright\ \gtp\hfill\else
\gt, Volume \thevolumenumber\ (\thevolumeyear)\hfill\fi}
\def\@evenfoot{\@oddfoot}
\makeatother
\fi


\newwrite\gtoutfile
\long\gdef\makeheadfile{  
{\def\\{, }\def\s{ }
\immediate\openout\gtoutfile head.xxx
\immediate\write\gtoutfile{To: math@arxiv.org}
\immediate\write\gtoutfile{Subject: put or rep NNNNN:pppp}
\immediate\write\gtoutfile{--text follows this line--}
\immediate\write\gtoutfile{Proxy-for: \ifx\theasciiauthors\relax
\theauthors\else\theasciiauthors\fi\s<\ifx\theasciiemail\relax\theemail\else\theasciiemail\fi>}
\immediate\write\gtoutfile{\noexpand\\}
\immediate\write\gtoutfile{Authors: \ifx\theasciiauthors\relax
\theauthors\else\theasciiauthors\fi}
\immediate\write\gtoutfile{Title: \ifx\theasciititle\relax
\thetitle\else\theasciititle\fi}
\immediate\write\gtoutfile{Subj-class: GT or SG or MG etc}
\immediate\write\gtoutfile{MSC-class: \theprimaryclass\ifx\thesecondaryclass\relax\else, \thesecondaryclass\fi}
\immediate\write\gtoutfile{Journal-ref: Geom. Topol. \thevolumenumber
(\thevolumeyear) \startpage-\finishpage}
\immediate\write\gtoutfile{Comments: Published by Geometry and Topology at}
\immediate\write\gtoutfile{\s\s http://www.maths.warwick.ac.uk/gt/GTVol\thevolumenumber/paper\thepapernumber.abs.html}
\immediate\write\gtoutfile{\noexpand\\}
\immediate\write\gtoutfile{}
\ifx\theasciiabstract\relax
\immediate\write\gtoutfile{\theabstract}\else
\immediate\write\gtoutfile{\theasciiabstract}\fi
\immediate\write\gtoutfile{}
\immediate\write\gtoutfile{\noexpand\\}
\immediate\write\gtoutfile{}
\immediate\closeout\gtoutfile}}  

\def\maketitlepage{\maketitlep\makeheadfile}
\let\maketitle\maketitlepage


\def\ifplaintex{\expandafter\ifx\csname documentclass\endcsname\relax}


\ifplaintex 
\hoffset 14truemm
\voffset 31truemm
\else
\headsep 23pt
\footskip 35pt
\hoffset -4truemm
\voffset 12.5truemm
\fi

\expandafter\ifx\csname beginpicture\endcsname\relax
\expandafter\ifx\csname documentclass\endcsname\relax
\input pictex \else\font\fiverm=cmr5
\input prepictex \input pictex \input postpictex \fi\fi

\def\gt{{\mathsurround=0pt\it $\cal G\mskip-2mu$eometry \&\ 
$\cal T\!\!$opology}}        

\def\gtp{{\mathsurround=0pt\it $\cal G\mskip-2mu$eometry \&\ 
$\cal T\!\!$opology $\cal P\!$ublications}}  


\def\lognumber#1{\def\thelognumber{#1}}
\def\volumenumber#1{\def\thevolumenumber{#1}}
\def\papernumber#1{\def\thepapernumber{#1}}
\def\volumeyear#1{\def\thevolumeyear{#1}}

\def\pagenumbers#1#2{\def\startpage{#1}\def\finishpage{#2}}
\def\published#1{\def\publishdate{#1}}
\def\proposed#1{\def\theproposer{#1}}
\def\seconded#1{\def\theseconders{#1}}
\def\received#1{\def\receiveddate{#1}}
\def\revised#1{\def\reviseddate{#1}}
\def\accepted#1{\def\accepteddate{#1}}

\long\def\asciiabstract#1{\long\def\theasciiabstract{#1}}


\let\\\par\let\thelognumber\relax
\let\thevolumenumber\relax\let\thepapernumber\relax
\let\thevolumeyear\relax\let\thesamplenumber\relax\let\startpage\relax
\let\finishpage\relax\let\publishdate\relax\let\receiveddate\relax
\let\reviseddate\relax\let\accepteddate\relax\let\theasciititle\relax
\let\theasciiauthors\relax
\let\theasciiabstract\relax
\let\theasciiemail\relax\let\theshortauthors\relax\let\theshorttitle\relax

\long\def\maketitlep{   

\count0=\startpage

\gt\hfill      
\beginpicture
\setcoordinatesystem units <0.33truein, 0.33truein> point at 2.2 0.9
\setplotsymbol ({$\cal G$})
\plotsymbolspacing=9truept
\circulararc 315 degrees from 0 1 center at 0 0
\setplotsymbol ({$\cal T$})
\circulararc 315 degrees from 1 -1 center at 1 0
\endpicture
%
\break
{\small\ifx\thesamplenumber\relax 
Volume \else Sample
\fi\thevolumenumber\ (\thevolumeyear)
\startpage--\finishpage\nl
Published: \publishdate}
\vglue 0.5truein plus 0.4fil minus 0.1truein

{\parskip=0pt\leftskip 0pt plus 1fil\def\\{\par\smallskip}{\ifplaintex\large
\else\Large\fi\bf\thetitle}\par\medskip}   

\vglue 0pt plus 0.1fil 

{\parskip=0pt\leftskip 0pt plus 1fil\def\\{\par}{\sc\theauthors}
\par\medskip}

\vglue 0pt plus 0.1fil 

{\small\parskip=0pt\let\newline\\
{\leftskip 0pt plus 1fil\def\\{\par}{\sl\theaddress}\par}
\expandafter\ifx\theemail\relax    
\relax\else\vglue 5pt plus 0.02fil minus 2pt\def\\{\stdspace{\rm 
and}\stdspace} 
\cl{Email:\stdspace\tt\theemail}\fi
\ifx\theurl\relax                  
\relax\else\vglue 5pt plus 0.02fil minus 2pt\def\\{\stdspace{\rm 
and}\stdspace}
\cl{URL:\stdspace\tt\theurl}\fi\par}

\vglue 7pt plus 0.3fil minus 3pt

{\bf Abstract}
\vglue 5pt plus 0.1fil minus 2pt

\theabstract

\vglue 7pt plus 0.3fil minus 3pt

{\bf AMS Classification numbers}\quad Primary:\quad \theprimaryclass

Secondary:\quad \thesecondaryclass

\vglue 5pt plus 0.3fil minus 2pt

{\bf Keywords}\quad \thekeywords

\vglue 10pt plus 0.5fil minus 5pt

{\small  Proposed: \theproposer\hfill Received: \receiveddate\nl
Seconded: \theseconders\hfill 
\ifx\reviseddate\relax                         
Accepted: \accepteddate                        
\else
Revised: \reviseddate                          
\fi}
\eject
}       

\let\maketitlepage\maketitlep
\let\maketitle\maketitlepage


\font\phead=cmsl9 scaled 950
\font\lhead=cmsl9 scaled 1050
\font\pnum=cmbx10 scaled 913
\font\lnum=cmbx10 
\font\pfoot=cmsl9 scaled 950
\font\lfoot=cmsl9 scaled 1050
\ifplaintex
\headline{\vbox to 0pt{\vskip -4.5mm\line{\small\phead\ifnum
\count0=\startpage ISSN 1364-0380 (on line)
1465-3060 (printed) \hfill {\pnum\folio}\else\ifodd\count0\def\\{ }%
\ifx\theshorttitle\relax\thetitle\else\theshorttitle\fi\hfill{\pnum\folio}
\else\def\\{ and }{\pnum\folio}\hfill\ifx\theshortauthors\relax\theauthors
\else\theshortauthors\fi\fi\fi}\vss}}
\footline{\vbox to 0pt{\vglue 0mm\line{\small\pfoot\ifnum\count0=\startpage
\copyright\ \gtp\hfill\else
\gt, Volume \thevolumenumber\ (\thevolumeyear)\hfill\fi}\vss
}}
\else
\makeatletter
\def\@oddhead{{\small\lhead\ifnum\count0=\startpage ISSN 1364-0380 (on line)
1465-3060 (printed) \hfill {\lnum\number\count0}\else\ifodd\count0
\def\\{ }\ifx\theshorttitle\relax \thetitle \else\theshorttitle\fi\hfill
{\lnum\number\count0}\else\def\\{ and }{\lnum\number\count0}
\hfill\ifx\theshortauthors\relax 
\theauthors\else\theshortauthors\fi\fi\fi}}\def\@evenhead{\@oddhead}
\def\@oddfoot{\small\lfoot\ifnum\count0=\startpage\copyright\ \gtp\hfill\else
\gt, Volume \thevolumenumber\ (\thevolumeyear)\hfill\fi}
\def\@evenfoot{\@oddfoot}
\makeatother
\fi


\newwrite\gtoutfile
\long\gdef\makeheadfile{  
{\def\\{, }\def\s{ }
\immediate\openout\gtoutfile head.xxx
\immediate\write\gtoutfile{To: math@arxiv.org}
\immediate\write\gtoutfile{Subject: put or rep NNNNN:pppp}
\immediate\write\gtoutfile{--text follows this line--}
\immediate\write\gtoutfile{Proxy-for: \ifx\theasciiauthors\relax
\theauthors\else\theasciiauthors\fi\s<\ifx\theasciiemail\relax\theemail\else\theasciiemail\fi>}
\immediate\write\gtoutfile{\noexpand\\}
\immediate\write\gtoutfile{Authors: \ifx\theasciiauthors\relax
\theauthors\else\theasciiauthors\fi}
\immediate\write\gtoutfile{Title: \ifx\theasciititle\relax
\thetitle\else\theasciititle\fi}
\immediate\write\gtoutfile{Subj-class: GT or SG or MG etc}
\immediate\write\gtoutfile{MSC-class: \theprimaryclass\ifx\thesecondaryclass\relax\else, \thesecondaryclass\fi}
\immediate\write\gtoutfile{Journal-ref: Geom. Topol. \thevolumenumber
(\thevolumeyear) \startpage-\finishpage}
\immediate\write\gtoutfile{Comments: Published by Geometry and Topology at}
\immediate\write\gtoutfile{\s\s http://www.maths.warwick.ac.uk/gt/GTVol\thevolumenumber/paper\thepapernumber.abs.html}
\immediate\write\gtoutfile{\noexpand\\}
\immediate\write\gtoutfile{}
\ifx\theasciiabstract\relax
\immediate\write\gtoutfile{\theabstract}\else
\immediate\write\gtoutfile{\theasciiabstract}\fi
\immediate\write\gtoutfile{}
\immediate\write\gtoutfile{\noexpand\\}
\immediate\write\gtoutfile{}
\immediate\closeout\gtoutfile}}  

\def\maketitlepage{\maketitlep\makeheadfile}
\let\maketitle\maketitlepage


\def\ifplaintex{\expandafter\ifx\csname documentclass\endcsname\relax}


\ifplaintex 
\hoffset 14truemm
\voffset 31truemm
\else
\headsep 23pt
\footskip 35pt
\hoffset -4truemm
\voffset 12.5truemm
\fi

\expandafter\ifx\csname beginpicture\endcsname\relax
\expandafter\ifx\csname documentclass\endcsname\relax
\input pictex \else\font\fiverm=cmr5
\input prepictex \input pictex \input postpictex \fi\fi

\def\gt{{\mathsurround=0pt\it $\cal G\mskip-2mu$eometry \&\ 
$\cal T\!\!$opology}}        

\def\gtp{{\mathsurround=0pt\it $\cal G\mskip-2mu$eometry \&\ 
$\cal T\!\!$opology $\cal P\!$ublications}}  


\def\lognumber#1{\def\thelognumber{#1}}
\def\volumenumber#1{\def\thevolumenumber{#1}}
\def\papernumber#1{\def\thepapernumber{#1}}
\def\volumeyear#1{\def\thevolumeyear{#1}}

\def\pagenumbers#1#2{\def\startpage{#1}\def\finishpage{#2}}
\def\published#1{\def\publishdate{#1}}
\def\proposed#1{\def\theproposer{#1}}
\def\seconded#1{\def\theseconders{#1}}
\def\received#1{\def\receiveddate{#1}}
\def\revised#1{\def\reviseddate{#1}}
\def\accepted#1{\def\accepteddate{#1}}

\long\def\asciiabstract#1{\long\def\theasciiabstract{#1}}


\let\\\par\let\thelognumber\relax
\let\thevolumenumber\relax\let\thepapernumber\relax
\let\thevolumeyear\relax\let\thesamplenumber\relax\let\startpage\relax
\let\finishpage\relax\let\publishdate\relax\let\receiveddate\relax
\let\reviseddate\relax\let\accepteddate\relax\let\theasciititle\relax
\let\theasciiauthors\relax
\let\theasciiabstract\relax
\let\theasciiemail\relax\let\theshortauthors\relax\let\theshorttitle\relax

\long\def\maketitlep{   

\count0=\startpage

\gt\hfill      
\beginpicture
\setcoordinatesystem units <0.33truein, 0.33truein> point at 2.2 0.9
\setplotsymbol ({$\cal G$})
\plotsymbolspacing=9truept
\circulararc 315 degrees from 0 1 center at 0 0
\setplotsymbol ({$\cal T$})
\circulararc 315 degrees from 1 -1 center at 1 0
\endpicture
%
\break
{\small\ifx\thesamplenumber\relax 
Volume \else Sample
\fi\thevolumenumber\ (\thevolumeyear)
\startpage--\finishpage\nl
Published: \publishdate}
\vglue 0.5truein plus 0.4fil minus 0.1truein

{\parskip=0pt\leftskip 0pt plus 1fil\def\\{\par\smallskip}{\ifplaintex\large
\else\Large\fi\bf\thetitle}\par\medskip}   

\vglue 0pt plus 0.1fil 

{\parskip=0pt\leftskip 0pt plus 1fil\def\\{\par}{\sc\theauthors}
\par\medskip}

\vglue 0pt plus 0.1fil 

{\small\parskip=0pt\let\newline\\
{\leftskip 0pt plus 1fil\def\\{\par}{\sl\theaddress}\par}
\expandafter\ifx\theemail\relax    
\relax\else\vglue 5pt plus 0.02fil minus 2pt\def\\{\stdspace{\rm 
and}\stdspace} 
\cl{Email:\stdspace\tt\theemail}\fi
\ifx\theurl\relax                  
\relax\else\vglue 5pt plus 0.02fil minus 2pt\def\\{\stdspace{\rm 
and}\stdspace}
\cl{URL:\stdspace\tt\theurl}\fi\par}

\vglue 7pt plus 0.3fil minus 3pt

{\bf Abstract}
\vglue 5pt plus 0.1fil minus 2pt

\theabstract

\vglue 7pt plus 0.3fil minus 3pt

{\bf AMS Classification numbers}\quad Primary:\quad \theprimaryclass

Secondary:\quad \thesecondaryclass

\vglue 5pt plus 0.3fil minus 2pt

{\bf Keywords}\quad \thekeywords

\vglue 10pt plus 0.5fil minus 5pt

{\small  Proposed: \theproposer\hfill Received: \receiveddate\nl
Seconded: \theseconders\hfill 
\ifx\reviseddate\relax                         
Accepted: \accepteddate                        
\else
Revised: \reviseddate                          
\fi}
\eject
}       

\let\maketitlepage\maketitlep
\let\maketitle\maketitlepage


\font\phead=cmsl9 scaled 950
\font\lhead=cmsl9 scaled 1050
\font\pnum=cmbx10 scaled 913
\font\lnum=cmbx10 
\font\pfoot=cmsl9 scaled 950
\font\lfoot=cmsl9 scaled 1050
\ifplaintex
\headline{\vbox to 0pt{\vskip -4.5mm\line{\small\phead\ifnum
\count0=\startpage ISSN 1364-0380 (on line)
1465-3060 (printed) \hfill {\pnum\folio}\else\ifodd\count0\def\\{ }%
\ifx\theshorttitle\relax\thetitle\else\theshorttitle\fi\hfill{\pnum\folio}
\else\def\\{ and }{\pnum\folio}\hfill\ifx\theshortauthors\relax\theauthors
\else\theshortauthors\fi\fi\fi}\vss}}
\footline{\vbox to 0pt{\vglue 0mm\line{\small\pfoot\ifnum\count0=\startpage
\copyright\ \gtp\hfill\else
\gt, Volume \thevolumenumber\ (\thevolumeyear)\hfill\fi}\vss
}}
\else
\makeatletter
\def\@oddhead{{\small\lhead\ifnum\count0=\startpage ISSN 1364-0380 (on line)
1465-3060 (printed) \hfill {\lnum\number\count0}\else\ifodd\count0
\def\\{ }\ifx\theshorttitle\relax \thetitle \else\theshorttitle\fi\hfill
{\lnum\number\count0}\else\def\\{ and }{\lnum\number\count0}
\hfill\ifx\theshortauthors\relax 
\theauthors\else\theshortauthors\fi\fi\fi}}\def\@evenhead{\@oddhead}
\def\@oddfoot{\small\lfoot\ifnum\count0=\startpage\copyright\ \gtp\hfill\else
\gt, Volume \thevolumenumber\ (\thevolumeyear)\hfill\fi}
\def\@evenfoot{\@oddfoot}
\makeatother
\fi


\newwrite\gtoutfile
\long\gdef\makeheadfile{  
{\def\\{, }\def\s{ }
\immediate\openout\gtoutfile head.xxx
\immediate\write\gtoutfile{To: math@arxiv.org}
\immediate\write\gtoutfile{Subject: put or rep NNNNN:pppp}
\immediate\write\gtoutfile{--text follows this line--}
\immediate\write\gtoutfile{Proxy-for: \ifx\theasciiauthors\relax
\theauthors\else\theasciiauthors\fi\s<\ifx\theasciiemail\relax\theemail\else\theasciiemail\fi>}
\immediate\write\gtoutfile{\noexpand\\}
\immediate\write\gtoutfile{Authors: \ifx\theasciiauthors\relax
\theauthors\else\theasciiauthors\fi}
\immediate\write\gtoutfile{Title: \ifx\theasciititle\relax
\thetitle\else\theasciititle\fi}
\immediate\write\gtoutfile{Subj-class: GT or SG or MG etc}
\immediate\write\gtoutfile{MSC-class: \theprimaryclass\ifx\thesecondaryclass\relax\else, \thesecondaryclass\fi}
\immediate\write\gtoutfile{Journal-ref: Geom. Topol. \thevolumenumber
(\thevolumeyear) \startpage-\finishpage}
\immediate\write\gtoutfile{Comments: Published by Geometry and Topology at}
\immediate\write\gtoutfile{\s\s http://www.maths.warwick.ac.uk/gt/GTVol\thevolumenumber/paper\thepapernumber.abs.html}
\immediate\write\gtoutfile{\noexpand\\}
\immediate\write\gtoutfile{}
\ifx\theasciiabstract\relax
\immediate\write\gtoutfile{\theabstract}\else
\immediate\write\gtoutfile{\theasciiabstract}\fi
\immediate\write\gtoutfile{}
\immediate\write\gtoutfile{\noexpand\\}
\immediate\write\gtoutfile{}
\immediate\closeout\gtoutfile}}  

\def\maketitlepage{\maketitlep\makeheadfile}
\let\maketitle\maketitlepage


\def\ifplaintex{\expandafter\ifx\csname documentclass\endcsname\relax}


\ifplaintex 
\hoffset 14truemm
\voffset 31truemm
\else
\headsep 23pt
\footskip 35pt
\hoffset -4truemm
\voffset 12.5truemm
\fi

\expandafter\ifx\csname beginpicture\endcsname\relax
\expandafter\ifx\csname documentclass\endcsname\relax
\input pictex \else\font\fiverm=cmr5
\input prepictex \input pictex \input postpictex \fi\fi

\def\gt{{\mathsurround=0pt\it $\cal G\mskip-2mu$eometry \&\ 
$\cal T\!\!$opology}}        

\def\gtp{{\mathsurround=0pt\it $\cal G\mskip-2mu$eometry \&\ 
$\cal T\!\!$opology $\cal P\!$ublications}}  


\def\lognumber#1{\def\thelognumber{#1}}
\def\volumenumber#1{\def\thevolumenumber{#1}}
\def\papernumber#1{\def\thepapernumber{#1}}
\def\volumeyear#1{\def\thevolumeyear{#1}}

\def\pagenumbers#1#2{\def\startpage{#1}\def\finishpage{#2}}
\def\published#1{\def\publishdate{#1}}
\def\proposed#1{\def\theproposer{#1}}
\def\seconded#1{\def\theseconders{#1}}
\def\received#1{\def\receiveddate{#1}}
\def\revised#1{\def\reviseddate{#1}}
\def\accepted#1{\def\accepteddate{#1}}

\long\def\asciiabstract#1{\long\def\theasciiabstract{#1}}


\let\\\par\let\thelognumber\relax
\let\thevolumenumber\relax\let\thepapernumber\relax
\let\thevolumeyear\relax\let\thesamplenumber\relax\let\startpage\relax
\let\finishpage\relax\let\publishdate\relax\let\receiveddate\relax
\let\reviseddate\relax\let\accepteddate\relax\let\theasciititle\relax
\let\theasciiauthors\relax
\let\theasciiabstract\relax
\let\theasciiemail\relax\let\theshortauthors\relax\let\theshorttitle\relax

\long\def\maketitlep{   

\count0=\startpage

\gt\hfill      
\beginpicture
\setcoordinatesystem units <0.33truein, 0.33truein> point at 2.2 0.9
\setplotsymbol ({$\cal G$})
\plotsymbolspacing=9truept
\circulararc 315 degrees from 0 1 center at 0 0
\setplotsymbol ({$\cal T$})
\circulararc 315 degrees from 1 -1 center at 1 0
\endpicture
%
\break
{\small\ifx\thesamplenumber\relax 
Volume \else Sample
\fi\thevolumenumber\ (\thevolumeyear)
\startpage--\finishpage\nl
Published: \publishdate}
\vglue 0.5truein plus 0.4fil minus 0.1truein

{\parskip=0pt\leftskip 0pt plus 1fil\def\\{\par\smallskip}{\ifplaintex\large
\else\Large\fi\bf\thetitle}\par\medskip}   

\vglue 0pt plus 0.1fil 

{\parskip=0pt\leftskip 0pt plus 1fil\def\\{\par}{\sc\theauthors}
\par\medskip}

\vglue 0pt plus 0.1fil 

{\small\parskip=0pt\let\newline\\
{\leftskip 0pt plus 1fil\def\\{\par}{\sl\theaddress}\par}
\expandafter\ifx\theemail\relax    
\relax\else\vglue 5pt plus 0.02fil minus 2pt\def\\{\stdspace{\rm 
and}\stdspace} 
\cl{Email:\stdspace\tt\theemail}\fi
\ifx\theurl\relax                  
\relax\else\vglue 5pt plus 0.02fil minus 2pt\def\\{\stdspace{\rm 
and}\stdspace}
\cl{URL:\stdspace\tt\theurl}\fi\par}

\vglue 7pt plus 0.3fil minus 3pt

{\bf Abstract}
\vglue 5pt plus 0.1fil minus 2pt

\theabstract

\vglue 7pt plus 0.3fil minus 3pt

{\bf AMS Classification numbers}\quad Primary:\quad \theprimaryclass

Secondary:\quad \thesecondaryclass

\vglue 5pt plus 0.3fil minus 2pt

{\bf Keywords}\quad \thekeywords

\vglue 10pt plus 0.5fil minus 5pt

{\small  Proposed: \theproposer\hfill Received: \receiveddate\nl
Seconded: \theseconders\hfill 
\ifx\reviseddate\relax                         
Accepted: \accepteddate                        
\else
Revised: \reviseddate                          
\fi}
\eject
}       

\let\maketitlepage\maketitlep
\let\maketitle\maketitlepage


\font\phead=cmsl9 scaled 950
\font\lhead=cmsl9 scaled 1050
\font\pnum=cmbx10 scaled 913
\font\lnum=cmbx10 
\font\pfoot=cmsl9 scaled 950
\font\lfoot=cmsl9 scaled 1050
\ifplaintex
\headline{\vbox to 0pt{\vskip -4.5mm\line{\small\phead\ifnum
\count0=\startpage ISSN 1364-0380 (on line)
1465-3060 (printed) \hfill {\pnum\folio}\else\ifodd\count0\def\\{ }%
\ifx\theshorttitle\relax\thetitle\else\theshorttitle\fi\hfill{\pnum\folio}
\else\def\\{ and }{\pnum\folio}\hfill\ifx\theshortauthors\relax\theauthors
\else\theshortauthors\fi\fi\fi}\vss}}
\footline{\vbox to 0pt{\vglue 0mm\line{\small\pfoot\ifnum\count0=\startpage
\copyright\ \gtp\hfill\else
\gt, Volume \thevolumenumber\ (\thevolumeyear)\hfill\fi}\vss
}}
\else
\makeatletter
\def\@oddhead{{\small\lhead\ifnum\count0=\startpage ISSN 1364-0380 (on line)
1465-3060 (printed) \hfill {\lnum\number\count0}\else\ifodd\count0
\def\\{ }\ifx\theshorttitle\relax \thetitle \else\theshorttitle\fi\hfill
{\lnum\number\count0}\else\def\\{ and }{\lnum\number\count0}
\hfill\ifx\theshortauthors\relax 
\theauthors\else\theshortauthors\fi\fi\fi}}\def\@evenhead{\@oddhead}
\def\@oddfoot{\small\lfoot\ifnum\count0=\startpage\copyright\ \gtp\hfill\else
\gt, Volume \thevolumenumber\ (\thevolumeyear)\hfill\fi}
\def\@evenfoot{\@oddfoot}
\makeatother
\fi


\newwrite\gtoutfile
\long\gdef\makeheadfile{  
{\def\\{, }\def\s{ }
\immediate\openout\gtoutfile head.xxx
\immediate\write\gtoutfile{To: math@arxiv.org}
\immediate\write\gtoutfile{Subject: put or rep NNNNN:pppp}
\immediate\write\gtoutfile{--text follows this line--}
\immediate\write\gtoutfile{Proxy-for: \ifx\theasciiauthors\relax
\theauthors\else\theasciiauthors\fi\s<\ifx\theasciiemail\relax\theemail\else\theasciiemail\fi>}
\immediate\write\gtoutfile{\noexpand\\}
\immediate\write\gtoutfile{Authors: \ifx\theasciiauthors\relax
\theauthors\else\theasciiauthors\fi}
\immediate\write\gtoutfile{Title: \ifx\theasciititle\relax
\thetitle\else\theasciititle\fi}
\immediate\write\gtoutfile{Subj-class: GT or SG or MG etc}
\immediate\write\gtoutfile{MSC-class: \theprimaryclass\ifx\thesecondaryclass\relax\else, \thesecondaryclass\fi}
\immediate\write\gtoutfile{Journal-ref: Geom. Topol. \thevolumenumber
(\thevolumeyear) \startpage-\finishpage}
\immediate\write\gtoutfile{Comments: Published by Geometry and Topology at}
\immediate\write\gtoutfile{\s\s http://www.maths.warwick.ac.uk/gt/GTVol\thevolumenumber/paper\thepapernumber.abs.html}
\immediate\write\gtoutfile{\noexpand\\}
\immediate\write\gtoutfile{}
\ifx\theasciiabstract\relax
\immediate\write\gtoutfile{\theabstract}\else
\immediate\write\gtoutfile{\theasciiabstract}\fi
\immediate\write\gtoutfile{}
\immediate\write\gtoutfile{\noexpand\\}
\immediate\write\gtoutfile{}
\immediate\closeout\gtoutfile}}  

\def\maketitlepage{\maketitlep\makeheadfile}
\let\maketitle\maketitlepage

\lognumber{188}
\volumenumber{5}\papernumber{23}\volumeyear{2001}
\pagenumbers{719}{760}
\received{15 June 2001}
\revised{6 September 2001}
\accepted{5 October 2001}
\proposed{Yasha Elaishberg}
\seconded{Joan Birman, Robion Kirby}
\published{11 October 2001}

\input amsnames
\input amstex
\input epsf     
\chardef\newinsCatAt\the\catcode `\@
\catcode `\@=11
%
%
%
\newskip\insertskipamount\newskip\inserthardskipamount
\insertskipamount 12pt plus2pt     
\inserthardskipamount 4pt          
\def\insertskip{\vskip\insertskipamount}
%
%
\newskip\LastSkip
\def\SaveLastSkip{\LastSkip\lastskip}
\def\RestoreLastSkip{\nobreak\vskip-\LastSkip\vskip\LastSkip}
%
%
\newcount\SplitTest
\def\SetSplitTest{\SplitTest\insertpenalties
  \insert\topins{\floatingpenalty1}%
  \advance\SplitTest-\insertpenalties}
%
%
\def\midinsert{\par
 \SaveLastSkip\penalty-150\SetSplitTest\RestoreLastSkip
 \ifnum\SplitTest=-1
  \@midfalse\p@gefalse\else\@midtrue\fi\@ins}
\def\@ins{\par\begingroup\setbox\z@\vbox\bgroup%
  \vglue\inserthardskipamount}
\def\endinsert{\egroup 
  \if@mid \dimen@\ht\z@ \advance\dimen@\dp\z@
    \advance\dimen@\insertskipamount
    \advance\dimen@\pagetotal\advance\dimen@-\pageshrink
    \ifdim\dimen@>\pagegoal\@midfalse\p@gefalse\fi\fi
  \if@mid%
    \ifdim\lastskip<\insertskipamount\removelastskip\insertskip\fi
    \nointerlineskip\box\z@\penalty-200\insertskip
  \else%
    \SaveLastSkip
    \insert\topins{\penalty100 
    \splittopskip\z@skip
    \splitmaxdepth\maxdimen \floatingpenalty\z@
    \ifp@ge \dimen@\dp\z@
    \vbox to\vsize{\unvbox\z@\kern-\dimen@}
    \else \box\z@\nobreak\insertskip\fi}
    \RestoreLastSkip
   \fi\endgroup}
%
\catcode `\@=\newinsCatAt

\let\cal\Cal    
\catcode`\@=12  
\let\\\par
\def\topmatter{\relax}
\def\endtopmatter{\maketitlepage}
\let\gttitle\title
\def\title#1\endtitle{\gttitle{#1}}
\let\gtauthor\author
\def\author#1\endauthor{\gtauthor{#1}}
\let\gtaddress\address
\def\address#1\endaddress{\gtaddress{#1}}
\def\affil#1\endaffil{\gtaddress{#1}}
\let\gtemail\email
\def\email#1\endemail{\gtemail{#1}}
\def\subjclass#1\endsubjclass{\primaryclass{#1}}
\let\gtkeywords\keywords
\def\keywords#1\endkeywords{\gtkeywords{#1}}
\def\heading#1\endheading{{\def\S##1{\relax}\def\\{\relax\ignorespaces}
    \section{#1}}}
\def\head#1\endhead{\heading#1\endheading}
\def\subheading#1{\sh{#1}}
\def\subhead#1\endsubhead{\sh{#1}}
\def\subsubhead#1\endsubsubhead{\sh{#1}}
\def\specialhead#1\endspecialhead{\sh{#1}}
\def\demo#1{\rk{#1}\ignorespaces}
\def\enddemo{\ppar}
\let\remark\demo
\def\endremark{}
\let\definition\demo
\def\enddefinition{\ppar}
\let\example\demo
\def\endexample{\ppar}
\def\qed{\ifmmode\quad\sq\else\hbox{}\hfill$\sq$\par\goodbreak\rm\fi}  
\def\proclaim#1{\rk{#1}\sl\ignorespaces}
\def\endproclaim{\rm\ppar}
\def\cite#1{[#1]}
\newcount\itemnumber
\def\roster{\items\itemnumber=1}
\def\endroster{\enditems}
\let\itemold\item
\def\item{\itemold{{\rm(\number\itemnumber)}}%
\global\advance\itemnumber by 1\ignorespaces}
\def\S{section~\ignorespaces}  
\def\date#1\enddate{\relax}
\def\thanks#1\endthanks{\relax}   
\def\dedicatory#1\enddedicatory{\relax}  
\let\footnote\plainfootnote

\def\botcaption#1#2\endcaption{\vglue 10pt \cl{\small #1: #2}}
\def\lbotcaption#1#2\endcaption{\vglue 10pt \small #1: #2}
%
%
\def\Refs{\ppar{\large\bf References}\ppar\bgroup\leftskip=25pt
\frenchspacing\parskip=3pt plus2pt\small}       
\def\endRefs{\egroup}
\def\widestnumber#1#2{\relax}
\def\endrefitem{}
\def\refdef#1#2#3{\def#1{\leavevmode\unskip\endrefitem#2\def\endrefitem{#3}}}
\def\ref{\par}
\def\endref{\leavevmode\unskip\endrefitem\par\def\endrefitem{}}
\refdef\key{\noindent\llap\bgroup[}{]\ \ \egroup}
\refdef\no{\noindent\llap\bgroup[}{]\ \ \egroup}
\refdef\by{\bf}{\rm, }
\refdef\manyby{\bf}{\rm, }
\refdef\paper{\it}{\rm, }
\refdef\book{\it}{\rm, }
\refdef\jour{}{ }
\refdef\vol{}{ }
\refdef\yr{(}{) }
\refdef\ed{(}{, Editor) }
\refdef\publ{}{ }
\refdef\inbook{from: ``}{'', }
\refdef\pages{}{ }
\refdef\page{}{ }
\refdef\paperinfo{}{ }
\refdef\bookinfo{}{ }
\refdef\publaddr{}{ }
\refdef\moreref{}{ }
\refdef\eds{(}{, Editors) }
\refdef\bysame{\hbox to 3 em{\hrulefill}\thinspace,}{ }
\refdef\toappear{(to appear)}{ }
\refdef\issue{no.\ }{ }
\newcount\refnumber\refnumber=1
\def\refkey#1{\expandafter\xdef\csname cite#1\endcsname{\number\refnumber}%
\global\advance\refnumber by 1}
\def\cite#1{[\csname cite#1\endcsname]}
\def\Cite#1{\csname cite#1\endcsname}  
\def\key#1{\noindent\llap{[\csname cite#1\endcsname]\ \ }}
%
%
\refkey{A}
\refkey{B}
\refkey{Cp}
\refkey{Ck1}
\refkey{Ck2}
\refkey{Co}
\refkey{E1}
\refkey{E2}
\refkey{EG}
\refkey{EGH}
\refkey{F}
\refkey{G}
\refkey{N}
\refkey{NT}
\refkey{Ta}
\refkey{Th1}
\refkey{Th2}
\refkey{Tr1}
\refkey{Tr2}
\refkey{Tr3}
\refkey{V1} 
\refkey{V2} 
%
%
%
\define\dss{d\.s\.s\.\,}
\define\ind{\operatorname{ind}}
\define\pd#1#2{\frac{\partial#1}{\partial#2}}
\define\jetS{\Cal J^1(S^1)}
\define\gqi{g\.q\.i\.\, }

\define\jet0{j^1(0)}

\topmatter
\title
Generating function polynomials for\\legendrian links
\endtitle
\author
Lisa Traynor
\endauthor
\date June 2001 \enddate
\address{Mathematics Department, Bryn Mawr College\\Bryn 
Mawr, PA 19010, USA}\endaddress
\email{ltraynor@brynmawr.edu}\endemail
\abstract{It is
shown that, in the $1$--jet space of the circle, the swapping and the flyping
procedures, which produce topologically equivalent links, can
produce  nonequivalent legendrian
links.  Each component of
the links considered is legendrian isotopic to the $1$--jet of the
$0$--function, and thus cannot be distinguished by the
classical rotation number or Thurston--Bennequin invariants.
The links are distinguished by calculating
invariant polynomials defined via
homology groups associated to the links through the
 theory of generating functions. The many calculations of
these generating function polynomials
support the belief that these polynomials
carry the same information as a refined version of
Chekanov's first order polynomials which are defined via
the theory of holomorphic curves.
}
\asciiabstract{
It is shown that, in the 1-jet space of the circle, 
the swapping and the flyping
procedures, which produce topologically equivalent links, can
produce  nonequivalent legendrian
links.  Each component of 
the links considered is legendrian isotopic 
to the 1-jet of the
0-function, and thus cannot be distinguished by the
classical rotation number or Thurston-Bennequin invariants. 
The links are distinguished by calculating
invariant polynomials defined via  
homology groups associated to the links through the 
theory of generating functions. The many calculations of
these generating function polynomials 
support the belief that these polynomials
carry the same information as a refined version of
Chekanov's first order polynomials which are defined via
the theory of holomorphic curves.
}
\endabstract
\keywords
Contact topology, contact homology, generating functions,
legendrian links, knot polynomials
\endkeywords
\primaryclass{53D35}
\secondaryclass{58E05}
\endtopmatter

\document

\heading{Introduction}\endheading

 The $1$--jet space of the circle, $\jetS$, is a manifold diffeomorphic
to $S^1
\times
\Bbb R^2$:
$$\jetS= T^*(S^1) \times \Bbb R = \{ (q,p,z) \: q \in S^1, \  p,z \in
\Bbb R
\}.$$
Viewing $S^1$ as a quotient of the unit interval,
$S^1 = \{ q \in [0,1] \: 0 \sim 1 \}$, it is easy to
visualize knots in $\jetS$ as quotients of arcs in $I \times \Bbb
R^2$.
This paper focuses on   two-component links in
$\jetS$ that satisfy a geometrical condition imposed by
the standard contact structure.
This standard contact structure on $\jetS$ is a field of
hyperplanes $\xi$ given by
 $\xi = \ker (dz-pdq).$ There are
  no integral surfaces of this hyperplane distribution; however, there
are integral curves.
An immersed curve $\Cal L$ in $(\jetS, \xi)$ is
{\sl legendrian} if
it is tangent to $\xi$, $T\Cal L \i \xi$.
Functions on $S^1$ with values in $\Bbb R$ give rise to legendrian
knots in
$\left( \jetS, \xi \right)$:
the graph of a smooth $1$--periodic function $f$ in the $(q,z)$--plane
has a lift to an embedded legendrian curve
$j^1(f) := \left\{ (q,p,z) \: z = f(q),\  p = \frac{d}{dq}f(q)
\right\}.$ Notice that $j^1(f)$
is the {\sl 1--jet of the function
$f$}.   Similarly, the graphs of two smooth functions $f,g \: S^1 \to
\Bbb R$ will lift to a legendrian link as long as they have no
points of tangency. In \cite{Tr2}, a link of the form
$j^1(f) \amalg j^1(g)$ is considered where
$f,g\: S^1 \to \Bbb R$ satisfy
$f(q) > 0$ and $g(q)< 0$, for all $q \in S^1$. In particular, it is
shown that
$j^1(f) \amalg j^1(g)$ is an ``ordered" legendrian link.
 A legendrian link $\Lambda_1 \amalg \Lambda_0$ is {\sl ordered}
if $ \Lambda_1 \amalg \Lambda_0$ is not legendrianly equivalent
to $\Lambda_0 \amalg \Lambda_1$:  there is not a smooth $1$--parameter family
of legendrian links $\Theta_t$, $t \in [0,1]$, such that
$\Theta_0 =  \Lambda_1 \amalg \Lambda_0$ and
$\Theta_1 =  \Lambda_0 \amalg \Lambda_1$.  $\Lambda_0 \amalg \Lambda_1$ will be
called the {\sl swap} of $\Lambda_1 \amalg \Lambda_0$.

\midinsert
\cl{\epsfxsize.6\hsize\epsfbox{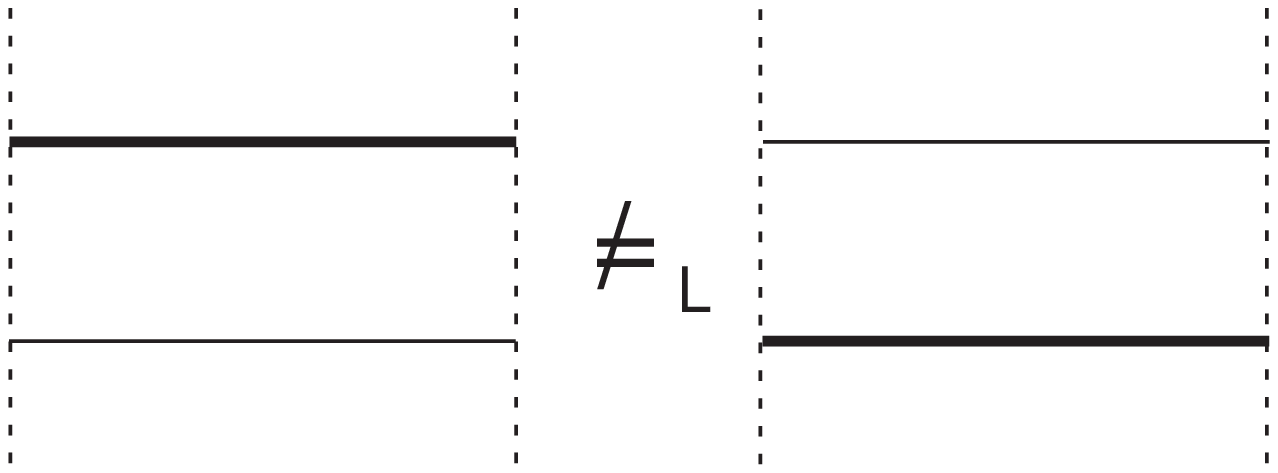}} 
\botcaption{Figure 1.1}{The legendrian link $(0)$ and its 
nonequivalent swap} 
\endcaption
\endinsert

 In this paper, more complicated legendrian links will be studied.
All these  links will be topologically unordered, but many will be
legendrianly ordered.  In addition, ``flyping" moves applied
to these links will  produce
topologically equivalent links that are
legendrianly distinct.
These more complicated links are
constructed as the lifts of  graphs of ``multivalued" $1$--periodic functions.
Namely, in $\{ (q,z) \:   q\in I, z \in \Bbb R\}$,
consider a piecewise smooth arc
with semicubical cusps at the nonsmooth points.  Then the arc has a
well-defined tangent line at each point.  As long as the tangent
line
is never vertical, and there are no self-tangencies, this arc
will have
a lift, with $p$ specified by the slope of the tangent line,
 to an embedded legendrian arc in
$I \times \Bbb R^2$.  As long as the appropriate boundary conditions
are satisfied, this arc will lift to a legendrian knot in $\jetS$.
The graphs of such cusped curves can be seen in Figures 1.5, 1.12.  Notice
that it is not necessary to use broken curves to indicate which is
the upperstrand in the lift:  the third coordinate is determined by
slope.

The basic links considered in this paper can be thought of as closures
of rational tangles, defined by Conway in \cite{Co}, and so it will be
convenient
to label the links
 using the  rational tangle nomenclature
developed in \cite{Tr3}, which is similar to the nomenclature in \cite{Co} and
\cite{A}. With this notation, the ordered link
$j^1(f) \amalg j^1(g)$ considered above will be called the link
$(0)$. More generally, for any $h \geq 0$,  an ``integral" link
$(2h)$ is constructed as the lift  of the graphs of two functions $f,g\:S^1
\to
\Bbb R$ that cross transversally at $2h$ points.
More complicated
links will  be described by length $2n-1$ vectors of the
form
$$(2h_n, v_{n-1}, \dots, 2h_2, v_1, 2h_1), \quad
h_n, v_{n-1}, \dots, h_2, v_1 \geq
1, \text{ and } h_1 \geq 0.
\tag{1.2} $$  The {\sl standard rational legendrian link $(2h_n,
v_{n-1}, \dots, 2h_2, v_1, 2h_1)$}
can be constructed recursively:  for $n =1$, these are the
integral links $(2h)$   described above, and for $n \geq 2$, the
$(2n-1)$--length
link
$(2h_n, v_{n-1}, \dots, 2h_2, v_1, 2h_1)$
is formed from ``vertical and horizontal additions" to the $(2n-3)$--length
 link
$(2h_n,\dots, v_2, 2h_2)$, as shown in Figure 1.3.
This rational tangle nomenclature
is extremely valuable in labeling knots and links in topological knot tables.
Developing a
legendrian version of such nomenclature is obviously useful for labeling
legendrian knots and links
in $\jetS$, but also in $\Bbb R^3$ since, as described by Ng in \cite{N},
a satellite construction glues these links in $\jetS$ into a tubular
neighborhood
of a knot in $\Bbb R^3$ to produce links in $\Bbb R^3$.
\midinsert
\cl{\epsfxsize.9\hsize\epsfbox{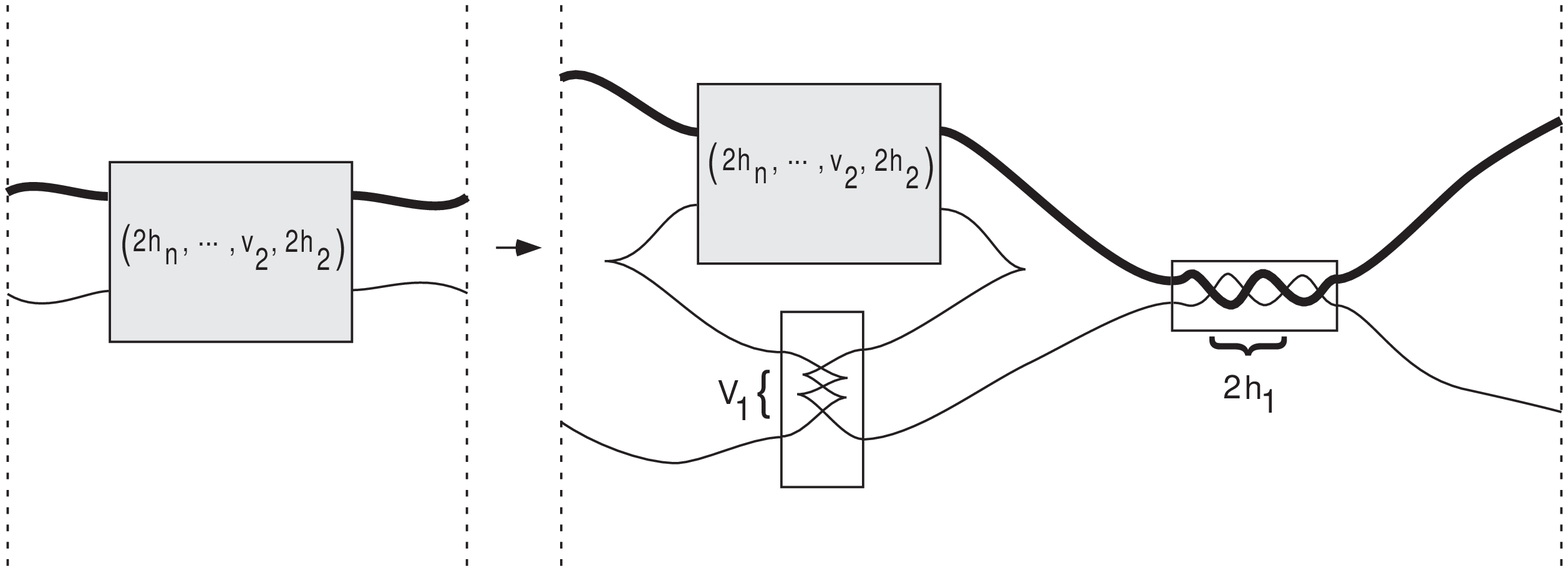}}
\smallskip
\lbotcaption
{Figure 1.3}{The recursive construction of the link 
 $(2h_n, v_{n-1}, \dots, 2h_2, v_1, 2h_1)$  from the link $(2h_n, v_{n-1}, \dots,
2h_2)$}
\endcaption
\endinsert
When the link $(2h_n, \dots, v_1, 2h_1)$ is thought of as a subset of  $I
\times
\Bbb R^2$ without the boundary identification, the legendrian arcs
topologically
form  a rational tangle
that is alternatively known
as the
nonnegative rational number
$$ q := 2h_1 + 1/(v_1 + 1/(2h_2 + 1/(v_2 + \dots + 1/(v_{n-1} + 1/
2h_n)\dots))) \in \Bbb Q.
\tag{1.4} $$
In fact, by Conway's construction,
there is a topological tangle  associated to every
  $q \in \Bbb Q$.   Many of these tangles
do not close up to two-component links in
$\jetS$.  Other $q \in \Bbb Q$
 will correspond to two-component topological links, but will not be
considered here,
since this paper will discuss only
legendrian links where each component is
legendrian isotopic to $j^1(0)$, the $1$--jet of  the $0$--function.  Such
legendrian links will be
called {\sl minimal} legendrian links.  It is shown in \cite{Tr3}
that a minimal legendrian link version of $q \in \Bbb Q$
exists  if and only if $q$ corresponds to a vector of the form in (1.2).
Minimal legendrian   links cannot be
distinguished by examining the legendrian invariants for
the strands
known as the Thurston--Bennequin invariant and  the rotation number.
Background on these classic invariants can be found, for example, in \cite{B},
\cite{E1}, \cite{Ta}.

\midinsert
\cl{\epsfxsize.9\hsize\epsfbox{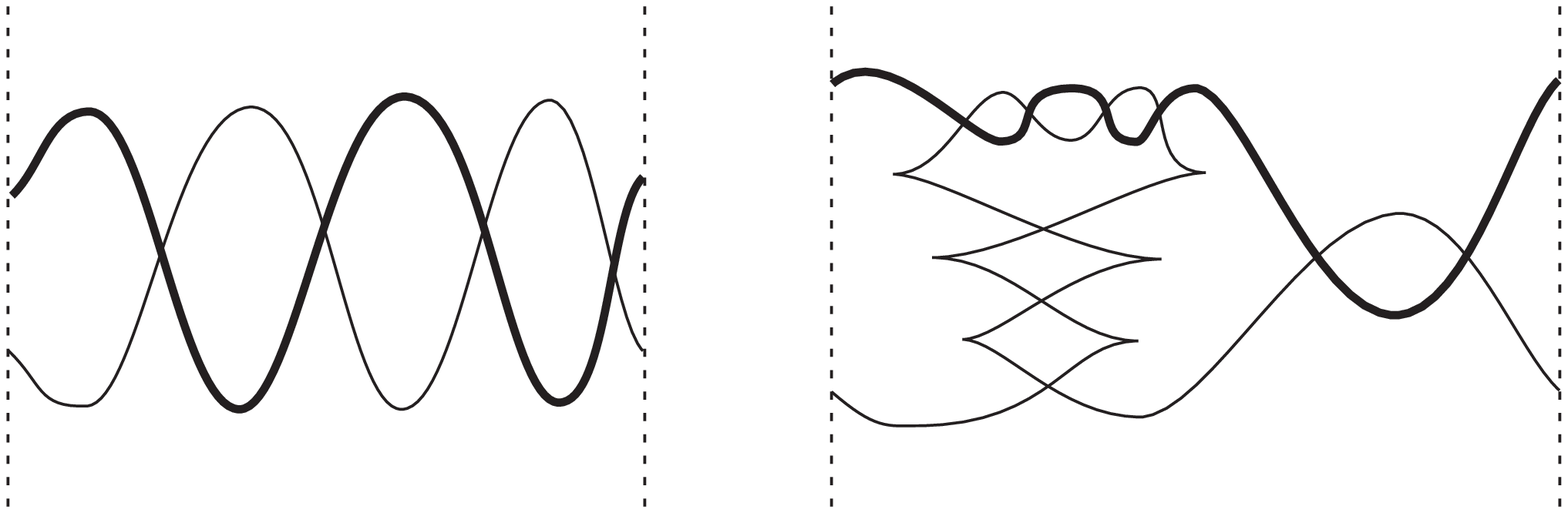}}
\smallskip
\botcaption{Figure 1.5}  Minimal legendrian links $4$ and 
$(4,3,2)= \frac{30}{13}$ 
\endcaption
\endinsert

For all the rational links, changing the order of the components produces
topologically equivalent links.
Although the legendrian link $(0)$ is ordered, it is easy to verify that
the  legendrian  link $(2h)$, for $h \geq 1$, is not ordered.  However, this
is the exception among the  rational legendrian links.

\proclaim{Theorem 1.6} {\rm (See Theorem 7.2)}\qua  Consider the legendrian
link \newline
$L = (2h_n,v_{n-1}, 2h_{n-1}, \dots, v_1,2h_1)$.
 For $n = 1$,  $L = (2h_1)$ is ordered iff $h_1 = 0$.
 For $n \geq 2$,  $L$
is ordered for
all choices of $h_i, v_i$.
\endproclaim

In \cite{Tr2}, the link $(0)$  is shown to be ordered by studying
Viterbo's invariants known as $c_\pm$.  These invariants arise as
``canonical critical values" of a difference of generating functions
associated to the strands of the link.  The links of Theorem 1.6
with $h_1 = 0$ can also be distinguished by studying
$c_\pm$.  All the links will be proven to be ordered by doing a more
in depth analysis of the generating functions:  rather than merely
studying the critical values,
homology groups  $H_k^\pm(L)$, $k \in \Bbb Z$,
for a minimal legendrian link $L$ will be constructed
by examining the relative homology groups of canonical sublevel sets.  This
construction is explained in  Section 3.
From $H_k^\pm[L]$,  {\sl positive}
and {\sl negative homology polynomials}  are legendrian invariants
associated to each
link
$L$:
$$ \Gamma^+(\lambda)[L] = \sum_{k= -\infty}^\infty \dim H_k^+(L) \lambda^k,
\qquad
\Gamma^-(\lambda)[L] = \sum_{k= -\infty}^\infty \dim H_k^-(L) \lambda^k.
\tag{1.7} $$

 A comparison of $\Gamma^+(\lambda)[L]$ and $\Gamma^-(\lambda)[L]$  can
detect that
the  legendrian link $L$ is ordered.
It will be said that  polynomials
$\alpha(\lambda) = \sum_{k=-\infty}^\infty a_k \lambda^k$ and
$\beta(\lambda) = \sum_{k=-\infty}^\infty b_k \lambda^k$ are {\sl 1--shift
palindromic} if
$ \alpha(\lambda) = \lambda \cdot \overline{\beta(\lambda)},$
where $\overline{\beta(\lambda)}$ denotes the palindrome of $\beta$:
$\overline{\beta(\lambda)} = \sum_{k=-\infty}^\infty b_k \lambda^{-k}$.

\proclaim{Theorem 1.8} {\rm (See Corollary 3.17)}\qua  If
$\Gamma^+(\lambda)[L]$ and
$\Gamma^-(\lambda)[L]$ are not $1$--shift palindromic, then the link $L$ is
ordered.
\endproclaim

Theorem 1.6 then follows easily from Theorem 1.8 and the following
calculation which
is proven in Section 6 after developing algebraic topology tools in Section 4
and methods to calculate the indices of critical points in Section 5.

\proclaim{Theorem 1.9} {\rm (See Theorem 6.1)}\qua Consider the
legendrian link \newline $L = (2h_n,v_{n-1}, \dots, v_1,2h_1)$.
Then
$$\align
\Gamma^-(\lambda)\left[  L\right] &=
h_1 + h_2\lambda^{-v_1} + h_3\lambda^{-v_1-v_2} + \dots +
h_n\lambda^{-v_1-v_2-\dots-v_{n-1}},
\\
\Gamma^+(\lambda)\left[  L \right] &=
\cases  \lambda \cdot \Gamma^-(\lambda)\left[
 L \right], & h_1 \geq 1 \\
 (1 + \lambda)  +  \lambda \cdot \Gamma^-(\lambda)\left[ L \right], & h_1 = 0.
\endcases
\endalign$$
\endproclaim

In \cite{NT}, the author and Lenny
Ng find similar calculations of
 a refined version of the Chekanov first order polynomials; these
Chekanov polynomials are invariants
obtained from the differential algebras obtained from the theory of
holomorphic curves,
\cite{Ck2}, \cite{E2}, \cite{EGH}.

For the topological version  of the link $(2h_n, \dots, v_1, 2h_1)$,
$n \geq 2$, in addition to changing the order of the components, there
 are ``flyping" moves that do not change the topological type of the link.
A topological flype occurs when  a portion of the link, represented
by the circle labeled with ``F" in Figure 1.10, is rotated
$180^\circ$ about a vertical axis (a vertical flype), or about a
horizontal axis (a horizontal flype).   For background on topological
flypes, see, for
example, \cite{A}.
This motivates the definition of a {\sl legendrian flype}:  when a crossing is
formed by two   edges emanating from  a legendrian
tangle, represented by the box labeled with ``F" in Figure 1.10, a
legendrian vertical
(horizontal) flype occurs when the tangle is rotated $180^\circ$ about
vertical (horizontal) axis and the crossing  is ``transferred" to the
opposite edges.  This rotation action is not a legendrian
isotopy; so, although the resulting legendrian links are topologically
equivalent, they are potentially not legendrianly equivalent.

\midinsert
\cl{\epsfxsize.9\hsize\epsfbox{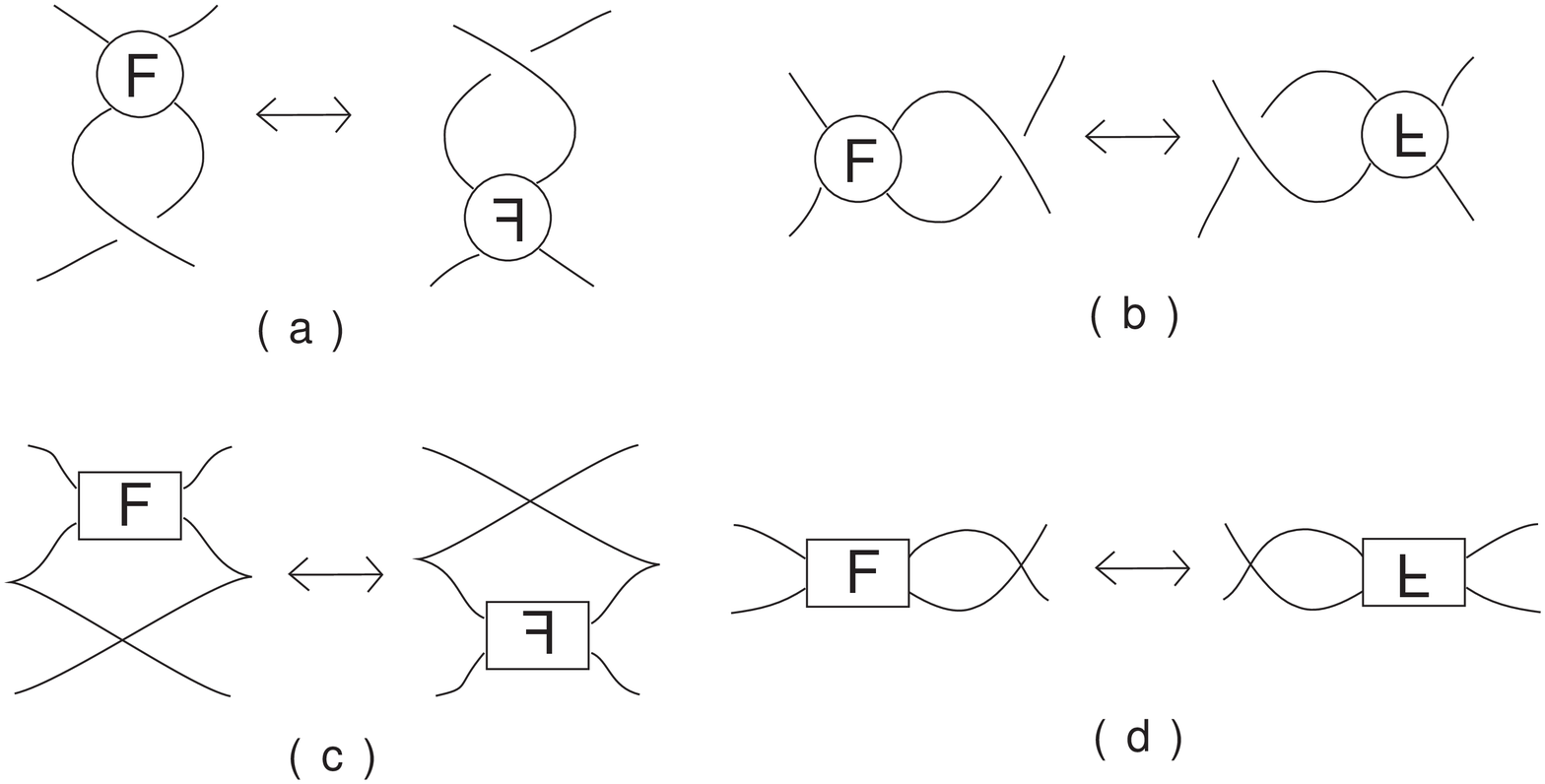}}
\smallskip
\lbotcaption
{Figure 1.10} (a) a topological vertical flype, (b) a
topological horizontal flype, (c) a legendrian vertical flype, (d) a
legendrian horizontal flype.
\endcaption
\endinsert

For each positive
horizontal entry $2h_i$, $i \neq n$, in the legendrian link \newline
$(2h_n, \dots, 2h_2, v_1, 2h_1)$,  it is possible to perform
$0$, $1$, \dots, or $2h_i$ successive horizontal flypes; for each
vertical entry $v_i$, it is possible to perform
$0$, $1$, \dots, or $v_i$ successive vertical flypes.
A flype
at the $2h_i$ entry horizontally flips the tangle  constructed from entries
$2h_n, \dots, 2h_{i+1},
\text{and } v_i$, while a flype at $v_i$ vertically flips the tangle made
from entries $2h_n, \dots,
2h_{i+1}$.
  The flyping
procedure preserves the minimality of the links. The nomenclature
$$ \aligned
\left(2h_n, v_{n-1}^{q_{n-1}}, 2h_{n-1}^{p_{n-1}}, \dots, 2h_2^{p_2},
v_1^{q_1},
2h_1^{p_1}\right),
\quad &q_i \in \{ 0, \dots, v_{i}\}, \\
&   p_i \in \{ 0, \dots, 2h_i\},
\endaligned
\tag{1.11}
$$
will be used to denote the modification of the
 standard  link $(2h_n, \dots, 2h_2, v_1,$\break
$2h_1)$ by $p_i$ horizontal flypes in the $i^{th}$ horizontal
component, and $q_i$
vertical flypes in the $i^{th}$ vertical component.
With this notation, the standard link
$(2h_n,\dots, 2h_2, v_1, 2h_1)$  is written as
$(2h_n, v_{n-1}^0, 2h_{n-1}^0, \dots, 2h_2^0, v_1^0, 2h_1^0)$.  If no
superscript is
specified for an entry of the vector, it will be assumed to be $0$.

First consider a link $L$ that is obtained by applying vertical flypes
to  a standard rational link:
$L = (2h_n, v_{n-1}^{q_{n-1}}, 2h_{n-1}, \dots, v_1^{q_1}, 2h_1),$
$q_i \in \{ 0, \dots, v_i \}.$
Figure 1.12 illustrates some links that  differ
 by vertical flypes. In fact, any link $L$ that is obtained from $L_0 = (2h_n,
v_{n-1},
\dots, v_1, 2h_1)$ by vertical flypes is legendrianly equivalent to $L_0$.
\midinsert   
\cl{\epsfxsize\hsize\epsfbox{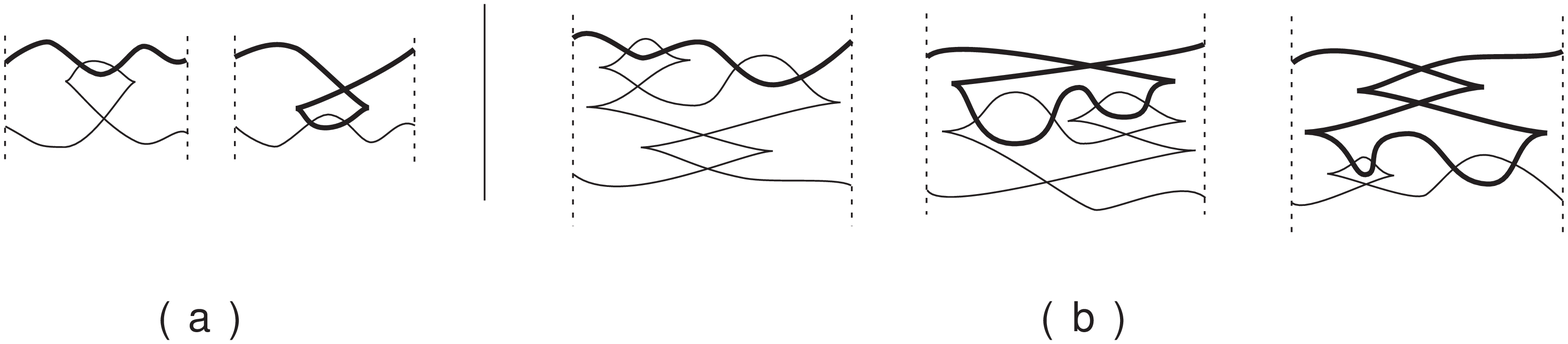}}
\lbotcaption{Figure 1.12}   (a) The equivalent 
legendrian links $(2,1,0)$ and $(2,1^1,0)$;  
(b) The equivalent
legendrian links $(2,1,2,2,0)$, $(2,1,2,2^1,0)$, and 
$(2,1,2,2^2,0)$.
\endcaption
\endinsert

\proclaim{Theorem 1.13} {\rm (See Theorem 2.1)}\qua  Consider the  legendrian
links \newline
$L_0\!=(2h_n, v_{n-1}, 2h_{n-1}, \dots, v_1, 2h_1)$ and
$L_1\!=(2h_n, v_{n-1}^{q_{n-1}}, 2h_{n-1}, \dots, v_1^{q_1}, 2h_1)$.  Then $L_0$
and $L_1$ are legendrianly equivalent.
\endproclaim

Theorem 1.13 is proved in Section 2 by showing that
the $(q,z)$--projections of $L_1$ and $L_0$ are
equivalent through  a
 sequence of ``legendrian planar isotopies"
and ``legendrian Reidemeister moves".

Next consider a link $L$ that is obtained by applying horizontal flypes
to  a standard rational link:
$L = (2h_n, v_{n-1}, 2h_{n-1}^{p_{n-1}}, \dots, v_1, 2h_1^{q_1})$,
$p_i \in \{ 0, \dots, 2h_i \}.$
Figure 1.14 illustrates some links that differ  by horizontal flypes.  In
contrast to the
vertical flyping  situation, it is possible to obtain distinct legendrian links
by horizontal flypes.  For example, the legendrian links
$(2,1,2)$ and $(2,1,2^1)$ are topologically equivalent but
not legendrianly equivalent.  This is a consequence of
calculating the
$\Gamma^-(\lambda)$ or $\Gamma^+(\lambda)$ polynomials.

\midinsert
\cl{\epsfxsize.9\hsize\epsfbox{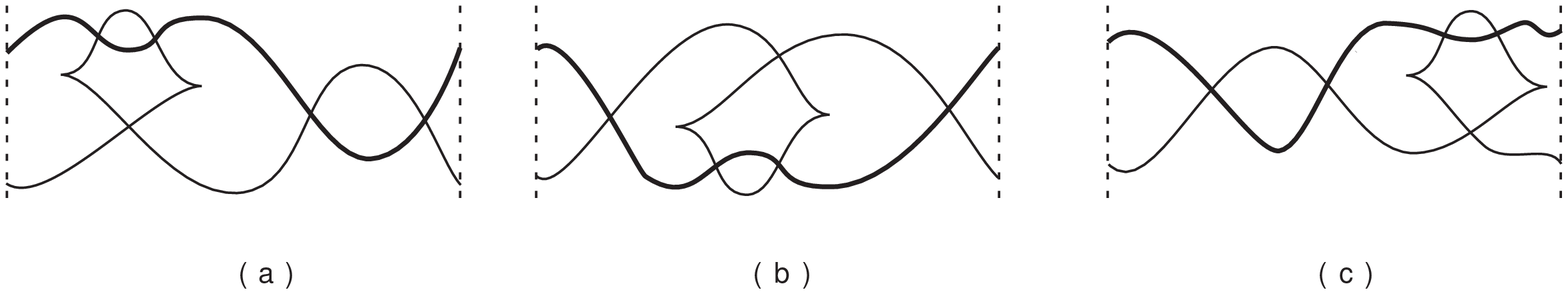}}
\smallskip
\lbotcaption{Figure 1.14} The legendrian links
(a) $(2,1,2^0)$, (b) $(2,1,2^1)$, and (c) $(2,1,2^2)$.  
The link $(2,1,1^0)$ is equivalent to $(2,1,2^2)$, but
distinct from $(2,1,2^1)$. 
\endcaption
\endinsert

\proclaim{Theorem 1.15} {\rm (See Theorem 6.2)}\qua   Consider the  legendrian
link
$$
L = \left(2h_n,v_{n-1}, 2h_{n-1}^{p_{n-1}}, \dots, v_1,2h_1^{p_1}\right),
 \quad p_i \in \{0, \dots, 2h_i \}.
$$
  For $j = 1, \dots, n-1$, let
$\sigma(j) = 1 + \sum_{i = 1}^j p_i \mod 2$.
Then
$$\align
\Gamma^-(\lambda)\left[   L  \right] &= h_1 + \sum_{i=2}^n
h_i\lambda^{(-1)^{\sigma(1)}v_1 + (-1)^{\sigma(2)} v_2+ \dots +
(-1)^{\sigma(i-1)}v_{i-1}} ,
\\
\Gamma^+(\lambda)\left[  L  \right] &=
\cases
\lambda \cdot  \Gamma^-(\lambda)\left[ L \right], & h_1 \geq 1 \\
(1+\lambda) +  \lambda \cdot  \Gamma^-(\lambda)\left[L \right], & h_1 = 0.
\endcases
\endalign$$
\endproclaim

\remark{\bf Remark/Question 1.16}  Notice that
given
$L_0 = ( 2h_n, v_{n-1}, 2h_{n-1}^{p_{n-1}}, \dots, v_1,$\break 
$2h_1^{p_1})$ and
$L_1 = \left( 2h_n, v_{n-1}, 2h_{n-1}^{w_{n-1}}, \dots, v_1, 2h_1^{w_1}
\right)$,
if
$p_i \equiv w_i \mod 2$, for all $i$, then $\Gamma^-(\lambda)[L_0] =
\Gamma^-(\lambda)[L_1]$, and
$\Gamma^+(\lambda)[L_0] = \Gamma^+(\lambda)[L_1] $  .
 This is a natural condition on $p_1$:
if  $w_1 \equiv p_1 \mod 2$, and $w_i = p_i$ when $i \neq 1$,   then $L_0$
and $L_1$
are, in fact, equivalent.  Due to the boundary identification,
``double" horizontal flypes  are eqivalent to a rotation. It would
be interesting to know if only the parity is important for other $p_i$, $i
\neq 1$.
For example, are the links $(2,1,2,1,0)$ and $(2,1,2^2,1,0)$   illustrated in
Figure 1.17 legendrianly equivalent? 
\endremark

\midinsert
\cl{\epsfxsize.6\hsize\epsfbox{ 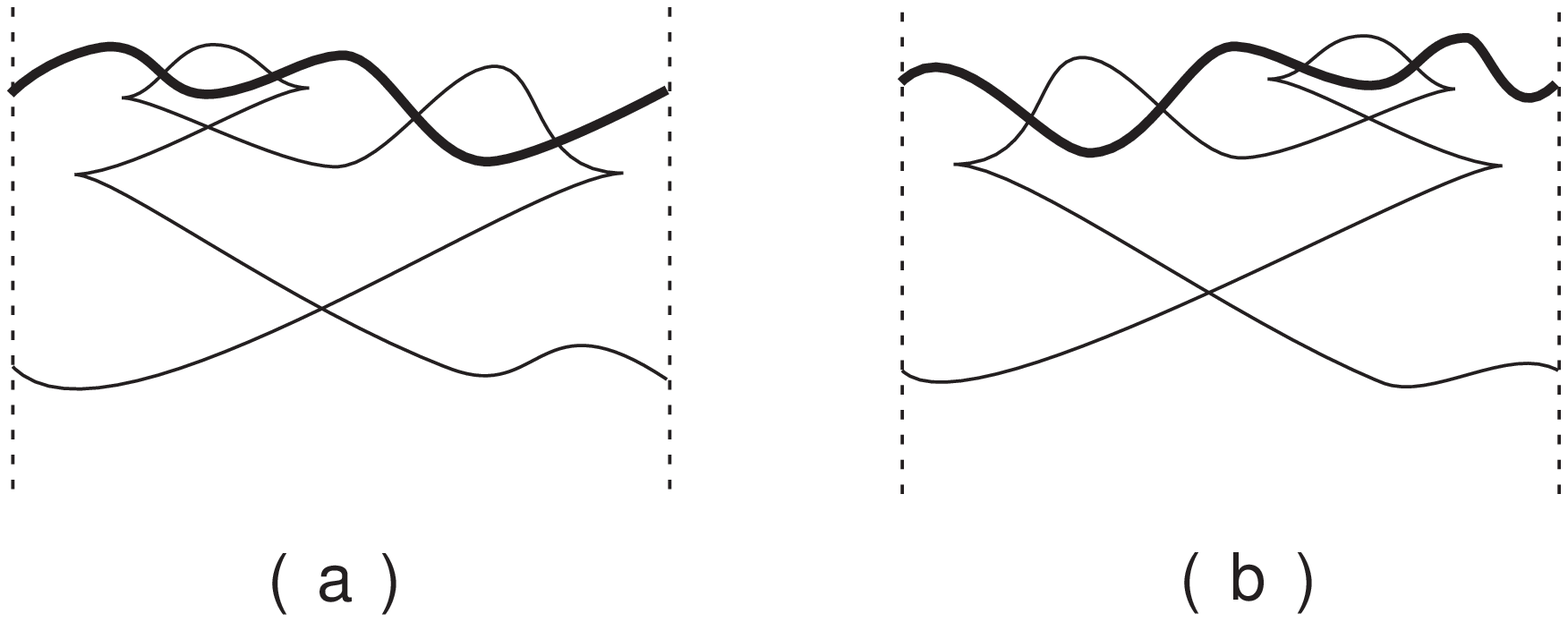}}
\lbotcaption{Figure 1.17}  The legendrian links (a) $(2,1,2,1,0)$, and (b)
$(2,1,2^2,1,0)$.  They have the same polynomials.  Are they equivalent?
\endcaption
\endinsert

\remark{\bf Remark/Question 1.18} In view of Theorem 1.6, it is natural to
ask if, for $n \geq 2$, all the possible horizontal flypes of
$L_n = (2h_n, \dots, v_1, 2h_1)$ are also ordered.  This is true
when $n = 2$: for $n = 2$,
 a flype is possible only when $h_1 \geq 1$, and, in this case,
 it is easy to check
that $L_2^1 = (2h_2, v_1, 2h_1^1)$ is equivalent to the
swap of $L_2^0 = (2h_2, v_1, 2h_1^0)$.
However,   for $n \geq 3$, there are examples
where the polynomials cannot detect if the link is ordered.  For
example,  the link $L_3 = (2,2,2^1,1,2^1)$, illustrated in Figure 1.19,
has $\Gamma^-(\lambda)[L_3] = \lambda^{-1} + 1 + \lambda$, and, since
$\Gamma^-(\lambda)$ is palindromic, Proposition 7.1 implies that  $L_3$ is
potentially
unordered.
 {\it Is $L_3$ equivalent to its swap?}
It is interesting to note that ``slight modifications" of $L_3$
result in links that are ordered.  For example,
$L_3^\prime = (4,2,2^1, 1, 2^1)$ and
$L_3^{\prime\prime} = (2,1,2^1, 1, 2^1)$ satisfy
$\Gamma^-(\lambda)[L_3^\prime] = 2\lambda^{-1} + 1 + \lambda$,
$\Gamma^-(\lambda)[L_3^\prime] = 2 + \lambda$, and, since these polynomials
are not palindromic,  $L_3^\prime$ and
$L_3^{\prime\prime}$ are each ordered. 
\endremark
\midinsert
\cl{\epsfxsize.9\hsize\epsfbox{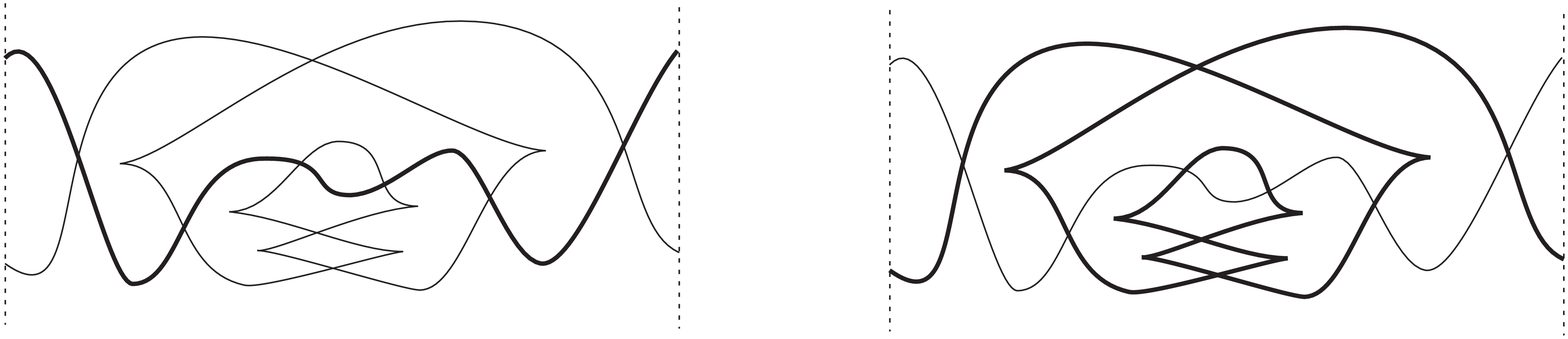}}
\lbotcaption{Figure 1.19} The link $(2,2,2^1,1,2^1)$ and its swap
$\overline{(2,2,2^1,1,2^1)}$.  Their polynomials are
the same.  Are the links equivalent?
\endcaption
\endinsert

\remark{\bf Remark/Question 1.20} {\sl See Corollary 2.3 and Proposition
7.3}\qua The flyping procedure can be
thought of as a generalization of the swapping procedure: as shown
in Corollary 2.3, for $h_1 \geq 1$, if
$$\gather
 L_0 = (2h_n, v_{n-1}, \dots, 2h_2, v_1, 2h_1) \quad\text{and}\quad
 L_1 = (2h_n, v_{n-1}, \dots, 2h_2, v_1, 2h_1^{1}), \\
\text{ where }
v_1, \dots, v_{n-2} \equiv 0 \mod 2,
\endgather
$$
then   $L_0$ and $L_1$ are swaps of one another. This motivates
two questions.  {\it Is this statement true without the hypothesis on the
parity of $v_i$?};  in particular, { is $(2,1,2,1,2^1)$ equivalent to the
swap of $(2,1,2,1,2)$?}  In fact, $(2,1,2,1,2^1)$ is equivalent to the swap
of $(2,1,2^2,1,2)$, and thus this question is closely related to
Remark/Question 1.16.
Secondly, {\it Is the analog of this statement true for horizontal flypes
of $L_0$?}
Namely, for $h_1 \geq 1$, consider
$$\align
 M_0 &= (2h_n, v_{n-1}, 2h_{n-1}^{p_{n-1}}, \dots, 2h_2^{p_2}, v_1,
2h_1^{p_1}), \\
   M_1 &= (2h_n, v_{n-1}, 2h_{n-1}^{p_{n-1}}, \dots, 2h_2^{p_2}, v_1,
2h_1^{p_1+1}).
\endalign
$$
Let $\overline{M_0}$ be the swap of $M_0$. It is shown in Proposition 7.3 that
$\Gamma^-(\lambda)[\overline{M_0}] = \Gamma^-(\lambda)[M_1]$,  and
$\Gamma^+(\lambda)[\overline{M_0}] =
\Gamma^+(\lambda)[M_1].$    {\it Are $\overline{M_0}$ and
$M_1$ equivalent?}
In particular,  are  $\overline{(2,2,2^1,2,2)}$ and
$(2,2,2^1,2,2^1)$ equivalent?  See Figure 1.21.  
 \endremark
\midinsert
\cl{\epsfxsize.9\hsize\epsfbox{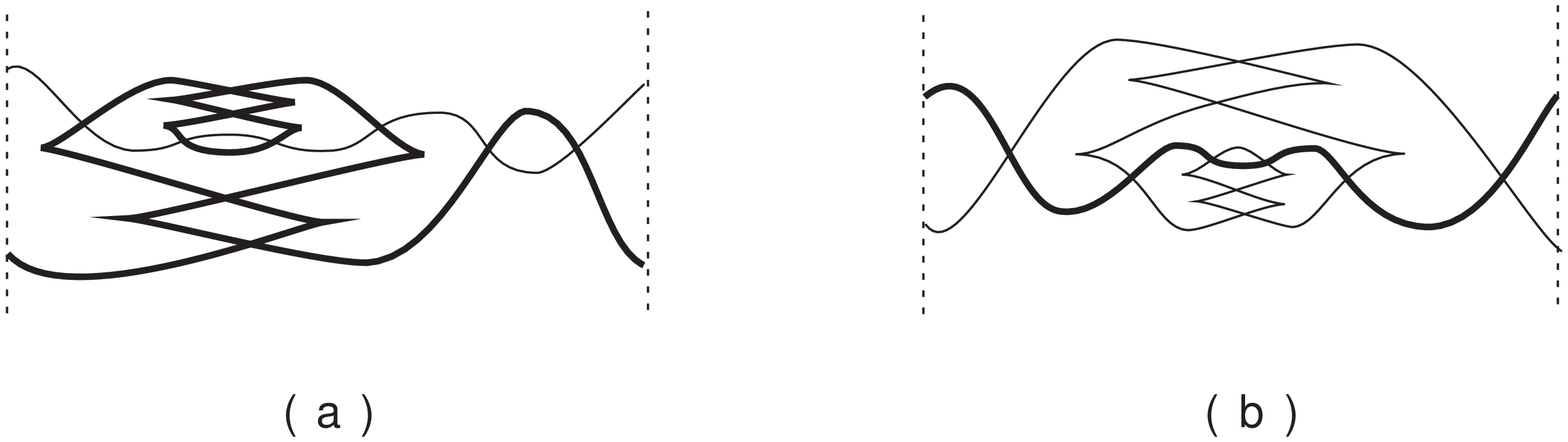}}
\lbotcaption{Figure 1.21}  The legendrian links
(a) $\overline{(2,2,2^1,2,2)}$, (b) $(2,2,2^1,2,2^1)$.  These links
have the same   polynomials.  Are they equivalent?   
\endcaption
\endinsert

\remark{\bf Remark/Question 1.22}
The question of the equivalence of the links\nl $\overline{(2,2,2^1,2,2)}$ and
$(2,2,2^1,2,2^1)$  mentioned in the previous remark is closely related to the
question of the equivalence of the links
$L_0 = (2,1,2^1, 1, 0)$ and $L_1 = (2,1,2^1, 1^1, 0)$;
see Figure 1.23.  Notice that
$L_1$ only differs from $L_0$ by
a vertical flype, but $L_0$ is not standard, and thus
Theorem 1.13 does not imply they are equivalent.  Remark 6.3
explains that the $\Gamma^+(\lambda)$ and $\Gamma^-(\lambda)$ polynomials will
never be able to distinguish two links that only differ by
vertical flypes.
It would be interesting to know if
it is ever possible to obtain distinct legendrian links by
a vertical flype: {\it  Do there exist links of the
form
$L_0=(2h_n, v_{n-1}^{q_{n-1}}, 2h_{n-1}^{p_{n-1}}, \dots, v_1^{q_1},
2h_1^{p_1})$
and
$L_1=(2h_n, v_{n-1}^{w_{n-1}}, 2h_{n-1}^{p_{n-1}}, \dots, v_1^{w_1},
2h_1^{p_1})$
such that $L_0$ and $L_1$ are not legendrianly equivalent?}  Lemma 2.1.1
implies that if
$w_i = q_i$,
when $i \neq n-1$, then $L_0$ and $L_1$ are equivalent.

\endremark

\midinsert
\cl{\epsfxsize.6\hsize\epsfbox{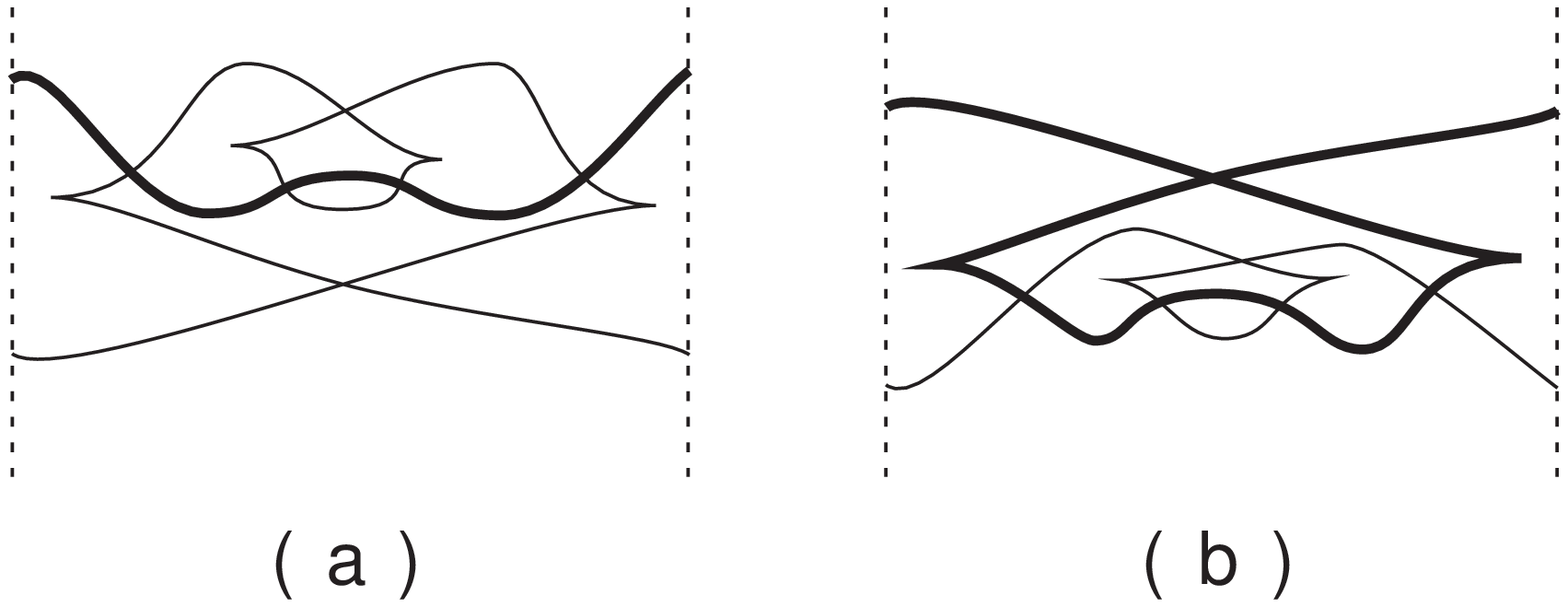}}
\lbotcaption{Figure 1.23}  The links (a) $(2,1,2^1, 1, 0)$, and
(b) $(2,1,2^1, 1^1, 0)$.  They have the same polynomials.  Are
they equivalent?
\endcaption
\endinsert

The following summarizes how many different legendrian representations
of a given rational link type can be constructed from
the swap and the flype operations.  It would be interesting to know
if there are other minimal legendrian versions of these links.

\proclaim{Theorem 1.24} {\rm (See  Theorems 7.2, 7.4, 7.6)}\qua
 Consider the topological link
$$L_n = (2h_n, v_{n-1}, 2h_{n-1}, \dots, 2h_2, v_1, 2h_1), \quad
h_n, v_{n-1}, \dots, h_2, v_1 \geq 1, \quad h_1 \geq 0.$$
\roster
\item
If $n = 2$, there are at least $2$ legendrianly distinct minimal
links that are topologically
equivalent to $L_2$.
\item If $n = 3$, and either $h_1 = 0$, $h_2 \neq h_3$, or $v_2 \neq 2v_1$,
then there are at least $4$  legendrianly distinct minimal
links that are topologically
equivalent to $L_3$.
\item For $n \geq 4$,  if $\max\{ h_1, 1 \} \cup \{ h_i \}_{i=2}^n $
form a set of order
$n$   such that the sums of all its $2^n$ subsets are distinct,
there are at least $2^{n-1}$ minimal legendrian
links that are topologically equivalent  to $L_n$.
\endroster
\endproclaim

When $h_1 \geq 1$, these $2^{n-1}$
different legendrian links arise by looking at
$$(2h_n, v_{n-1}, 2h_{n-1}^{p_{n-1}}, \dots, 2h_2^{p_2}, v_1, 2h_1^{p_1})
\tag{1.25} $$
for  $p_i \in \{ 0, 1 \}$,
$i = 1, \dots, n-1$.  When $h_1 = 0$, the variations
are obtained by the original and the swap  of
each of the $2^{n-2}$ links in (1.25).

\remark{\bf Remark/Question 1.26}  The condition that $h_1 = 0$, $h_2 \neq
h_3$, or
$v_2 \neq 2v_1$ when $n=3$, or that, for $n \geq 4$, $\{ h_i \}$ or $\{ h_i
\} \cup
\{ 1 \}$ form a set of order $n$ with ``distinct subset sums" guarantees that
distinct polynomials are associated to the
$2^{n-1}$   links in (1.25).
In contrast, consider
$$L^{(0,1,0)} = (2,1,2^0,1,2^1,1,2^0), \quad
L^{(0,1,1)} = (2,1,2^0,1,2^1,1,2^1);
$$
see Figure 1.27.
By Theorem 1.15,
$$\Gamma^-(\lambda)[ L^{(0,1,0)}] = \lambda^{-1} + 2 + \lambda =
\Gamma^-(\lambda)[L^{(0,1,1)} ].$${\it Are $L^{(0,1,0)}$ and
$L^{(0,1,1)} $ legendrianly  equivalent?}
It is interesting to note  that  variations of $L^{(0,1,0)}$ and
$L^{(0,1,1)}$ with   $h_2 \neq h_4$ will produce different
polynomials.  For example, if
$$M^{(0,1,0)} = (2,1,2^0,1,4^1,1,2^0), \quad
M^{(0,1,1)} = (2,1,2^0,1,4^1,1,2^1),
$$
then $\Gamma^-(\lambda)[ M^{(0,1,0)}] = 2\lambda^{-1} + 2 + \lambda$, while
$\Gamma^-(\lambda)[M^{(0,1,1)} ] = \lambda^{-1} + 2 + 2\lambda$. 
\endremark
\midinsert
\cl{\epsfxsize.9\hsize\epsfbox{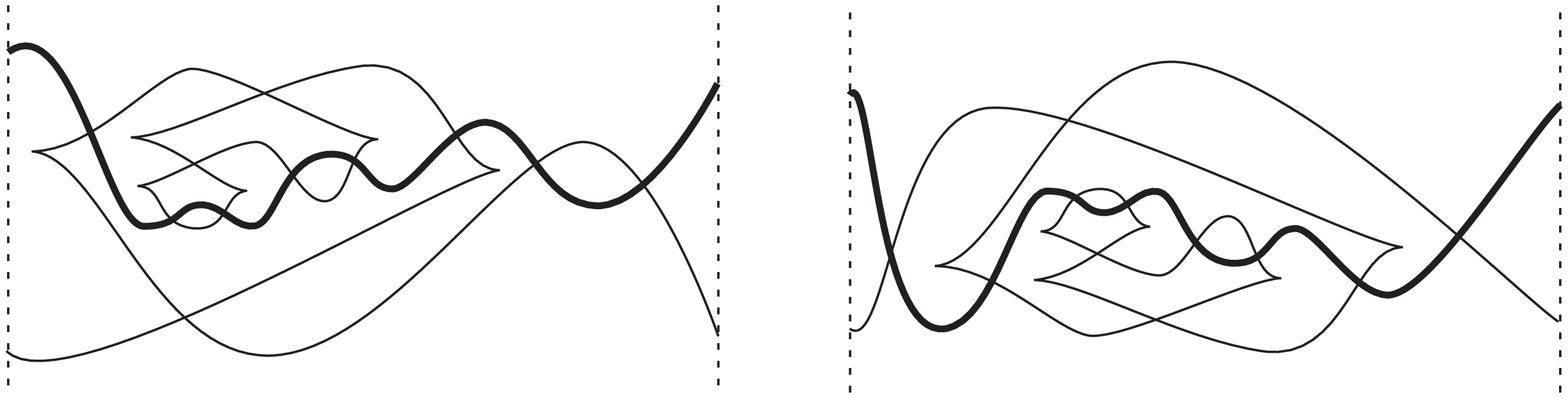}}
\lbotcaption{Figure 1.27}  The legendrian links
$L^{(0,1,0)} = (2,1,2^0,1,2^1,1, 2^0)$ and
 $L^{(0,1,1)} =(2,1,2^0,1,2^1,1, 2^1)$.  They have the
same   polynomials.  Are they equivalent?
\endcaption
\endinsert

\remark{\bf Remark 1.28}  For certain choices of $v_i$, it is possible that
the $2^{n-1}$  links \newline
$L_n = (2h_n, v_{n-1}, 2h_{n-1}^{p_{n-1}}, \dots, v_2, 2h_1^{p_1})$ have
distinct polynomials without the hypothesis that $\{h_1, \dots, h_n\}$
form a    set of order $n$ with distinct subset sums.  For example,  consider
$L_4 = (2,3,2,2,2,1,2).$
In this example, each flype gives a polynomial containing
a  different set of powers of $t$.  More generally, given $v_1, \dots, v_n$,
if the $2^{n-1}$ sets
$$\align
\{ &(-1)^{\sigma(1)}v_1,\   (-1)^{\sigma(1)}v_1 +  (-1)^{\sigma(2)}v_2,
\dots, \\
 &(-1)^{\sigma(1)}v_1 +  (-1)^{\sigma(2)}v_2 +  \dots  +
(-1)^{\sigma(n-1))}v_{n-1} \},
\quad \sigma\: \{ 1, \dots, n-1 \} \to \Bbb Z_2
\endalign$$
are distinct, then any choice of $h_i$ will produce $2^{n-1}$ different
$\Gamma^-(\lambda)$ polynomials.

\endremark

In Section 6,   the polynomials are calculated for rational
links, their flypes, and for the usually nonrational
  ``connect sums" of such links.  Since, up to legendrian
isotopy, the connect sum may depend on the choice of where
the links are cut into tangles, a standard position for
 cutting the links will be chosen.  Namely, the connect sum
$L_1 \# L_2$  is defined as the closure of the connect sum of the
legendrian rational tangles  $:L_1:$ and $:L_2:$, which
are constructed analogously to the links $L_i$.  This construction is
illustrated in Figure 1.29
where, if $L_1$ denotes the link $(2h_n, \dots, 2h_1)$, then $\: L_1 \:$
corresponds
to  Figure 1.3 except considered as a tangle rather than closed to a link.

\midinsert
\cl{\epsfxsize.9\hsize\epsfbox{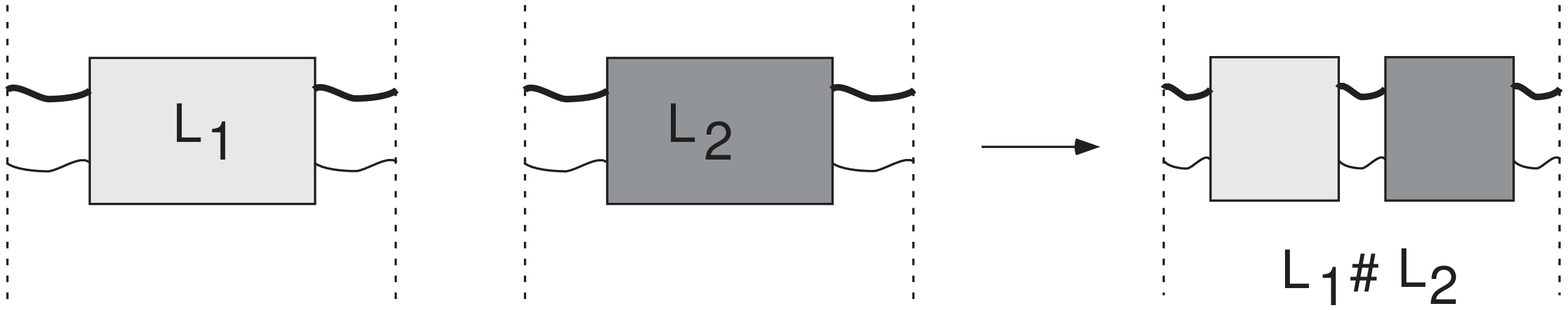}}
\botcaption{Figure 1.29}  The construction of the connect sum $L_1 \# L_2$
\endcaption
\endinsert

\proclaim{Theorem 1.30} {\rm (See Theorem 6.4)}\qua  Consider the  legendrian
links
$$
  L_1 = (2h_n, v_{n-1}, 2h_{n-1}^{p_{n-1}}, \dots, v_1, 2h_1^{p_1}),$$
$$ L_2 = (2k_m, u_{m-1}, 2k_{m-1}^{w_{m-1}}, \dots, u_1, 2k_1^{w_1}).
 $$
Then
$$\align
\Gamma^-(\lambda)[L_1 \# L_2] &= \Gamma^-(\lambda)[L_1] +
\Gamma^-(\lambda)[L_2];\\
\Gamma^+(\lambda)[L_1 \# L_2] &=
\cases
\Gamma^+(\lambda)[L_1] + \Gamma^+(\lambda)[L_2], & h_1, k_1 \geq 1 \\
\Gamma^+(\lambda)[L_1] + \Gamma^+(\lambda)[L_2] - (1+\lambda), & \text{ else.}
\endcases
\endalign$$
\endproclaim

\remark{\bf Remark/Question 1.31}  Theorem 1.30
gives many examples of topogically equivalent, nonrational minimal links that are legendrianly distinct; it also
raises some interesting questions.  For example,  {\it Are the
legendrian links\nl
$ (2,1,2) \# (2,1,2)$  and  $(2,1,2^2) \# (2,1,2)$
equivalent?}  See Figure 1.32. Notice that the rotation that made
$(2,1,2)$ and $(2,1,2^2)$ equivalent is no longer possible. 
\endremark

\midinsert
\cl{\epsfxsize.9\hsize\epsfbox{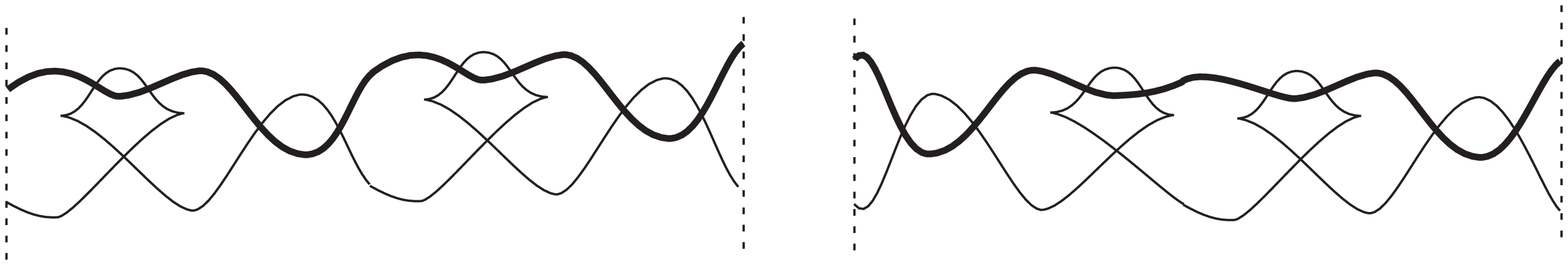}}
\lbotcaption{Figure 1.32}  The legendrian links
$(2,1,2) \# (2,1,2)$, and  $(2,1,2^2) \# (2,1,2)$.
They have the same polynomials.  
Are they equivalent?
\endcaption
\endinsert

Lastly, notice that the nomenclature for links in $\jetS$
easily lends itself to nomenclature for legendrian {\it knots} in
$\jetS$.   In analog with
$(1.2)$, length $2n-1$ vectors of the
form
$$ \aligned
\left(2h_n, v_{n-1}^{q_{n-1}}, 2h_{n-1}^{p_{n-1}}, \dots, 2h_2^{p_2},
v_1^{q_1},
2h_1^{p_1} - 1\right), \\  h_n, v_{n-1}, \dots, h_2, v_1, h_1 \geq
1, \quad
&q_i \in \{ 0, \dots, v_{i}\}, \\
&   p_i \in \{ 0, \dots, 2h_i\},
\endaligned
\tag{1.33}
$$
give rise to legendrian knots.  The technique of
generating functions, as used  in this paper, no longer applies.  The above
results about links
raise many interesting questions about  such knots.  For example,
{\it Are the legendrian knots $(2,1,2,1,1)$ and $(2,1,2^1,1,1)$
legendrianly equivalent?}
See Figure 1.34.

\midinsert
\cl{\epsfxsize.7\hsize\epsfbox{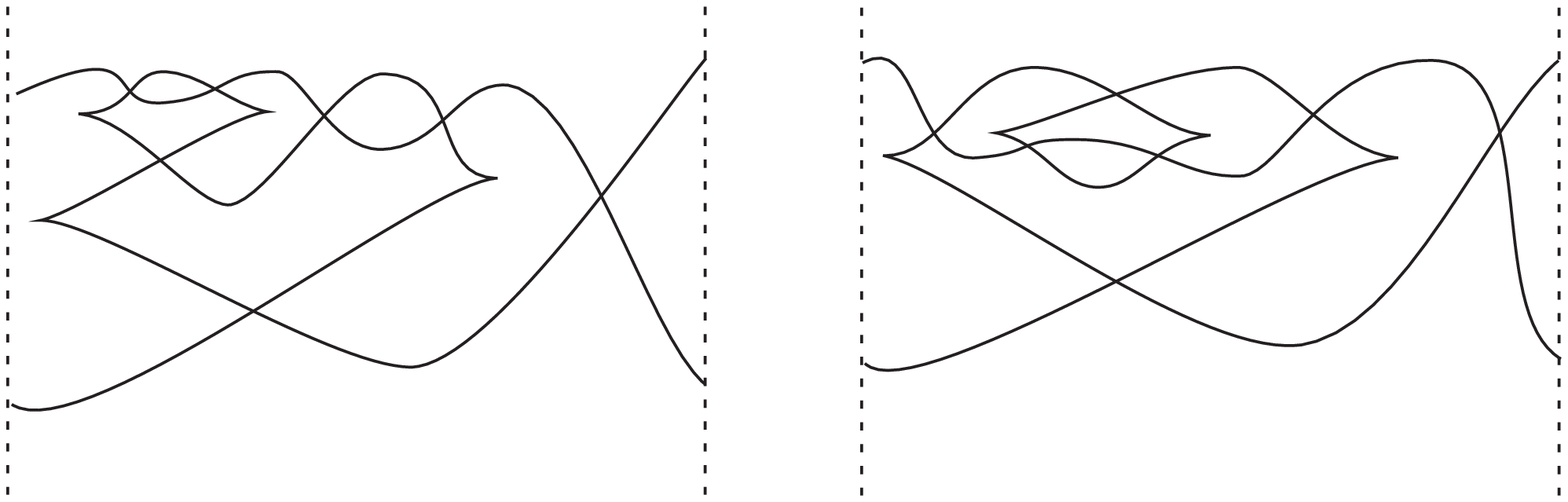}}
\lbotcaption{Figure 1.34} The legendrian knots
$(2,1,2,1,1)$ and $(2,1,2^1,1,1)$.  They have the same
classical invariants.  Are they equivalent?
\endcaption
\endinsert

In \cite{NT}, it is shown that sometimes the differential algebra
approach can say something about some of the above questions.
It is a topic for further study to understand if the generating
function approach can be further refined to capture as many
invariants as the holomorphic curve approach.


\heading{Equivalent Vertical Flypes}\endheading

In this section, it is shown that it is not possible to
produce a nonequivalent legendrian link  by performing
vertical flypes to a standard link.  Recall the terminology
$(2h_n, v_{n-1}^{q_{n-1}}, 2h_{n-1}, \dots, v_1^{q_1}, 2h_1)$,
$q_i \in \{ 0, \dots, v_i \}$,
introduced in (1.11).

\proclaim{Theorem 2.1}  Consider the legendrian links
$L_0=(2h_n, v_{n-1}, 2h_{n-1}, \dots, v_1,$\break $2h_1)$ and
$L_1=(2h_n, v_{n-1}^{q_{n-1}}, 2h_{n-1}, \dots, v_1^{q_1}, 2h_1)$.  Then $L_0$
and $L_1$ are legendrianly equivalent.
\endproclaim

 To prove Theorem 2.1, it suffices to show that the $(q,z)$--projections
of the links  are equivalent by a sequence of ``legendrian planar
isotopies"
and ``legendrian (Reidemeister) moves"; see Figure 2.2. (For background on
the topological Reidemeister moves, see, for example, \cite{A}.) A
{\sl legendrian planar isotopy} is a planar isotopy that does not introduce
cusps or vertical tangents. Each of the  {\sl legendrian type 1$^\pm$ moves}
are analogous to  one of the  type I topological Reidemeister moves:  one
additional crossing and two additional cusps are introduced into the
projection.  The  {\sl legendrian type 2 moves} are analogous to  the  type II
Reidemeister moves:  two new crossings are introduced into the  projection
after a cusp crosses a noncusped segment.  Note that the relative slopes of
the cusp and the segment determine if the cusped segment passes
over or under the noncusped segment.  Lastly, the {\sl legendrian type 3 move}
is analogous to one of the type III Reidemeister moves:  a strand is slid
from one side of a crossing to the other.

\midinsert
\cl{\epsfxsize.9\hsize\epsfbox{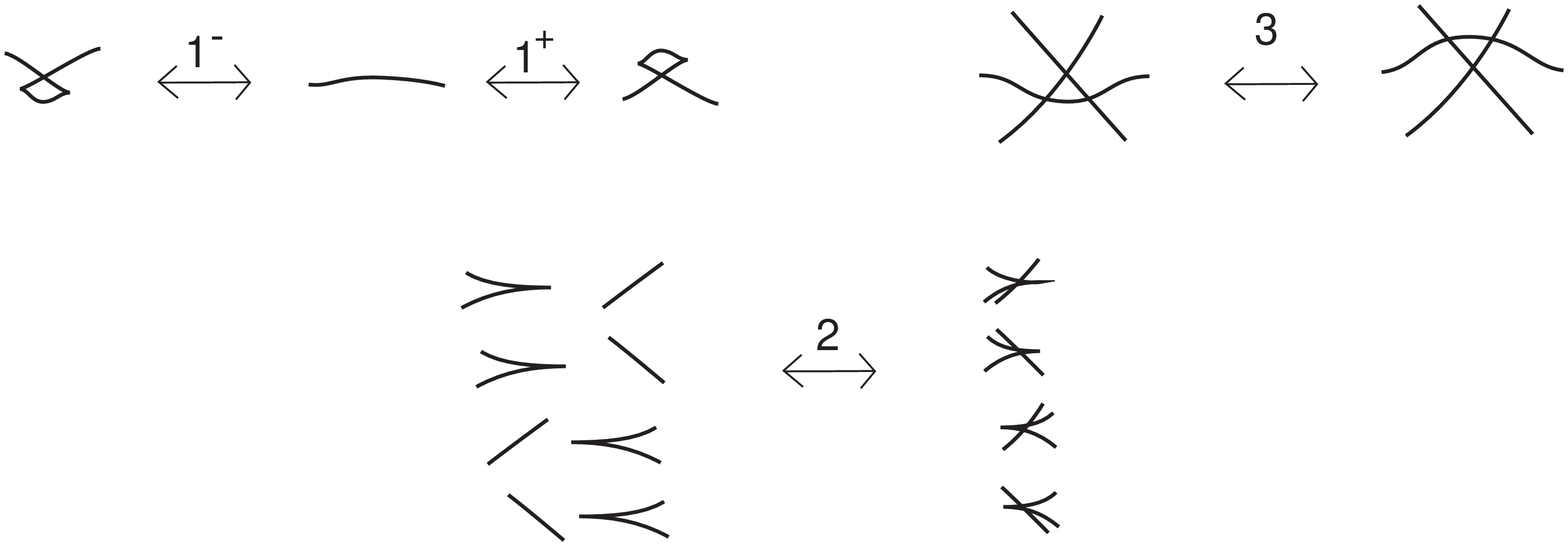}}
\botcaption{Figure 2.2}  The legendrian Reidemeister moves
\endcaption
\endinsert

For the proof of Theorem 2.1, it will be useful to
introduce the notion of a {\sl legendrian tangle}, \cite{Tr3}.
A legendrian tangle consists of two disjoint legendrian arcs
$\Lambda_1, \Lambda_0 \subset \Cal J^1([0,1])$, where
$\Cal J^1([0,1])$ denotes  the $1$--jet space of
the interval $[0,1]$, with
 $\partial \Lambda_1, \partial \Lambda_0 \subset
\{q = 0 \} \cup \{q = 1 \}$.  Legendrian tangles $T_1, T_2$
are equivalent if their $(q,z)$--projections
are equivalent by a sequence of
legendrian planar
isotopies
and legendrian  moves  supported
in $(0,1) \times \Bbb R$.  The nomenclature
$:2h_n, v_{n-1}, \dots, 2h_1:$ will be used
to denote the legendrian tangle constructed using
the same recursive procedure   used to construct
the legendrian link $(2h_n, v_{n-1}, \dots, 2h_1)$.

\demo{Proof of Theorem 2.1}
For $n \geq 2$, the desired equivalence
of the links\newline $(2h_n, v_{n-1}, \dots, v_1, 2h_1)$ and
$(2h_n, v_{n-1}^{q_{n-1}}, \dots, v_1^{q_1}, 2h_1)$
will follow  from a proof that the
tangles
$:2h_n, v_{n-1}, \dots,  v_1, 2h_1:$ and
$:2h_n, v_{n-1}^{q_{n-1}}, \dots,  v_1^{q_1}, 2h_1:$
are equivalent for all choices of $q_i \in \{ 0, \dots, v_i\}$,
$i = 1, \dots, n-1$.  The equivalence of the tangles
will be proved by  induction on $n$.  Lemma 2.1.1 proves
the base case  of $n = 2$.

\proclaim{Lemma 2.1.1}
The legendrian tangles $:2h_2, v_1, 2h_1:$ and
$:2h_2, v_1^{q_1}, 2h_1:$ are equivalent for any $q_1 \in \{ 0, \dots,
v_1\}$.
\endproclaim

\demo{Proof}  It  suffices to show that
 the tangles
$:2h, 1, 0:$  and
$:2h, 1^{1}, 0:$ are equivalent for all $h \geq 1$.  This
will be shown
by an induction argument on $h$.  Figure 2.1.1.1
outlines the legendrian moves that prove
the base case of $h=1$.
\midinsert
\cl{\epsfxsize.9\hsize\epsfbox{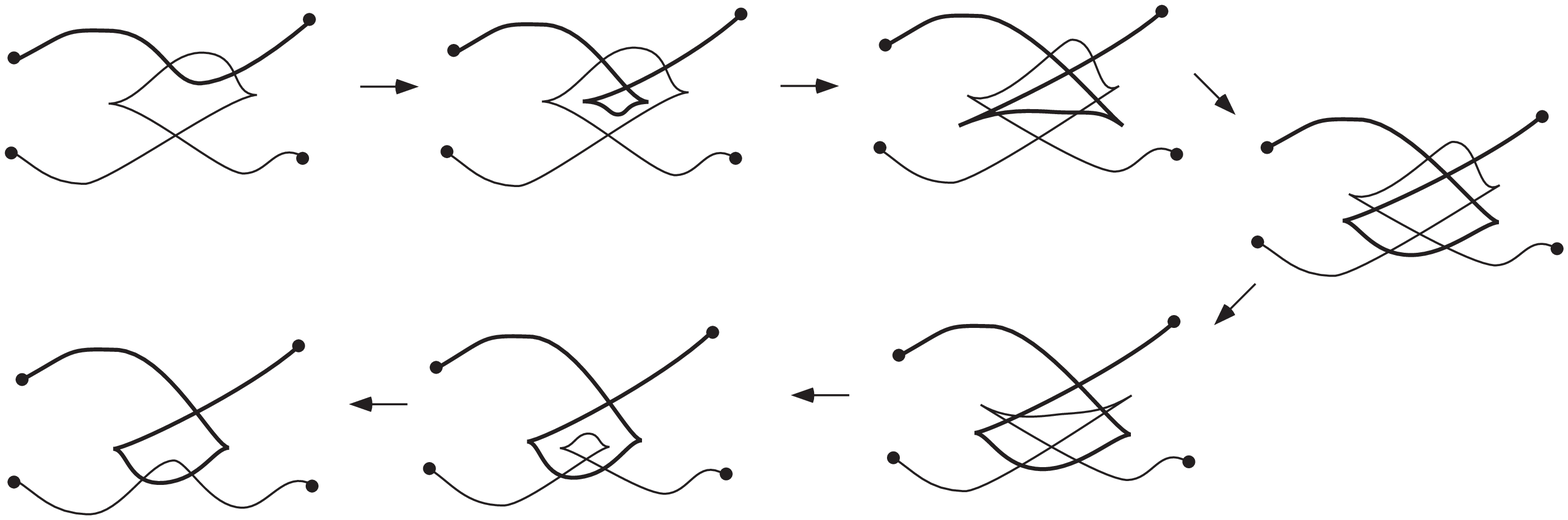}}
\botcaption
{Figure 2.1.1.1}  The equivalence of the tangles
$:2,1,0:$ and $:2,1^1, 0:$ 
\endcaption
\endinsert
As the induction step, assume $:2k, 1,0:$ and $:2k, 1^1,0:$ are equivalent.
Figure 2.1.1.2 then outlines the legendrian  moves that demonstrate
the equivalence of
$:2k+2,1,0:$ and $:2k+2,1^1, 0:$.
\midinsert
\cl{\epsfxsize.9\hsize\epsfbox{ 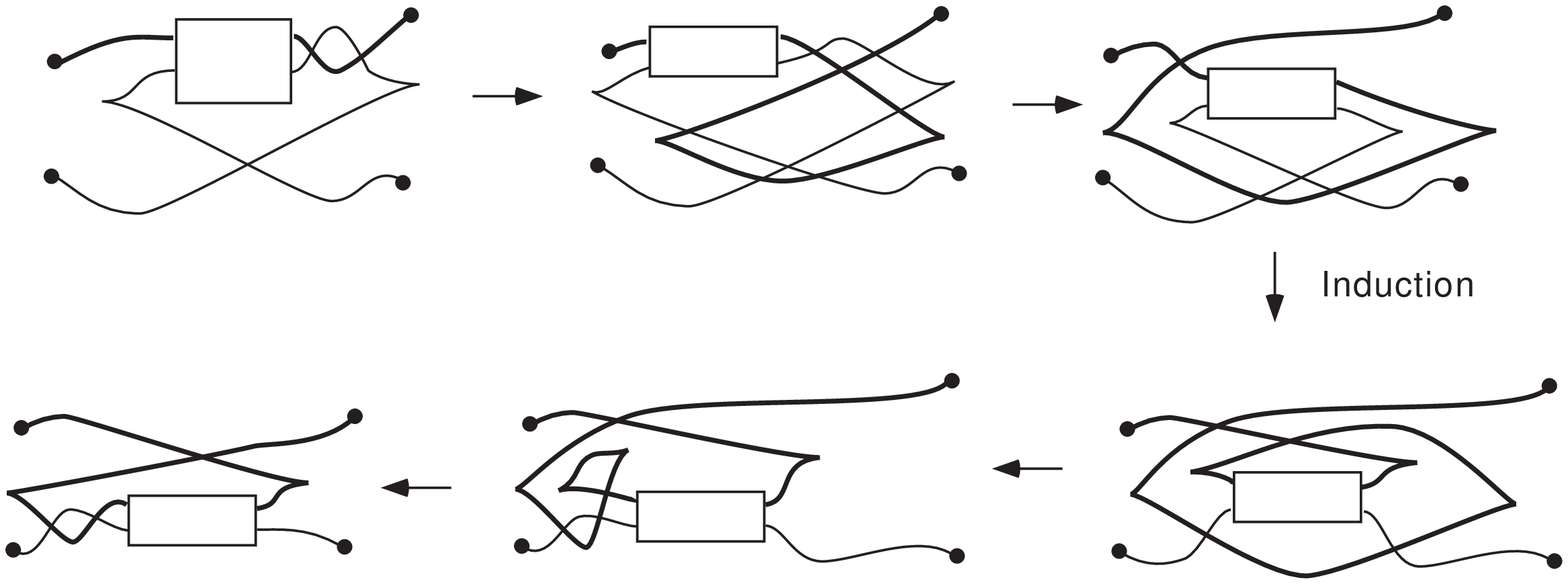}}
\lbotcaption
{Figure 2.1.1.2}  The equivalence of the tangles
$:2k+2,1,0:$ and $:2k+2,1^1, 0:$ assuming the equivalence of
$:2k,1,0:$ and $:2k,1^1, 0:$
\endcaption
\endinsert
This completes the proof of Lemma 2.1.1.  \qed
\enddemo

The induction step in the proof of Theorem 2.1
is to show that the length
$2n-1$ tangles  $:2h_n, v_{n-1}, \dots, v_2,
2h_2, v_1, 2h_1:$ and $:2h_n, v_{n-1}^{q_{n-1}}, \dots, v_2^{q_2}, 2h_2,
v_1^{q_1}, 2h_1:$ are equivalent for all choices of
$q_i$, $i = 1, \dots, n-1$,
assuming  that
the length $2n-3$ tangles
$:2h_n, v_{n-1}, \dots, v_2, 2h_2:$ and
$:2h_n, v_{n-1}^{q_{n-1}}, \dots, v_2^{q_2}, 2h_2:$ are
equivalent for all choices of $q_i$, $i = 2, \dots, n-1$.
   For this, it suffices to prove that
$:2h_n, v_{n-1}, \dots, v_2, 2h_2, 1, 0:$
and $:2h_n, v_{n-1}, \dots, v_2, 2h_2, 1^{1}, 0:$ are
equivalent.
This will   be proved by induction on $h_2$. Figure 2.1.2
outlines the moves that prove the base case  statement that
$:2h_n, v_{n-1}, \dots, v_2, 2, 1, 0:$
and $:2h_n, v_{n-1}, \dots, v_2, 2, 1^{1}, 0:$ are
equivalent when $v_2 \geq 2$.  The figure is easily
modified to prove the case $v_2 = 1$.
\midinsert
\cl{\epsfxsize.9\hsize\epsfbox{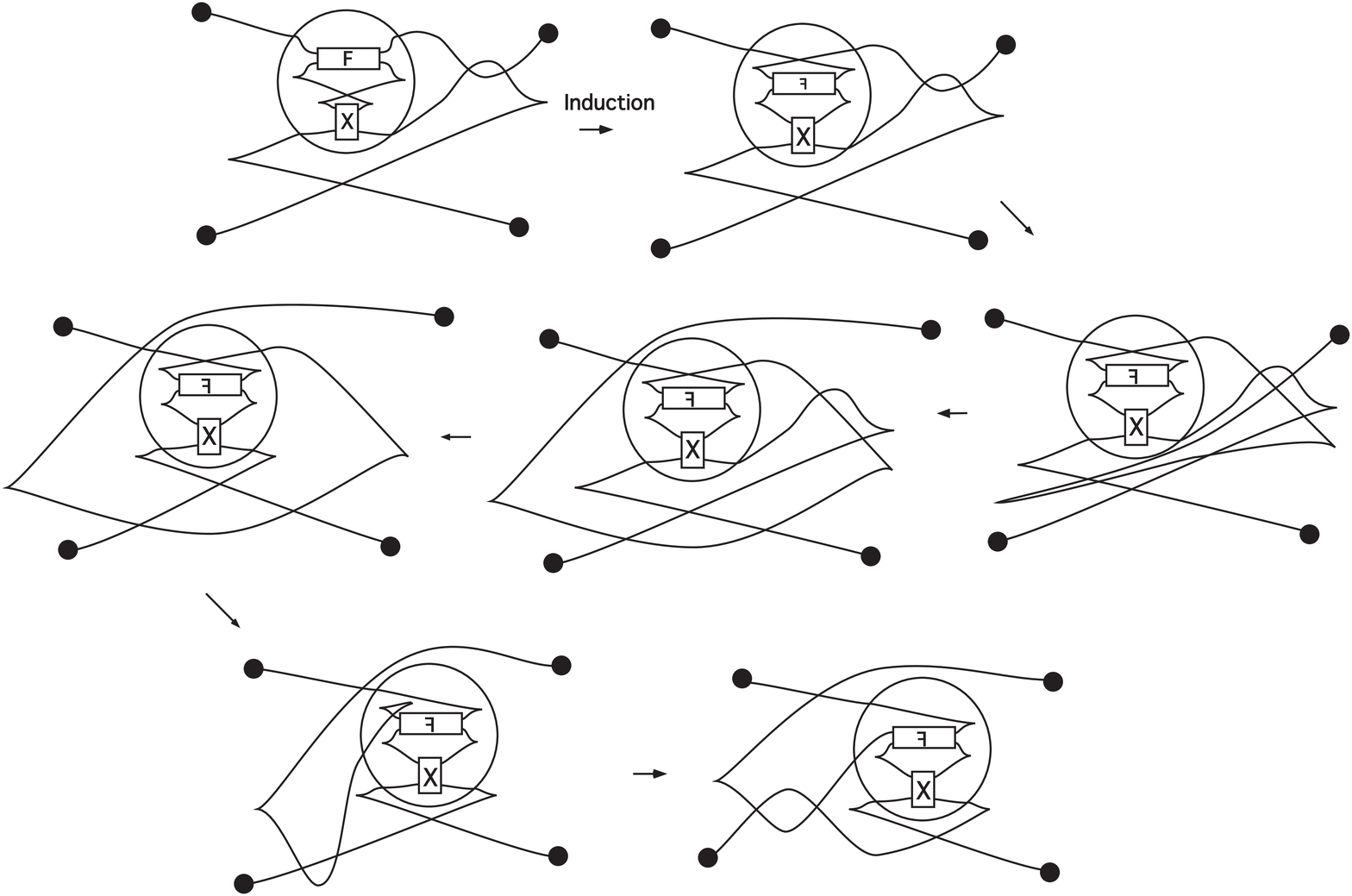}}
\lbotcaption{Figure 2.1.2} The equivalence of 
the tangles $:2h_n, v_{n-1}, \dots, v_2, 2, 1, 0:$
and\nl $:2h_n, v_{n-1}, \dots, v_2, 2, 1^{1}, 0:$
assuming the equivalence of the tangles\nl
$:2h_n, v_{n-1}, \dots, v_2, 0:$
and $:2h_n, v_{n-1}, \dots, v_2^{q_2}, 0:$  
\endcaption
\endinsert

The induction statement is that the equivalence of the tangles \newline
 $:2h_n, v_{n-1}, \dots, v_2, 2k+2, 1, 0:$
and $:2h_n, v_{n-1}, \dots, v_2, 2k+2, 1^{1}, 0:$
follows from the equivalence of the tangles
$:2h_n, v_{n-1}, \dots, v_2, 2k, 1, 0:$ and \newline
$:2h_n, v_{n-1}^{q_{n-1}}, \dots, v_2^{q_2}, 2k, 1^{1}, 0:$
This is proven using a sequence of moves similar to those shown
in Figure 2.1.1.2.  \qed
\enddemo

A nice consequence of Theorem 2.1 is that
the flyping procedure can be seen as a generalization of
the swapping procedure.

\proclaim{Corollary 2.3}  For $h_1 \geq 1$, consider the
legendrian links
$$
 L_0 = (2h_n, v_{n-1},   \dots, 2h_2, v_1, 2h_1), \quad
 L_1 = (2h_n, v_{n-1},  \dots, 2h_2, v_1, 2h_1^{1}),
$$
when $v_1, \dots, v_{n-2}$ are even.
If $\overline{L_0}$ denotes the swap of $L_0$,   then
$\overline{L_0}$ and $L_1$ are equivalent legendrian links.
\endproclaim

\demo{Proof}  By a rotation, $L_1$ is equivalent to the
swap of
$$\left(2h_n, v_{n-1}^{v_{n-1}}, 2h_{n-1}^{p_{n-1}}, v_{n-2}^{v_{n-2}}, \dots,
2h_2^{p_2}, v_1^{v_1}, 2h_1\right),$$
where, for $i = 2, \dots, n-1$,
$$
p_i = \cases
0, & \tau(i) \equiv 0 \mod 2 \\
2h_i, & \tau(i) \equiv 1 \mod 2
\endcases, \quad \text{where }\tau(i) = \sum_{k=1}^{i-1} v_k.
$$
Thus, if $v_1, \dots, v_{n-2}$ are even,
$L_1$ is equivalent to the swap of \newline
$\left(2h_n, v_{n-1}^{v_{n-1}}, 2h_{n-1}, v_{n-2}^{v_{n-2}}, \dots,
2h_2, v_2^{v_2}, 2h_1\right)$, and thus, by Theorem 2.1, $L_1$ is
equivalent to the swap of ${L_0}$.
\qed
\enddemo


\heading{Generating Function Theory}\endheading

Recall that the links in $\jetS$ under consideration are minimal, and
thus, by definition, each strand is legendrian isotopic
to $j^1(0)$, the $1$--jet of the $0$--function.  This condition guarantees
that   each
component of the link has an essentially unique ``generating function".
 The technique of generating functions is an
extension of the fact that
the $1$--jet of a smooth
function, $f\: \{ q \in S^1 \}  \to  \{ z \in \Bbb R
\}$,  is a legendrian submanifold.
More generally, if
$F \: S^1 \times \Bbb R^N \to \Bbb R$
has fiber derivatives
$\pd{F}{e}$ transverse to $0$, then
$$ \Lambda := \left\{ \left( q_0,
\pd{F}{q}(q_0,e_0),\  F(q_0,e_0) \right) \:
\frac{\partial{F}}{\partial e} (q_0,e_0) = 0
\right\}, \tag{3.1} $$
is an immersed legendrian submanifold of $\jetS$, and
$F$ is called a {\sl generating function} for
$\Lambda$. Critical points of
a generating function correspond to points where $\Lambda$
intersects   $\{ p = 0 \}$.
A function $F\: S^1 \times \Bbb R^N \to \Bbb R$ is said to
be
{\sl quadratic at infinity}
if there exists a fiberwise quadratic, nondegenerate form
$Q(q,e)$ such that
$F(q,e) - Q(q,e)$ has compact support.
The {\sl index} of a  quadratic
at infinity function  will refer
to the index of the associated quadratic form.
 The abbreviation
{\gqi function} will be used to denote a  generating and
quadratic at infinity function.   There is a
parallel definition of quadratic at infinity
generating functions for
lagrangian submanifolds of cotangent bundles; see, for
example, \cite{V1}, \cite{Th1}, \cite{V2}, \cite{Tr1}, \cite{EG}.

The following existence theorem
is proved, for a more general situation,
 by Chaperon in \cite{Cp}, by Chekanov in \cite{Ck1}, and in the
appendix to
\cite{Tr2}.

\proclaim{Existence (3.2)}
If $\Lambda_t\i \jetS$, $t \in [0,1]$, is
a smooth $1$--parameter family of legendrian submanifolds such
that $\Lambda_0 = \jet0$,
then there
exists an $N \in \Bbb N$, and a smooth 1--parameter family
  $F_t \: S^1 \times \Bbb R^N \to \Bbb R$, $t \in [0,1]$, of \gqi functions
for $\Lambda_t $.
\endproclaim

\example{Example 3.3}  For the strands of the links under consideration in this
paper,   generating functions can be
explicitly described.   For example, if the $(q,z)$--projection
of $\Lambda$, $\pi_{q,z}(\Lambda)$, is the
graph of a function
$f$, then $f \: S^1 \to \Bbb R$ is a  \gqi function for
$\Lambda$. Notice that if $Q(e)$ is a quadratic form,
$ Q(e) = \sum \alpha_{ij} e_i e_j,$
then $F(q,e) = f(q) + Q(e)$ is also a   \gqi function for $\Lambda$.
Next consider $\Lambda$ given as the ``nongraph" strand the link
$(2,1^1,0)$ as pictured in Figure 1.12.
  Construct a  \gqi function $F \: S^1 \times \Bbb R \to
\Bbb R$ for
this strand with a ``bubble" as follows.   On
$\{ q = 0 \}$, let the fiber function $F(0, \cdot) \: \Bbb R \to
\Bbb R$ be a  quadratic function of index
$1$ with critical point $a$ with value given by $\pi_{q,z}(\Lambda) \cap \{ q
= 0
\}$. For $q_0$ in a neighborhood of $0$, $F(q_0, \cdot) \:
\Bbb R \to \Bbb R$ continues to be a quadratic function, and there
is a path of critical points  $a(q_0) \in \{q_0\} \times \Bbb R$ of $F(q_0,
\cdot)$
with values given by $\pi_{q,z}(\Lambda) \cap \{ q = q_0 \}$.
At the $q$--coordinate where
  the left cusp occurs, the fiber function experiences
 a birth of a degenerate critical point.
This can be accomplished by a compact
  perturbation of the function.
 As $q$ increases, this degenerate critical
point bifurcates into  two paths of   nondegenerate
critical points $b(q_0), c(q_0) \in \{ q_0 \} \times \Bbb R$
of $F(q_0, \cdot)$ of indices $1, 0$ with
$b(q_0)$ of index $1$ having a larger critical value.
As $q$ increases further, the critical values of the critical points
$a(q_0), b(q_0), c(q_0)$
 of $F(q_0, \cdot)$
are traced out by $\pi_{q,z}(\Lambda) \cap \{ q = q_0\}$.
Eventually, the  critical value of $b(q_0)$  is larger
than the critical value of $a(q_0)$, and the
critical value of $a(q_0)$ approaches the value of the
critical value of $c(q_0)$.   At the $q$--coordinate
where the right cusp occurs, the critical points $a(q_0)$
and $c(q_0)$  merge
to form a degenerate critical point which dies as
$q$ increases further.  After this point, $F(q_0, \cdot)$ is
again a quadratic function of index $1$.  After applying
fiber preserving diffeomorphisms, $F$ will be
quadratic at infinity of index $1$. This procedure can be
generalized to construct a  \gqi function for a strand $L_2$ from
a \gqi function for a strand $L_1$, when $L_2$ has an additional bubble
resulting from a legendrian type $1^\pm$ move applied to $L_1$.

\endexample

As can be seen from the previous example, there are
 choices in the domain and in the location of
the critical points, but not in the critical values.
If $\Lambda$ is defined by a
\gqi function $F \: S^1 \times \Bbb R^N \to \Bbb R$, Theorem 3.4
shows that
all other  \gqi functions for $\Lambda$ arise
from the following \lq\lq natural modifications"
of $F$:
\roster
\item{\bf  Fiber Preserving Diffeomorphism}\qua Given a fiber preserving
diffeomorphism $\Phi\: S^1 \times \Bbb R^N \to S^1
\times \Bbb R^N$, consider
$\widetilde F = F \circ \Phi$;
\item{\bf Stabilization}\qua Let $Q \: \Bbb R^M \to \Bbb R$ be a
nondegenerate  quadratic form, $Q(f) = \sum \alpha_{ij}f_if_j$,
and consider
$\widetilde F \: S^1 \times \Bbb R^N \times \Bbb R^M \to \Bbb
R$ defined by
$\widetilde F(q,e,f) = F(q,e) + Q(f)$.
\endroster

The following uniqueness theorem parallels   a
  uniqueness result for  quadratic at
infinity generatating functions
of lagrangian submanifolds,  proved by
Viterbo and Th\'eret,  \cite{V1},\cite{Th1}, \cite{V2}.
The following can be proved by modifying Th\'eret's  careful proof to
the legendrian setting, replacing references to Sikorav's lagrangian
\gqi function existence
results
by the above mentioned legendrian \gqi function
existence results, \cite{Cp}.

\proclaim{Uniqueness Theorem 3.4}{\rm (Th\'eret)}\qua Let $\Lambda \subset \jetS$ be
legendrian isotopic to $\jet0$.
If
$F_1, F_2$ are both \gqi functions for $\Lambda$, then there exist
nondegenerate quadratic forms
$Q_1, Q_2$ and a fiber preserving diffeomorphism $\Phi$ so that
$ F_2 + Q_2 = \left( F_1 + Q_1 \right) \circ \Phi.$
\endproclaim

\definition{Definition 3.5}
Consider a minimal legendrian link $L = \Lambda_1 \amalg \Lambda_0$.
  Let
$$F_1 \: S^1 \times \Bbb R^{N_1} \to \Bbb R, \quad
F_0 \: S^1 \times \Bbb R^{N_0} \to \Bbb R$$
be \gqi functions
for $\Lambda_1, \Lambda_0$.
Then the associated (quadratic at infinity) {\sl difference function} of $L$,
$\Delta \: S^1 \times \Bbb R^{N_1} \times \Bbb R^{N_0}
\to \Bbb R$, is defined as
$$ \Delta (q, e_1, e_0) :=
F_1(q, e_1) - F_0(q, e_0).$$
\enddefinition

\proclaim{Proposition 3.6} Suppose $L = \Lambda_1 \amalg \Lambda_0  \subset
\Cal
J^1(S^1)$ is a minimal legendrian link.   Then for any  difference function
$\Delta$ of $\Lambda_1 \amalg \Lambda_0$, critical points of $\Delta$ are
in $1--1$ correspondence with points
 $((q_0, p_0, z_1), (q_0, p_0, z_0)) \in \Lambda_1 \times \Lambda_0$, and
 $0$ is never a critical value of $\Delta$.
\endproclaim

\demo{Proof}  Using formula (3.1), it is easy to verify that
the function $\Delta\: S^1
\times
\Bbb R^{N_1} \times
\Bbb R^{N_0} \to \Bbb R$ generates the
legendrian $D :=
\{\, (q, p_1 - p_0, z_1 - z_0 ) \:
( q, p_i , z_i ) \in \Lambda_i
\}.$
Since critical points of $\Delta$ correspond to points where $D$ intersects
$\{ p = 0\}$,  there is a $1--1$ correspondence between critical points of
$\Delta$
and the specified points of $\Lambda_1 \times \Lambda_0$.
Furthermore, since $\Lambda_1 \cap \Lambda_0 = \emptyset$,
$0$  cannot be a critical value for $\Delta$. \qed
\enddemo

For $c \in \Bbb R$, a noncritical value of $\Delta\: S^1 \times \Bbb
R^{N_1} \times \Bbb R^{N_0}
\to \Bbb R$,
let
$$
\Delta^{c}  :=  \{ (q, e_1, e_0) \: \Delta(q, e_1,
e_0) \leq c \}. \tag{3.7}
 $$
For every link, there exists $M > 0$ so that
all critical values of $\Delta$ are contained in
$[-M + \epsilon, M - \epsilon]$, for some $\epsilon > 0$.  For such an $M$, let
 $$
 \Delta^{\pm \infty} := \Delta^{\pm M}.
\tag{3.8} $$

\definition{Definition 3.9}
The {\sl total, positive, and negative homology groups} of a minimal
legendrian link
$L= \Lambda_1 \amalg \Lambda_0 \subset \jetS$ are defined as
$$\align
H_k(L)& := H_{k+q}\left(\Delta^\infty,
\Delta^{-\infty}\right), \\
 H_k^+(L)& :=
H_{k+q}\left(\Delta^\infty, \Delta^0\right), \\
H_k^-(L)& := H_{k+q}\left(\Delta^0,
\Delta^{-\infty}\right), \qquad k \in \Bbb Z,
\endalign
$$
where $\Delta$ is a difference function for $L$,
$q $ is the index of
$\Delta$, and the relative homology groups are calculated
with coefficients in $\Bbb Z_2$.
\enddefinition

From the definitions of these homology groups, if $\Delta$
has index $q$, then
critical points of $\Delta$ of index $i$ will often contribute
to the $H_{i-q}^\pm(L), H_{i-q}(L)$.  For this reason,
when $\Delta$ is quadratic at
infinity of index $q$,
if $x$ is a critical point of $\Delta$
with index $i$, the {\sl shifted index}
of $x$ is defined as $i - q$.

The following is a classical result of Morse theory, but
for the reader's convenience, a proof will be given.
This lemma will be useful when showing that
$H_k(L), H_k^\pm(L)$ are well-defined invariants of $L$,
and in the Section 6 calculations of these homology groups.

\proclaim{Lemma 3.10}  Consider a smooth $1$--parameter family  of
 quadratic  at infinity functions $\Delta_t \: S^1 \times \Bbb
R^N \to \Bbb R$,
$t \in [0,1]$, of index $q$.
Given paths $\alpha, \beta \: [0,1] \to \Bbb R$ such that,
for all $t$, $\alpha(t), \beta(t)$
are noncritical values of $\Delta_t$ with $\alpha(t) < \beta(t)$, then
$$ H_{q+k}\left(\Delta_0^{\beta(0)}, \Delta_0^{\alpha(0)}\right) \simeq
H_{q+k}\left(\Delta_t^{\beta(t)}, \Delta_t^{\alpha(t)}\right), \qquad \forall
t
\in [0,1], \quad
\forall k \in \Bbb Z.
$$
\endproclaim

\demo{Proof}  By applying fiber preserving diffeomorphisms, which will not
change the calculation of the homology groups, it can be assumed that for
all $t_0, t_1
\in [0,1]$,
 $\Delta_{t_0} = \Delta_{t_1}$ outside a compact set.
It suffices to show that for all $t_0 \in [0,1]$, there
exists a neighborhood $U(t_0)$ of $t_0$ such that
$H_{q+k}(\Delta_t^{\beta(t)}, \Delta_t^{\alpha(t)}) \simeq
H_{q+k}(\Delta_{t_0}^{\beta(t_0)}, \Delta_{t_0}^{\alpha(t_0)})$.

First notice that since $\alpha(t)$ and $\beta(t)$ are noncritical
values for $\Delta_t$, it is possible to choose $\epsilon > 0$ such
that for all $t \in [0,1]$, there are no critical values of $\Delta_t$ in
$(\alpha(t) - 2\epsilon, \alpha(t) + 2\epsilon)$ or $(\beta(t)-2\epsilon,
\beta(t) + 2\epsilon)$.  This implies that if $|b-\beta(t)| < 2\epsilon$
and $|a - \alpha(t)| < 2\epsilon$, then
 $H_{q+k}(\Delta_t^{b}, \Delta_t^{a}) \simeq
H_{q+k}(\Delta_{t}^{\beta(t)}, \Delta_{t}^{\alpha(t)})$.

Next, using such an $\epsilon$, choose a neighborhood $U(t_0)$ of $t_0 \in
[0,1]$ such
that, for
$t \in U(t_0)$,
\roster
\item $ \sup
\left\{ |\Delta_t(x) - \Delta_{t_0} (x) | \: x \in S^1 \times \Bbb
R^N\right\} < \epsilon$, and
\item $|\beta(t) - \beta(t_0)| < \epsilon$, and $|\alpha(t) - \alpha(t_0)| <
\epsilon.$
\endroster
By (2), $H_{q+k}\left(\Delta_{t}^{\beta(t)}, \Delta_{t}^{\alpha(t)}\right)
\simeq
H_{q+k}\left(\Delta_{t}^{\beta(t_0) + \epsilon}, \Delta_{t}^{\alpha(t_0) +
\epsilon}\right)$.  By (1), for  $c = \beta(t_0)$ and $c  = \alpha(t_0)$, the
inclusions
$\Delta_{t}^{c-\epsilon} \i \Delta_{t_0}^c \i \Delta_{t}^{c+\epsilon}
\i \Delta_{t_0}^{c + 2 \epsilon}$
 induce homomorphisms
$$\align
 H_{q+k} \left( \Delta_{t}^{\beta(t_0)-\epsilon},
\Delta_{t}^{\alpha(t_0)-\epsilon} \right)
&\buildrel{\phi_1} \over \longrightarrow
H_{q+k} \left( \Delta_{t_0}^{\beta(t_0)}, \Delta_{t_0}^{\alpha(t_0)} \right)
\buildrel{\phi_2} \over \longrightarrow \\
 H_{q+k} \left( \Delta_{t}^{\beta(t_0)+\epsilon},
\Delta_{t}^{\alpha(t_0)+\epsilon} \right)
 &\buildrel{\phi_3} \over \longrightarrow
H_{q+k} \left( \Delta_{t_0}^{\beta(t_0) + 2\epsilon},
\Delta_{t_0}^{\alpha(t_0)+2\epsilon} \right).
\endalign $$
Since $|\beta(t_0) \pm \epsilon - \beta(t)|,
|\alpha(t_0) \pm \epsilon - \alpha(t)|, |\beta(t_0) \pm 2\epsilon -
\beta(t_0)| <
2\epsilon$,  the first and third groups are
isomorphic to $H_{q+k}\left( \Delta_{t}^{\beta(t)}, \Delta_{t}^{\alpha(t)}
\right)$, the second and fourth groups are isomorphic to
 $H_{q+k}\left( \Delta_{t_0}^{\beta(t_0)}, \Delta_{t_0}^{\alpha(t_0)}
\right)$, and
 $\phi_2 \circ \phi_1$  and
$\phi_3 \circ \phi_2 $ are isomorphisms.  Thus $\phi_2$ is an
isomorphism, and  it follows that
$H_{q+k}\left(\Delta_t^{\beta(t)}, \Delta_t^{\alpha(t)}\right) \simeq
H_{q+k}\left(\Delta_{t_0}^{\beta(t_0)}, \Delta_{t_0}^{\alpha(t_0)}\right)$.
\qed
\enddemo

\proclaim{Theorem 3.11}  $H_k(L)$, $H_k^+ (L)$, and
$H_k^- (L)$ are well-defined invariants of a minimal
legendrian link $L \i \jetS$.
\endproclaim

\demo{Proof} It must be shown that the homology groups
do not depend on the choice of generating functions $F_i$
for the strands,
and will not change as the link  undergoes a legendrian
isotopy.

Suppose $L = \Lambda_1 \amalg \Lambda_0$, and let $F_i \: S^1 \times \Bbb
R^{N_i} \to
\Bbb R$  be  \gqi function
for $\Lambda_i$.
It must be shown that if
$\widetilde F_i \: S^1 \times \Bbb R^{M_i} \to \Bbb R$
are other \gqi function for $\Lambda_i$,
  then the relative
homology groups of $\widetilde \Delta(q,\widetilde{e_1},
\widetilde{e_0}) := \widetilde F_1(q,\widetilde{e_1}) -
\widetilde F_0(q, \widetilde{f_0})$ agree, up to a shift
of appropriate index, with those of
$\Delta(q,e_1, e_0) := F_1(q,e_1, e_0) - F_0(q,e_1, e_0)$.
By the Uniqueness Theorem 3.4,  it is only necessary to check the cases where
$\widetilde \Delta$ differs from $\Delta$ by a fiber preserving diffeomorphism
or by  a stabilization.  If $\widetilde \Delta = \Delta \circ \Phi$,
then the associated sublevel sets are diffeomorphic:
$ \widetilde \Delta^c = \Phi^{-1}(\Delta^c)$.  Thus
the relative homology groups are unchanged.
Next, suppose that  $\widetilde \Delta \: S^1 \times \Bbb R^{N}
\times \Bbb R^{M}  \to \Bbb R$ is
defined by
$\widetilde \Delta (q,e,f) = \Delta(q,e) + Q(f)$,
where $\Delta$ is a \gqi function of index $q$, and $Q$ is a
nondegenerate quadratic form of index $j$.  It is easily
verified that the critical values of $\Delta$ agree with
those of $\tilde \Delta$, and that for a noncritical value
$v$, for all $k \in \Bbb Z$, there is an isomorphism $\phi\:
C_{q+k}^v(\Delta) \to
 C_{q+ j + k}^v(\widetilde \Delta)$, where
$C_i^v(\Delta)$ (respectively $C_i^v(\widetilde \Delta)$) denotes the
$i$--chains of $\Delta^v$ (respectively $\widetilde \Delta^v$).
For all noncritical values $a < b$, and all $k \in \Bbb Z$,
 the isomorphism $\phi$
induces chain isomorphisms
$\widetilde \phi \: C_{q+k}^b(\Delta)/C_{q+k}^a(\Delta) \to
C_{q+j+k}^b(\widetilde \Delta)/C_{q+j+k}^a(\widetilde\Delta),$
and thus $H_{q+k}(\Delta^b, \Delta^a) \simeq
H_{q+j+k}(\widetilde\Delta^b, \widetilde\Delta^a)$, as desired.

Suppose $L_t$, $ t \in [0,1]$, is a $1$--parameter family of minimal
legendrian links.  By Existence (3.2), there exist difference
functions
$\Delta_t \: S^1 \times \Bbb R^N \to \Bbb R$ of index $q$  for
$L_t$.  It must be shown
that, for all $t \in [0,1]$, and all $k \in \Bbb Z$,
$$\align
H_{k+q} (\Delta_t^{+\infty}, \Delta_t^{-\infty}) &\simeq
H_{k+q} (\Delta_0^{+\infty}, \Delta_0^{-\infty}),\\
H_{k+q} (\Delta_t^{+\infty}, \Delta_t^{ 0}) &\simeq
H_{k+q} (\Delta_0^{+\infty}, \Delta_0^{ 0}),\\
H_{k+q} (\Delta_t^{ 0}, \Delta_t^{-\infty}) &\simeq
H_{k+q} (\Delta_0^{ 0}, \Delta_0^{-\infty}).
\endalign$$
Choose paths $\alpha, \beta, \gamma \: [0,1] \to \Bbb R$
such that $\alpha(t)$ is negative and is less than all critical values
of $\Delta_t$, $\beta(t) = 0$, and $\gamma(t)$ is positive and is greater
than all critical values of $\Delta_t$. By construction and
Proposition 3.6, $\alpha(t), \beta(t), \gamma(t)$ are noncritical
values of $\Delta_t$, and thus the desired result holds
by Lemma 3.10.  \qed
\enddemo

\proclaim{Proposition 3.12} For any minimal legendrian link $L$,
 $H_k(L) \simeq H_k(S^1)$,  for all $k \in \Bbb Z$.
\endproclaim

\demo{Proof} By hypothesis, each strand of $L$ can be individually
isotoped so that it is the graph of a function.  Thus there exists
a $1$--parameter family of immersed legendrian links $L_t$, $t \in [0,1]$ with
$L_0 = L$, $L_1 = j^1(+f) \amalg j^1(-f)$, where $j^1(\pm f)$ is the
$1$--jet of $\pm f(q) = \pm(\cos(2\pi q) + 2)$. For the associated
 $1$--parameter family of difference functions $\Delta_t$ of
$L_t$, choose paths
$\alpha, \beta\: [0,1] \to \Bbb R$
such that  $\alpha(t)$ is negative and less than all critical values
of $\Delta_t$, and $\beta(t)$ is positive and greater than all critical values
of $\Delta_t$. Since it is easy to understand the total homology group of
$L_1$,  the desired result then follows by Lemma 3.10.
\qed
\enddemo

\proclaim{Lemma 3.13}  For a  function $\Delta \: S^1 \times \Bbb R^N \to
\Bbb R$,
and $a,b,c$ noncritical values of $\Delta$ with $a < b < c$, there is a
long exact sequence\smallskip
$\dots
\buildrel \partial_* \over \to H_{q+k}(\Delta^b, \Delta^a)
\buildrel i_* \over \to
 H_{q+k}( \Delta^c,\Delta^a)
\buildrel \pi_* \over \to$\nl
\hbox to 10pt{\hss}\hfill$H_{q+k}( \Delta^c, \Delta^b)
\buildrel \partial_* \over \to
 H_{q+ k-1}(\Delta^{b}, \Delta^{a} )
\buildrel i_* \over \to \dots.
$
\endproclaim

\demo{Proof}
 Given $\Delta \: S^1 \times \Bbb R^N \to \Bbb R$ of index $q$,
and a noncritical value
$v$ of
$\Delta$,
 $C_k^v(\Delta)$ will denote the $(q+k)$--chains of $\Delta^v$.
Given a triple
$ a < b < c $ of noncritical values of $\Delta$,
 for each $k$, there is a natural exact sequence
 $$ 0 \to  C_k^b(\Delta) / C_k^a(\Delta) \buildrel i \over \to
C_k^c(\Delta) / C_k^a(\Delta)
\buildrel \pi \over \to   C_k^c(\Delta) / C_k^b(\Delta) \to 0.
$$
Since $i, \pi$ are chain maps, it
follows that
   there is an exact sequence
\smallskip
$\dots
\buildrel \partial_* \over \to H_{q+k}(\Delta^b, \Delta^a)
\buildrel i_* \over \to
 H_{q+k}( \Delta^c,\Delta^a)
\buildrel \pi_* \over \to$\nl
\hbox to 100pt{\hss}\hfill$H_{q+k}( \Delta^c, \Delta^b)
\buildrel \partial_* \over \to
 H_{q+ k-1}(\Delta^{b}, \Delta^{a} )
\buildrel i_* \over \to \dots.
$\qed
\enddemo

 \proclaim{Corollary  3.14} For a minimal legendrian
link $L$,
there is an exact
sequence
$$ \dots
\buildrel \partial_* \over \to H_k^- (L)
\buildrel i_* \over \to H_k(L)
\buildrel \pi_* \over \to H_k^+(L)
\buildrel \partial_* \over \to H_{k-1}^-(L)
\buildrel i_* \over \to \dots .$$
\endproclaim

\definition{Definition 3.15}  For a minimal legendrian
link $L$,
form the {\sl positive} and {\sl negative homology polynomials}:
$$\Gamma^+(\lambda)[L] = \sum_{k= -\infty}^\infty \dim H_k^+(L) \cdot
\lambda^k, \qquad
\Gamma^-(\lambda)[L] = \sum_{k= -\infty}^\infty \dim H_k^-(L) \cdot
\lambda^k.$$
\enddefinition

A comparison of $\Gamma^+(\lambda)[L]$ and $\Gamma^-(\lambda)[L]$  may
detect that the legendrian link $L$ is ordered.
Recall from Section 1 that   polynomials
$\alpha(\lambda) = \sum_{k=-\infty}^\infty a_k \lambda^k$ and
$\beta(\lambda) = \sum_{k=-\infty}^\infty b_k \lambda^k$ are 1--shift
palindromic if
$ \alpha(\lambda) = \lambda \cdot \overline{\beta(\lambda)},$
where $\overline{\beta(\lambda)}$ denotes the palindrome of $\beta(\lambda)$.

\proclaim{Theorem 3.16} Let $L = \Lambda_1 \amalg \Lambda_0$ be a  minimal
legendrian
link, and let $\overline L$ denote its swap,
$\overline L = \Lambda_0 \amalg \Lambda_1$.
Then
$\Gamma^+(\lambda)[L]$ and $\Gamma^-(\lambda)[\overline L]$
are $1$--shift palindromic, and $\Gamma^-(\lambda)[L]$ and
$\Gamma^+(\lambda)[\overline L]$ are $1$--shift
palindromic.
\endproclaim

\demo{Proof} If $\Delta\: S^1 \times
\Bbb R^{N-1} \to \Bbb R$ is a  difference function
for  $L$, then $\overline \Delta
= -\Delta$ is a difference function for $\overline L$. If $Q$ and $\overline
Q$ denote the indices  of
$\Delta$ and $\overline \Delta$, then $\overline Q = N-1-Q$.
Notice
$$\align
H_k^+(L) &= H_{k+Q}(\Delta^{+\infty}, \Delta^0)
\simeq H^{N-(k+Q)}(\Delta^{+\infty}, \Delta^0)
\simeq  H_{N-(k+Q)}\left( \overline{\Delta}^0, \overline{\Delta}^{\ -\infty}
\right)\\
&= H_{N-(k+Q)-\overline{Q}}^-(\overline L)
= H_{N-(k+Q)-N+1+Q}^-(\overline L)
= H_{-k+1}^-(\overline L).
\endalign$$
Thus $\dim H_k^+(L) = \dim H_{-k+1}^-(\overline L)$, and it follows that
$\Gamma^+(\lambda)[L]$ and $\Gamma^-(\lambda)[\overline L]$ are $1$--shift
palindromic.\kern 0.5pt
A similar calculation shows that $\Gamma^-(\lambda)[L]$ and
$\Gamma^+(\lambda)[\overline L]$ are $1$--shift
palindromic. \qed
\enddemo

\proclaim{Corollary 3.17}  Let $L = \Lambda_1 \amalg \Lambda_0$ be a
   minimal legendrian
link.  If
$\Gamma^+(\lambda)[L]$ and
$\Gamma^-(\lambda)[L]$ are not $1$--shift palindromic, then the link $L$ is
ordered.
\endproclaim

\demo{Proof}   If $L$ is not ordered, then $\Gamma^+(\lambda)[L] =
\Gamma^+(\lambda)[\overline L]$, and thus $\Gamma^+(\lambda)[L]$ and
$\Gamma^-(\lambda)[L]$ must be
$1$--shift palindromic.  \qed
\enddemo


\heading{Algebraic Topology Tools}\endheading

In this section, some tools will be developed that will
aid in the Section 6 calculations  of the link polynomials of
rational links and  connect sums of rational
links.  The first result is called the
Additive Extension Lemma since it gives conditions when,
for $b_2 > b_1 > 0$ and $a_2 < a_1 < 0$,
$\dim H_{*}(\Delta^{b_2}, \Delta^0) =
\dim H_{*}(\Delta^{b_2}, \Delta^{b_1}) +
\dim\! H_{*}(\Delta^{b_1}, \Delta^0),$
and
$\dim\! H_{*}(\Delta^{0}, \Delta^{a_2}) =
\dim\! H_{*}(\Delta^{0}, \Delta^{a_1}) + \dim\!
H_{*}(\Delta^{a_1},
\Delta^{a_2})
$.

\proclaim{Additive Extension Lemma 4.1}  Suppose that
$c_1^\pm, c_2^\pm$ are critical values of a
quadratic at infinity function
$\Delta \: S^1 \times \Bbb R^N \to \Bbb R$ of  index $q$, and that
$a_2,a_1, b_1,b_2$ are noncritical values such that
$$a_2 < c_2^-< a_1 < c_1^- < 0 < c_1^+< b_1 < c_2^+ < b_2.$$
Suppose that
\roster
\item $c_2^+$ is the only critical value in $[b_1, b_2]$,
$c_2^-$ is the only critical value in $[a_2, a_1]$,
  all critical points with  values $c_2^+$
are nondegenerate and have the same index, and all critical points with value
$c_2^-$ are nondegenerate and have the same index,
\item $H_{q+k}(\Delta^{b_1}, \Delta^{a_1}) = 0, \quad \forall
k \in \Bbb Z$, and
\item $H_{q+k}(\Delta^0, \Delta^{a_1}) \simeq
H_{q+k+1}(\Delta^{b_1},
\Delta^0), \quad \forall k \in \Bbb Z$.
\endroster
Then
$$\align
\dim H_{q+k}(\Delta^{b_2}, \Delta^0) &=
\dim H_{q+k}(\Delta^{b_2}, \Delta^{b_1}) +
\dim H_{q+k}(\Delta^{b_1}, \Delta^0), \text{ and } \\
\dim H_{q+k}(\Delta^{0}, \Delta^{a_2}) &=
\dim H_{q+k}(\Delta^{0}, \Delta^{a_1}) + \dim
H_{q+k}(\Delta^{a_1},
\Delta^{a_2}), \quad \forall k\in \Bbb Z.
\endalign$$
\endproclaim

\demo{Proof} The statement about
$\dim H_{q+k}(\Delta^{b_2}, \Delta^0)$ follows
immediately if $\partial_* \equiv 0$ in the
long exact sequence
$$ \aligned
\dots
&\buildrel \partial_* \over \to
H_{q+k}(\Delta^{b_1},\Delta^0)
\buildrel i_* \over \to\\
 H_{q+k}( \Delta^{b_2},\Delta^0)
\buildrel \pi_* \over \to
H_{q+k}( \Delta^{b_2}, \Delta^{b_1})
&\buildrel \partial_* \over \to
 H_{q+ k-1}(\Delta^{b_1}, \Delta^{0} )
\buildrel i_* \over \to \dots.
\endaligned
\tag{4.1.1}
$$
If there
exists $k$ such that $0 \neq \operatorname{im}
\partial_*(H_{q+k+1}(\Delta^{b_2}, \Delta^{b_1}))$,
then by hypothesis (1), $H_
{q+k}(\Delta^{b_2}, \Delta^{b_1}) = 0$, and
thus
$$\dim H_{q+k}(\Delta^{b_1},\Delta^0) >
\dim H_{q+k}(\Delta^{b_2}, \Delta^0).\tag{4.1.2} $$
Since all critical points with values in $[b_1,
b_2]$  have index $q +k + 1$,
another
exact sequence argument, using the fact that
$H_*(\Delta^{b_1},
\Delta^{a_1}) = 0$ for all $*$,  proves
$$
 H_{q+k-1}(\Delta^{b_2},\Delta^{a_1}) = 0
.\tag{4.1.3} $$
The facts $4.1.2$ and $4.1.3$  combine with the hypotheses
$$\dim H_{q+k-1}(\Delta^0, \Delta^{a_1}) = \dim
H_{q+k}(\Delta^{b_1},\Delta^0)$$  to give a
contradiction  to the necessary surjectivity of $\partial_*$
in the exact  sequence
 $$ \aligned
\dots
\buildrel \pi_* \over \to
H_{q+k}( \Delta^{b_2}, \Delta^{0})
\buildrel \partial_* \over \to
 H_{q+ k-1}(\Delta^{0}, \Delta^{a_1} )
\buildrel i_* \over \to
& H_{q+k-1}(\Delta^{b_2},\Delta^{a_1})
 \buildrel \pi_* \over \to  \dots .
\endaligned
\tag{4.1.4} $$
This proves the claim about $\dim H_{q+k}(\Delta^{b_2},
\Delta^0)$. A proof of the claim about $\dim H_{q+k}(\Delta^{0},
\Delta^{a_2})$ can be proved similarly. \qed
\enddemo

The following proposition will be useful when calculating the
homology polynomials of legendrian links that are topologically
the  rational links $q$, for $q \geq 1$.

\proclaim{Positive Integral Proposition 4.2}  Suppose the minimal
legendrian link
$L=\Lambda_1 \amalg \Lambda_0$ has a
 difference function $\Delta \: S^1 \times \Bbb R^N \to \Bbb R$ with
critical values $c_1^\pm, \dots, c_n^\pm$, and noncritical
values $a_1, \dots, a_n, b_n, \dots, b_1$ satisfying
$$ \align
&a_1 < c_1^- < a_2 < c_2^- < \dots < a_{n-1} <  c_{n-1}^- < a_n < c_n^- < 0, \\
&0 < c_n^+ < b_n < c_{n-1}^+ < b_{n-1} < \dots < c_2^+ < b_2 < c_1^+ < b_1.
\endalign$$
If
\roster
\item for $k = 1, \dots, n$, there are $h_k$ nondegenerate critical points
with value $c_k^+$,
and $h_k$ nondegenerate critical points with value $c_k^-$; all critical
points with value $c_k^+$
have shifted index $i_k+1$, and all critical points with value $c_k^-$ have
shifted index $i_k$;
and
\item for $k = 2,\dots, n$, $H_*(\Delta^{b_k}, \Delta^{a_k}) = 0$, for all
$* \in \Bbb
Z$,
\endroster
then
$$
\Gamma^+(\lambda)[L] = \sum_{k=1}^n h_k \lambda^{i_k+1}, \qquad
\Gamma^-(\lambda)[L] = \sum_{k=1}^n h_k \lambda^{i_k}.
$$
\endproclaim

\demo{Proof}  By hypothesis $(1)$,
$$\gather
\dim H_{q+*}(\Delta^{0}, \Delta^{a_n}) = \dim H_{q+*+1}(\Delta^{b_n},
\Delta^{0}) =
\cases
h_n & * = i_n\\
0 & \text{else}
\endcases, \quad\text{ and }\\
\dim H_{q+*}(\Delta^{a_{k+1}}, \Delta^{a_{k}}) =
\dim H_{q+*+1}(\Delta^{b_{k}}, \Delta^{b_{k+1}}) =
\cases
h_k & * = i_k\\
0 & \text{else}
\endcases,
 \endgather
$$
for $k = 1, \dots, n-1$.
Thus it suffices to show that for $n-1 \geq k \geq 1$,
$$\align
 \dim H_{q+*} (\Delta^{b_k}, \Delta^0) &=
\dim H_{q+*} (\Delta^{b_{k+1}}, \Delta^0) +
\dim H_{q+*}(\Delta^{b_{k}}, \Delta^{b_{k+1}}), \text{ and }\\
\dim H_{q+*} (\Delta^{0}, \Delta^{a_k}) &=
\dim H_{q+*} (\Delta^{0}, \Delta^{a_{k+1}}) +
\dim H_{q+*}(\Delta^{a_{k+1}}, \Delta^{a_{k}}), \ \forall * \in \Bbb Z.
\endalign$$
This will be proven by repeatedly applying the Additive Extension Lemma 4.1.
To apply this lemma, it must be shown that for $n \geq k \geq 2$,
$$H_{q+*}(\Delta^0, \Delta^{a_k})  \simeq H_{q+*+1}(\Delta^{b_k},
\Delta^{0}).$$
As mentioned above, this is true when $k = n$.  Assume it is true when
$k = \ell$.  Then  the Additive Extension  Lemma  applied to
$$ a_{\ell-1} < c_{\ell-1}^- < a_\ell < c_\ell^- < 0 < c_\ell^+ < b_\ell <
c_{\ell-1}^+ < b_{\ell-1}$$
shows that
$$\align
\dim H_{q+*+1}(\Delta^{b_{\ell-1}}, \Delta^0) &=
\dim H_{q+*+1}(\Delta^{b_{\ell}}, \Delta^0) + \dim
H_{q+*+1}(\Delta^{b_{\ell-1}}, \Delta^{b_\ell})\\
&=
\dim H_{q+*}(\Delta^{0}, \Delta^{a_\ell}) + \dim H_{q+*}(\Delta^{a_\ell},
\Delta^{a_{\ell-1}})\\
&=\dim H_{q+*}(\Delta^{0}, \Delta^{a_{\ell-1}}).
\endalign$$
Thus the isomorphism holds when $k = \ell - 1$.  This completes the proof.
\qed
\enddemo

The next  proposition will be used  to calculate the
homology polynomials of minimal legendrian links that are topologically
the rational links $q$, for $q < 1 $.

\proclaim{Zero Integral Proposition 4.3}  Suppose  the minimal legendrian link
$L=\Lambda_1 \amalg \Lambda_0$ has  a  difference function
$\Delta\: S^1 \times \Bbb R^N \to \Bbb R$  with critical values $c_1^0,
c_1^1, c_2^\pm,$\break $
\dots, c_n^\pm$ and noncritical values $a_2, \dots, a_n, b_n, \dots, b_1$
satisfying
$$ \align
&a_2 < c_2^- < \dots < a_{n-1} <  c_{n-1} < a_n < c_n^- < 0 \\
&0 < c_n^+ < b_n < c_{n-1}^+ < b_{n-1} < \dots < c_2^+ < b_2 < c_1^0 <
c_1^1 < b_1.
\endalign$$
If, for $k = 2, \dots, n$,
\roster
\item  there are $h_k$ nondegenerate  critical
points with value
$c_k^+$,  and $h_k$ nondegenerate critical points with value $c_k^-$;
all critical points with value $c_k^+$ have shifted index $i_k+1$, while
all critical points with value $c_k^-$ have shifted index
$i_k$, and
\item  $H_*(\Delta^{b_k}, \Delta^{a_k}) = 0$, for all $* \in \Bbb Z$,
\endroster
then
$$
\Gamma^+(\lambda)[L] = (1 + \lambda) + \sum_{k=2}^n h_k \lambda^{i_k+1}, \qquad
\Gamma^-(\lambda)[L]
= \sum_{k=2}^n h_k \lambda^{i_k}.
$$
\endproclaim

\demo{Proof}  Using the hypothesis $H_*(\Delta^{b_k}, \Delta^{a_k}) = 0$ for
$k = 3, \dots, n$,  arguments as in the proof
of Proposition 4.2 prove that
$$\align
\Gamma^-(\lambda)[L] &= \dim H_{q+ i_n}(\Delta^0, \Delta^{a_n}) \lambda^{i_n}+
\dim H_{q+ i_{n-1}}(\Delta^{a_n}, \Delta^{a_{n-1}}) \lambda^{i_{n-1}} +
\dots \\
& \text{\hskip 2.5in} + \dim H_{q + i_2}(\Delta^{a_3},
\Delta^{a_2})\lambda^{i_{2}} \\
&= \sum_{k=2}^n h_k\lambda^{i_k}.
\endalign$$
To prove the claim about $\Gamma^+(\lambda)[L]$, it suffices to show that
$$\dim H_k^+(L) = \cases
\dim H_{k-1}^-(L) + 1, & k = 0,1\\
\dim H_{k-1}^-(L), & \text{ else}.
\endcases$$
By Proposition 3.12, $\dim H_k(L) =1$, when $k = 0,1$, and vanishes
otherwise.  Thus
the desired calculations of $H_k^+(L)$ will follow if it is shown
that $i_* = 0$ in the exact
sequence
$$ \dots
\buildrel \partial_* \over \to H_k^- (L)
\buildrel i_* \over \to H_k(L)
\buildrel \pi_* \over \to H_k^+(L)
\buildrel \partial_* \over \to H_{k-1}^-(L)
\buildrel i_* \over \to \dots .$$
The map $i_*$ is induced by the inclusion  map
$ i \: \frac{C_k^0(\Delta)}{C_k^{a_2}(\Delta)} \to
\frac{C_k^{b_1}(\Delta)}{C_k^{a_2}(\Delta)}.
$
Since $i = i_2 \circ i_1$, where
$$\frac{C_k^0(\Delta)}{C_k^{a_2}(\Delta)}
\buildrel{i_1}\over \to
\frac{C_k^{b_2}(\Delta)}{C_k^{a_2}(\Delta)}
\buildrel{i_2}\over \to
\frac{C_k^{b_1}(\Delta)}{C_k^{a_2}(\Delta)},
$$
$i_* = (i_2)_* \circ (i_1)_*$.  However, since
$H_*(\Delta^{b_2}, \Delta^{a_2}) = 0$,  $(i_1)_* = 0$, and
thus $i_* = 0$, as desired.  \qed
\enddemo


\heading{Index Calculations}\endheading

Propositions 4.2 and 4.3 will make it  easy to calculate
the homology polynomials of a minimal
legendrian link $L = \Lambda_1 \amalg \Lambda_0$ from a difference function
that is in
a ``nice" form.
To apply the propositions, it will be necessary to find
the critical points of $\Delta$, and calculate
the indices of these critical points. Critical points
of $\Delta$ are easily described in terms of $\Lambda_1 \amalg \Lambda_0$:
as was shown in
Proposition 3.6, they correspond to points
$\left((q_0, p_0, z_1), (q_0, p_0, z_0)\right) \in \Lambda_1 \times \Lambda_0$.
The indices of the critical points of $\Delta$ can easily
be calculated in terms of data from the
$(q,z)$--projection of the link: as will be shown in Proposition
5.5, the index of the critical point corresponding
to $\left((q_0, p_0, z_1), (q_0, p_0, z_0)\right)$
 can be calculated as a difference
of ``branch indices" of $(q_0, p_0, z_1)$ and of $(q_0, p_0, z_0)$, and
a Morse ``graph index" between $(q_0, p_0, z_1)$ and  $(q_0, p_0, z_0)$.

Throughout this section, $\pi_{q,z}, \pi_{q,p} \: \jetS \to S^1 \times \Bbb R$
will denote the projections to the $(q,z)$--, $(q,p)$--coordinates, respectively.

\definition{Definition 5.1}  Given a legendrian knot $\Lambda
\subset
\Cal J^1(S^1)$, the {\sl branches} of $\Lambda$ are the
connected components of $\Lambda \backslash C$, where $C$ denotes the set of
points that project to cusp points in
$\pi_{q,z}(\Lambda)$.  Branches $V_1, V_0$ of $\Lambda$ are {\sl
adjacent} if their closures $\overline V_i$ intersect.
If
$V_1, V_0$ are adjacent, $V_1 \bold> V_0$ if there exists $v
\in
\overline V_1 \cap \overline V_0 $ and
a path $\gamma \subset \pi_{q,z}(\Lambda)$ with $\gamma[0,1/2) \subset
\pi_{q,z}(V_0)$,
$\gamma(1/2) = \pi_{q,z}(v)$, and $\gamma(1/2, 1] \subset \pi_{q,z}(V_1)$,
 so that $\pi_{q,z}(v)$ is an up cusp along the path.
\enddefinition

\definition{Definition 5.2}
Suppose $\Lambda \subset \Cal J^1(S^1)$
is legendrian isotopic to  $\jet0$.
Let $\{ V_i\}$ denote the
set of branches of $\Lambda$. Suppose there exists a
branch $\Cal I \in \{ V_i \}$,
 $v \in \Cal I$, and  a contact isotopy $\kappa_t$,
$t
\in [0,1]$, so that
$\kappa_1(\Lambda) = \jet0$, and
$\pi_{q,z}(\kappa_t(v))$ is never a cusp.  Then the {\sl branch index}
$i_{B} \: \{ V_i \} \to \Bbb Z$ is defined by
$i_{B}(V_i) = 0$, if
$V_i = \Cal I$, and  $i_{B}(V_1) - i_B(V_0) =
1$,  if $V_1, V_0$ are adjacent with
$V_1 > V_0$. 
\enddefinition

From this definition, it appears that the above branch index
may depend on the choice of an initial branch $\Cal I$.  However,
the next proposition shows that this is not the case.

\proclaim{Proposition 5.3}  Suppose $\Lambda$ is legendrian
isotopic to $\jet0$, $V$ is a branch of $\Lambda$, and
$v \in V$.  If $F \: S^1 \times \Bbb R^N \to \Bbb R$ is a \gqi function
of index $Q$ for $\Lambda$, and $(q_v, e_v) \in S^1 \times \Bbb R^N$
corresponds to $v$, then
$ i_B(V) = \ind F(q_v, \cdot) (e_v) - Q.$
\endproclaim

\demo{Proof}  Let $\Lambda_t$ be a legendrian isotopy
with $\Lambda_0 = \jet0$, $\Lambda_1 = \Lambda$, and
let $\pi_{q,p}(\Lambda_t)$ denote the lagrangian
projections.  By a classical result,
see for example Appendix~B of \cite{Th2} or \cite{Tr1}, the shifted index of
the fiber function  $F(q_v, \cdot) (e_v) - Q$ is equal to the
Maslov index of a path $\gamma(t) \in \pi_{q,p}(\Lambda_t)$,
$t \in [0,1]$,
with $\gamma(1) = \pi_{q,p}(v)$.  By a morsification
procedure, see for example \cite{F}, the rotation calculations
necessary to calculate the Maslov index can be computed
in terms of cusps:  by an appropriate choice of isotopy,
the Maslov index equals the difference
between the number of up and down cusps along a path
$\pi_{q,z}(\Lambda_1)$ that starts at $\pi_{q,z}(\iota)$,
$\iota \in \Cal I$, and ends at $\pi_{q,z}(v)$.
\qed
\enddemo

\definition{Definition 5.4}  Given legendrian knots
$\Lambda_1, \Lambda_0 \subset \Cal J^1(S^1)$, suppose
$(q_0, p_0, z_1)$\break$ \in V_1$, $(q_0, p_0, z_0) \in V_0$, where
$V_1, V_0$ are branches of $\Lambda_1, \Lambda_0$, respectively.
Then there exist
functions $f_1, f_0 \: U \to \Bbb R$ such
that near $(q_0, p_0, z_1)$,
 $\pi_{q,z}(V_1) = \left\{ \left(q, f_1(q) \right)\right\}$,
and near $(q_0, p_0, z_0)$,
  $\pi_{q,z}(V_0) = \{ (q, f_0(q)) \} $.
It follows that $q_0$ is a critical point of $f_1 - f_0$.
Then $((q_0,p_0, z_1), (q_0,p_0, z_1)) \in \Lambda_1 \times \Lambda_0$
is {\sl nondegenerate} if $q_0$ is a nondegenerate critical
point of $f_1 - f_0$;  the {\sl graph index}, $i_\Gamma$,  of
a nondegenerate point is the Morse index of $f_1 - f_0$ at
$q_0$.
\enddefinition

\proclaim{Proposition 5.5}  Suppose $L = \Lambda_1 \amalg
\Lambda_0
\subset
\Cal J^1(S^1)$ is a minimal legendrian link.
Given a nondegenerate
 $((q_0, p_0, z_1), (q_0, p_0, z_0)) \in \Lambda_1 \times
\Lambda_0$, the corresponding nondegenerate critical point $x$ of $\Delta$
has  shifted index $\tilde i(x)$ given by
$$\tilde i(x) = i_B(V_1) - i_B(V_0) + i_\Gamma((q_0, p_0, z_1), (q_0, p_0,
z_0)),$$
where
$(q_0, p_0, z_i) $ lies in the branch $V_i$  of $\Lambda_i$.
\endproclaim

\demo{Proof} Suppose that $x= (q_0, e_1, e_0)$ is a critical point of
$\Delta$.
By fiber preserving
diffeomorphisms, it can be assumed that in a neighborhood $U$ of $(q_0,e_1,
e_0)$,
$$
\Delta(q, \xi_1, \xi_0)|_U =\left( F_1(q, \xi_1) - F_0(q, \xi_0)\right) |_U
= (G_1(q) + J(\xi_1)) - (G_0(q) + H(\xi_0)),$$
and thus it suffices to
show
$$\ind \Delta(x) - \ind \Delta =
\ind(G_1- G_0)(q_0) + \ind J (e_1) + \ind (-H) (e_0).$$
It is easy to
verify that
$$\ind( G_1  - G_0)(q_0) = i_\Gamma((q_0, p_0, z_1), (q_0,
p_0, z_0)).$$ Let $-\Lambda_0 := \{ (q, -p, -z) \:
(q,p,z)
\in \Lambda_0 \}$.  Then $-F_0$ is a \gqi function for $-\Lambda_0$
and $(q_0, -p_0, - z_0)$ lies in a branch $W_0$ of $-\Lambda_0$
with $i_B(W_0) = - i_B(V_0)$.  Since
$(q_0, e_1) \in S^1 \times \Bbb R^{N_1}$ corresponds
to $(q_0, p_0, z_1) \in V_1$, and
$(q_0, e_0) \in S^1 \times \Bbb R^{N_0}$ corresponds
to $(q_0, -p_0, -z_0) \in W_0$, Proposition 5.3 implies that
$$ \ind{J(e_1)} -   i_{Q_1} = i_B(V_1), \quad
\ind{(-H)(e_0)} -   i_{Q_0} = i_B(W_0) = -i_B(V_0), $$
where $i_{Q_1}$ is the index of $J(e_1)$, $i_{Q_0}$ is the index of
$-H(e_0)$.   Since $\ind \Delta =  i_{Q_1} + i_{Q_0}$, the desired
result follows.  \qed
\enddemo


\heading{Polynomial Calculations}  \endheading

In this section, the positive and negative polynomials
are calculated for the standard rational legendrian links, flypes of  these
standard links, and for connect sums of these flypes.

\proclaim{Theorem 6.1}  Consider the
legendrian link $L\! =\! (2h_n,v_{n-1}, \dots, v_1,2h_1)$.
Then
$$\align
\Gamma^-(\lambda)\left[  L\right] &=
h_1 + h_2\lambda^{-v_1} + h_3\lambda^{-v_1-v_2} + \dots +
h_n\lambda^{-v_1-v_2-\dots-v_{n-1}}, \\
\Gamma^+(\lambda)\left[  L \right] &=
\cases  \lambda \cdot \Gamma^-(\lambda)\left[
 L \right], & h_1 \geq 1 \\
 (1 + \lambda)  +  \lambda\cdot \Gamma^-(\lambda)\left[ L \right], & h_1 = 0.
\endcases
\endalign$$
\endproclaim

\demo{Proof}  First consider the case of $h_1 \geq 1$.
It is possible to legendrian isotop $L$
to a position so that it has a
difference function
$\Delta \: S^1 \times \Bbb R^N \to \Bbb R$
 with $2h_1 + \dots + 2h_n$ nondegenerate critical points:
for $i = 1, \dots, n$, $\Delta$ has
 $h_i$ critical points with, by Proposition 5.5, shifted index
$- v_1 - v_2 - \dots - v_{i-1}+1 $ and value
$c_i^+$, and
$h_i$ critical points of shifted index
$- v_1 - v_2 - \dots -v_{i-1}$ and  value
$c_i^-$,
where
$$c_1^- < c_2^- < \dots < c_n^- < 0 < c_n^+ < \dots < c_2^+ < c_1^+.$$
Figure 6.1.1 illustrates one such construction for the link $(2,2,4,1,2)$.
\midinsert
\cl{\epsfxsize.5\hsize\epsfbox{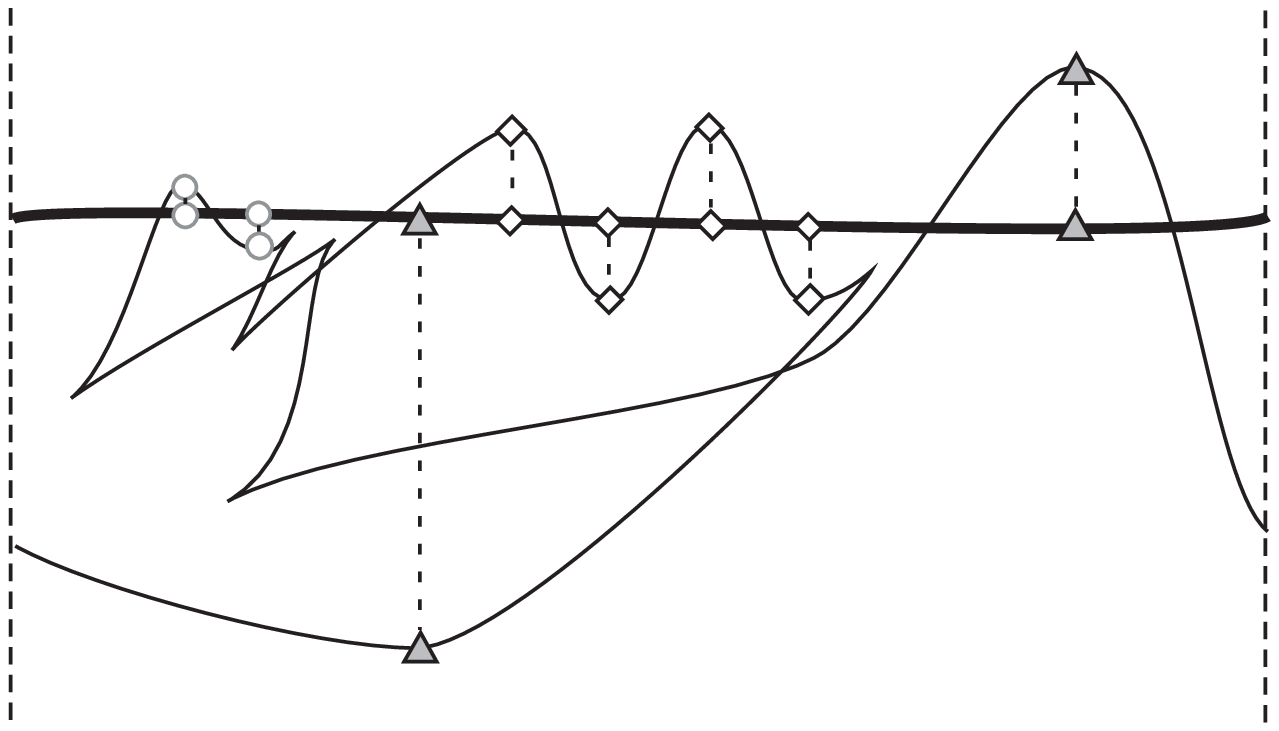}}
\lbotcaption{Figure 6.1.1}  The legendrian  link $(2,2,4,1,2)$   positioned
so that it has a difference function $\Delta$ with $(2+4+2)$ critical
points  which are represented as pairs of points
$((q_0,p_0,z_1), (q_0,p_0,z_0)) \in \Lambda_1 \times \Lambda_0$.  The
shifted indices of these critical points are calculated by Proposition 5.5.
\endcaption
\endinsert

For $n \geq k \geq 2$, choose $a_k, b_k$   so that
$$c_{k-1}^- < a_k < c_k^- < 0 < c_k^+ < b_k < c_{k-1}^+.$$
The claimed calculations will follow immediately from the Positive Integral
Proposition
4.2 if it is shown that
$H_{*}(\Delta^{b_k}, \Delta^{a_k}) = 0$, for all $* \in \Bbb Z$.
 By applying a deformation argument
as in  the proof of Proposition 3.12, it is possible to construct
a
$1$--parameter family of   quadratic at infinity functions
$\Delta_t \: S^1 \times \Bbb R^N \to \Bbb R$ such that
$a_k, b_k$ are noncritical values of $\Delta_t$ for all $t \in [0,1]$,
where $\Delta_0 = \Delta$,
and $\Delta_1$ has no critical points in $[a_k, b_k]$.  Thus by Lemma 3.10,
$ H_{*}(\Delta_0^{b_k}, \Delta_0^{a_k}) =
H_{*}(\Delta_1^{b_k}, \Delta_1^{a_k}) = 0$, for all $*\in \Bbb Z.$
This completes the proof in the case of $h_1 \geq 1$.

In the case of $h_1 = 0$,  it is possible to legendrian isotop $L$
to a position so
that it has a difference function
$\Delta \: S^1 \times \Bbb R^N \to \Bbb R$
that has
 $1$ critical point of shifted index $0$ with   value $c_1^0$,
 $1$ critical point of shifted index $1$ with  value $c_1^1$,
and for $i = 2, \dots, n$,
 $h_i$ critical points of shifted index $- v_1 - v_2 - \dots - v_{i-1} $ with
value $c_i^-$, and
 $h_i$ critical points of index $- v_1 - v_2 - \dots -v_{i-1} + 1$ with
value $c_i^+$,
where
$$c_2^- < \dots < c_n^- < 0 < c_n^+ < \dots < c_2^+ < c_1^0 < c_1^1.$$
For $k = 2, \dots, n$, choose $a_k, b_k$ so that
$$a_2 < c_2^- < c_2^+ < b_2 < c_1^0, \quad
  c_{k-1}^- < a_k < c_k^- < 0 < c_k^+ < b_k < c_{k-1}^+, \  n \geq k
\geq 3.$$
An argument as in the above paragraph proves
$H_{*}(\Delta^{b_k}, \Delta^{a_k}) = 0$ for all $*\in \Bbb Z$, for $k =
2,\dots,
n$.  Thus
  the Zero Integral Proposition 4.3 gives the desired calculation of
$\Gamma^+(\lambda)[L]$ and $\Gamma^-(\lambda)[L]$.  \qed
\enddemo

Next, the positive and negative polynomials will be calculated
for horizontal flypes of a standard rational link.  Recall the
notation
 $(2h_n,v_{n-1}, 2h_{n-1}^{p_{n-1}}, \dots, v_1,$\break $2h_1^{p_1})$,
$p_i \in \{ 0, \dots, 2h_i \}$,
introduced in (1.11).

\proclaim{Theorem 6.2}   Consider the  legendrian
link
$
L = (2h_n,v_{n-1}, 2h_{h-1}^{p_{n-1}}, \dots, v_1,$\break $2h_1^{p_1}).
$
  For $j = 2, \dots, n-1$, let
$\sigma(j) = 1 + \sum_{i=1}^j p_i \mod 2$.
Then
$$\align
\Gamma^-(\lambda)\left[   L  \right] &= h_1 + \sum_{i=2}^n
h_i\lambda^{(-1)^{\sigma(1)}v_1 + (-1)^{\sigma(2)} v_2+ \dots +
(-1)^{\sigma(i-1)}v_{i-1}} , \\
\Gamma^+(\lambda)\left[  L  \right] &=
\cases
\lambda \cdot  \Gamma^-(\lambda)\left[ L \right], & h_1 \geq 1 \\
(1+\lambda) +  \lambda \cdot  \Gamma^-(\lambda)\left[L \right], & h_1 = 0.
\endcases
\endalign$$
\endproclaim

\demo{Proof}  The claim will follow using arguments as in the
proof of Theorem 6.1 if it is shown that for
$L = \Lambda_1 \amalg \Lambda_0 = \left(2h_n,v_{n-1}, 2h_{h-1}^{p_{n-1}},
\dots,
v_1,2h_1^{p_1}\right)
$,  there exists a difference function $\Delta$ with
 $2h_1 + 2h_2 + \dots + 2h_n$ nondegenerate critical points, where for $k =
2, \dots, n$,
the $2h_k$ critical points correspond to $2h_k$ pairs of points on  branches
$W_1^k \subset \Lambda_1$, $W_0^k \subset \Lambda_0$ with branch
indices
$$
i_B(W_1^k) - i_B(W_0^k) = (-1)^{\sigma(1)}v_1 + (-1)^{\sigma(2)} v_2+ \dots +
(-1)^{\sigma(k-1)}v_{k-1},
$$
for $\sigma(j) = 1 +  \sum_{i=1}^j p_i \mod 2$.
As seen in the proof of Theorem 6.1, there exists such a difference
function  for
$L_0 = \left(2h_n,v_{n-1}, 2h_{h-1}^{p_{n-1}}, \dots, v_1,2h_1^{p_1}\right) $,
when $p_k = 0$ for all $k$.
For arbitrary $p_{n-1}, \dots, p_1$,
assume there exists such a difference
function for $L = \left(2h_n,v_{n-1}, 2h_{h-1}^{p_{n-1}},
\dots, v_1,2h_1^{p_1}\right)$.  Let
$\sigma(\ell) =  1 + \sum_{k=1}^\ell p_k \mod 2$.
Consider $L^\prime$ which differs from $L$
by one additional horizontal flype:
$L^\prime = \left(2h_n,v_{n-1}, 2h_{h-1}^{w_{n-1}}, \dots,
v_1,2h_1^{w_1}\right)$, $\exists j \: w_k = p_k$, for $k \neq j$,
and $w_j = p_j + 1$.
Notice that
$$\sigma^\prime(\ell) := 1 + \sum_{k=1}^\ell w_k \mod 2
\equiv_2
\cases
\sigma(\ell), & \ell \leq j-1 \\
\sigma(\ell) + 1, & \ell \geq j.
\endcases$$
As in the case of $L$,
there exists a difference function $\Delta^\prime$ for $L^\prime$ with
$2h_1 + \dots + 2h_n$ critical points; however,  now, due of the flype
 in the $h_j$ term, there is a change in the indices of the branches
containing points in the pairs associated to the
terms $2h_{j+1}, \dots, 2h_n$.
Figure 6.2.1 illustrates such a construction for
horizontal flypes of $(2,1,2,2,2)$.
More precisely, for $k = 2, \dots, n$,
the $2h_k$ critical points are associated to points on  branches
$Y_1^k \subset \Lambda_1^\prime$, $Y_0^k \subset \Lambda_0^\prime$ with branch
indices
$$\align
i_B(Y_0^k) - i_B(Y_1^k) &= (-1)^{\sigma(1)}v_1 + (-1)^{\sigma(2)} v_2+ \dots +
(-1)^{\sigma(j-1)}v_{j-1} + \\
&\text{\hglue-.3in}(-1)^{\sigma(j) + 1}v_{j} +
(-1)^{\sigma(j+1) + 1}v_{j+1} + \dots  +
(-1)^{\sigma(k-1) + 1}v_{k-1} \\
& = (-1)^{\mu(1)}v_1 + (-1)^{\mu(2)} v_2+ \dots +
(-1)^{\mu(k-1)}v_{k-1},
\endalign$$
where
$$\mu(\ell) \equiv_2 \cases
\sigma(\ell), & \ell \leq  j-1 \\
\sigma(\ell) + 1, & \ell \geq j,
\endcases$$
as desired.\qed
\enddemo

\midinsert
 \cl{\epsfxsize\hsize\epsfbox{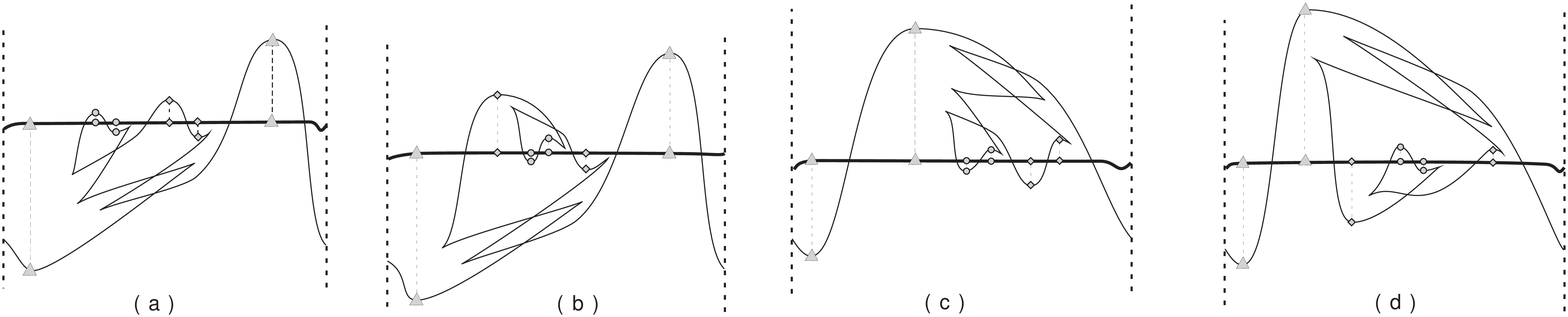}}
\lbotcaption{Figure 6.2.1}  Links that differ by horizontal
flypes: (a) $(2,1,2,2,2)$, (b)
$(2,1,2^1,2,2)$, (c) $(2,1,2,2,2^1)$, (d) $(2,1,2^1,2,2^1)$.  If $L$ and
$L^\prime$ differ by a horizontal flype from the $j^{th}$ component,
then there will be a change in the indices of the critical points associated to the terms
$h_{j+1}, \dots, h_{n}$. 
\endcaption
\endinsert

\remark{\bf Remark 6.3}  Applying vertical flypes to a link will leave the
 positive and negative homology polynomials unchanged.
This follows since if
$$L  = \left(2h_n,v_{n-1}^{q_{n-1}}, 2h_{h-1}^{p_{n-1}},
\dots, v_1^{q_1},2h_1^{p_1}\right),$$
 $$L^\prime =  \left(2h_n,v_{n-1}^{w_{n-1}}, 2h_{h-1}^{p_{n-1}},
\dots, v_1^{w_1},2h_1^{p_1}\right), $$
where for some $j$,  $ w_k = q_k$, for $ k \neq j$,  and
$w_j = q_j + 1$,
then there exist difference functions $\Delta$ and $\Delta^\prime$
for $L$ and $L^\prime$ with the same critical values and
the {\it same} shifted indices. Similar to the situation in
the  proof of
Theorem 6.2, a vertical flype from the $j^{th}$ component
will affect the indices of
the  branches that contain points of the pairs
associated to the
terms $2h_{j+1}, \dots, 2h_n$.  Although the branch
index associated to each point
in the pair  will change, the difference
  between the branch indices of this
pair is unchanged.  
\endremark

 \proclaim{Theorem 6.4}  Consider the legendrian links
$$
 L_1\! = (2h_n, v_{n-1}, 2h_{n-1}^{p_{n-1}}, \dots, v_1, 2h_1^{p_1}),\
 L_2\! = (2k_m, u_{m-1}, 2k_{m-1}^{w_{m-1}}, \dots, u_1, 2k_1^{w_1}).
$$
  Then
$$\align
\Gamma^-(\lambda)[L_1 \# L_2] &= \Gamma^-(\lambda)[L_1] +
\Gamma^-(\lambda)[L_2];\\
\Gamma^+(\lambda)[L_1 \# L_2] &=
\cases
\Gamma^+(\lambda)[L_1] + \Gamma^+(\lambda)[L_2], & h_1, k_1 \geq 1 \\
\Gamma^+(\lambda)[L_1] + \Gamma^+(\lambda)[L_2] - (1+\lambda), & \text{ else.}
\endcases
\endalign$$
\endproclaim

\demo{Proof} Let $: L_i :$ denote the tangles whose closure give $L_i$.
It is possible to isotop the tangles $: L_1 :$ and $:L_2:$ to the
 configurations similar to those
used to calculate $\Gamma^+(\lambda)[L_i]$ and  $\Gamma^-(\lambda)[L_i]$
but ``scaled" so there exists an $M
\in \Bbb N$ such that
 all critical values associated to
$:L_1:$ are contained in $[-M, M]$, while all critical values
associated to
$:L_2:$ are contained in $(-\infty, -M) \cup (M, \infty)$.
Figure 6.4.1 illustrates this construction when $p_1 = 0 = w_1$.
\midinsert
\cl{\epsfxsize.6\hsize\epsfbox{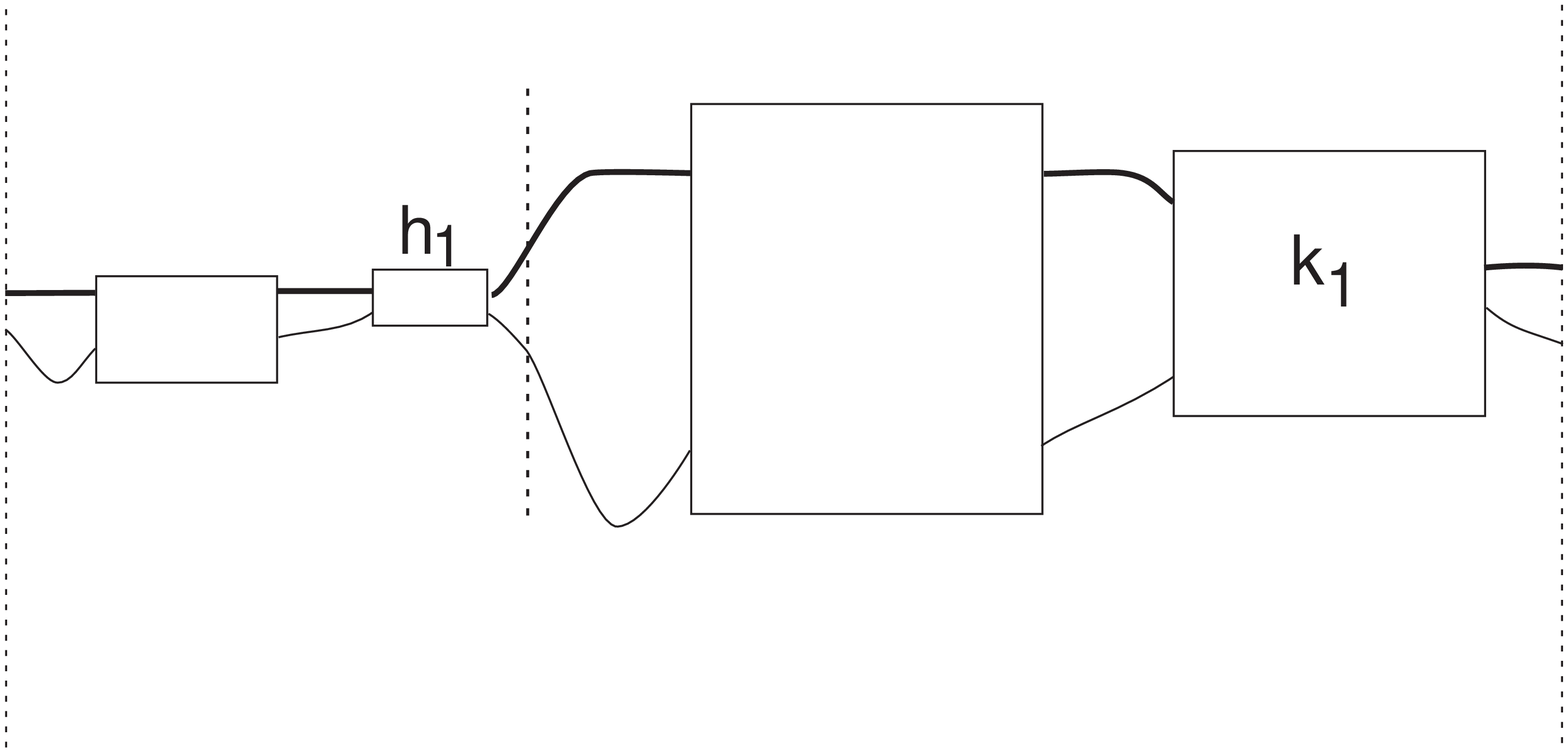}}
\lbotcaption{Figure 6.4.1} The connect sum of $L_1$ and $L_2$, when
$p_1 = w_1 = 0$,
conveniently scaled.   
\endcaption
\endinsert

When $h_1, k_1 \geq 1$, this procedure gives rise to a difference
function $\Delta$ with
$2(h_1 + \dots + h_n + k_1 + \dots +k_m)$ critical points:
for $j = 1, \dots, n$, $ \iota= 1, \dots, m$, $\Delta$ has
 $h_j$ critical points of shifted index
$i_j+1 $ and value $c_j^+$, $h_j$ critical points of shifted index
$i_j$ and  value $c_j^-$,
$k_\iota$ critical points of shifted index
$ \ell_\iota+1 $ and value
$d_\iota^+$, and $k_\iota$ critical points of shifted index
$ \ell_\iota$ and  value
$d_\iota^-$,
where
$$\align
d_1^- < d_2^- < \dots < d_m^- < c_1^- < c_2^-& < \dots < c_n^-  < 0  \\
 0 < &c_n^+ < \dots < c_2^+ < c_1^+< d_m^+ < \dots < d_2^+ < d_1^+.
\endalign$$
Arguments as in the proofs of Theorems 6.1 and 6.2 show that the  hypotheses
of Proposition 4.2 are satisfied.  Thus
$$\align
\Gamma^-(\lambda)[L_1 \# L_2] &= \sum_{j=1}^n h_j \lambda^{i_j} +
\sum_{\iota=1}^m k_\iota \lambda^{\ell_\iota}
= \Gamma^-(\lambda)[L_1]  + \Gamma^-(\lambda)[L_2], \\
\Gamma^+(\lambda)[L_1 \# L_2] &= \sum_{j=1}^n h_j \lambda^{i_j+1}+
\sum_{\iota=1}^m k_\iota
\lambda^{\ell_\iota+1} = \Gamma^+(\lambda)[L_1]  + \Gamma^+(\lambda)[L_2].
\endalign$$

For the case $h_1=0$ or $k_1=0$, since $L_1 \# L_2$ is equivalent to
$L_2 \# L_1$, it can be assumed $h_1 = 0$.  In this case, two positive critical
points that were necessary so that the tangles $: L_i :$ closed to  links
can be
eliminated.  Once  $L_1 \# L_2$ is positioned so that
these points are eliminated, it is possible to
construct a
$\Delta$ that satisfies all the hypotheses of either Proposition 4.2 or
Proposition 4.3.
More precisely, when $k_1 \geq 1$, it is possible to legendrian isotop the
link $L_1\# L_2$ to a position so that it has a generating function
$\Delta$
 with  $2(h_2 + \dots + h_n  + k_1 + \dots + k_m)$ critical points:
for $j = 2, \dots, n$, $\iota = 1, \dots, m$, $\Delta$ has
$h_j$ critical points of shifted index
$i_j+1 $ and value
$c_j^+$, $h_j$ critical points of shifted index
$i_j$ and  value
$c_j^-$,
$k_\iota$ critical points of shifted index
$ \ell_\iota+1 $ and value
$d_\iota^+$, and
$k_\iota$ critical points of shifted index
$\ell_\iota$ and  value
$d_\iota^-$,
where
$$\align
d_1^- < d_2^- < \dots < d_m^- < c_2^- <& \dots < c_n^- < 0  \\
 & 0 < c_n^+ < \dots < c_2^+ < d_m^+ < \dots < d_2^+  <  d_1^+.
\endalign$$
All hypotheses of Proposition 4.2 are satisfied, and thus
$$\align
\Gamma^-(\lambda)[L_1 \# L_2] &= \sum_{j=2}^n h_j \lambda^{i_j} +
\sum_{\iota=1}^m k_\iota \lambda^{\ell_\iota}
= \Gamma^-(\lambda)[L_1]  + \Gamma^-(\lambda)[L_2], \\
\Gamma^+(\lambda)[L_1 \# L_2] &= \sum_{j=2}^n h_j \lambda^{i_j+1}+
\sum_{\iota=1}^m k_\iota
\lambda^{\ell_\iota+1}\\ &= \left( \Gamma^+(\lambda)[L_1]- (1+\lambda) \right)
+ \Gamma^+(\lambda)[L_2] .
\endalign$$
When $k_1 = 0$, it is possible to legendrian isotop the
link $L_1\# L_2$  to a position so that it has a generating function
$\Delta$
with  $2(h_2 + \dots + h_n + k_2 + \dots + k_m) + 2$ critical points: $\Delta$
has $1$ critical point of shifted index $0$ with value $d_1^0$,
$1$ critical point of shifted index $1$ with value $d_1^1$,
for $j = 2, \dots, n$, $ \iota= 2, \dots, m$, $\Delta$ has
 $h_j$ critical points of shifted index
$i_j+1 $ and value $c_j^+$,  $h_j$ critical points of shifted index
$ i_j$ and  value
$c_j^-$,
$k_\iota$ critical points of shifted index
$ \ell_\iota+1 $ and value
$d_\iota^+$,
 and $k_\iota$ critical points of shifted index
$ \ell_\iota$ and  value
$d_\iota^-$, where
$$\align
d_2^- < \dots < d_m^- < c_2^- <& \dots < c_n^- < 0  \\
&0 < c_n^+ < \dots < c_2^+ < d_m^+ < \dots < d_2^+  < d_1^0 <  d_1^1.
\endalign$$
All hypotheses of Proposition 4.3 are satisfied, and thus
$$\align
\Gamma^-(\lambda)[L_1 \# L_2] &= \sum_{j=2}^n h_j \lambda^{i_j} +
\sum_{\iota=2}^m k_\iota \lambda^{\ell_\iota}
= \Gamma^-(\lambda)[L_1]  + \Gamma^-(\lambda)[L_2], \\
 \Gamma^+(\lambda)[L_1 \# L_2] &= (1+\lambda)  + \sum_{j=2}^n h_j
\lambda^{i_j+1} +
  \sum_{\iota=2}^m k_\iota \lambda^{\ell_\iota+1}  =
\Gamma^+(\lambda)[L_1]\\&\text{\hglue 1.8in}  + (\Gamma^+(\lambda)[L_2]-
(1+\lambda) ) .\qquad\qed
\endalign$$
\enddemo


\heading{Applications}\endheading

In this section, the polynomial calculations of Section 6 will
be applied to  show that ``most" of the legendrian
links $L_n = (2h_n,v_{n-1}, \dots, 2h_1)$ are  ordered.  In addition,
by analyzing the polynomials resulting from
the swap and the flype operations, lower bounds will be given for  the
number of
different minimal legendrian representations
of a given topological link type.

The following proposition shows that there is a simple
relation  between the  polynomials of
a rational link and its swap.

\proclaim{Proposition 7.1}  Let
$L = (2h_n, v_{n-1}, 2h_{n-2}^{p_{n-2}}, \dots, h_2^{p_2}, v_1, 2h_1^{p_1})$,
and let $\overline L$ denote the swap of $L$.  Then
$$\Gamma^-(\lambda)[\overline{L}] =
\cases
\overline{\Gamma^-(\lambda)[L]}, & h_1 \geq 1,\\
(1+\lambda) + \overline{\Gamma^-(\lambda)[L]}, & h _1 = 0,
\endcases
$$
where $\overline{\Gamma^-(\lambda)[L]}$ denotes the palindrome of
$\Gamma^-(\lambda)[L]$.
\endproclaim

\demo{Proof}  This follows easily from Theorem  3.16 and Theorem 6.2. For
$h_1
\geq 1$,
$$\lambda\cdot\Gamma^-(\lambda)[L] = \Gamma^+(\lambda)[L] = \lambda \cdot
\overline{\Gamma^-(\lambda)[\overline{L}]}.$$ For $h_1 = 0$,
$$(1+\lambda) + \lambda\cdot\Gamma^-(\lambda)[L] = \Gamma^+(\lambda)[L] =
\lambda \cdot
\overline{\Gamma^-(\lambda)[\overline{L}]}.$$
This implies, by examining palindromes, that
$$ \lambda^{-1} \Gamma^-(\lambda)[\overline L] = 1 + \lambda^{-1} +
\lambda^{-1}\overline{\Gamma^-(\lambda)[L]},$$ and the desired result
follows.  \qed
\enddemo

\proclaim{Theorem 7.2}  Consider the legendrian link
$L = (2h_n,v_{n-1}, 2h_{n-1}, \dots, v_1,$\break$2h_1)$.
If $h_1 \geq 1$,  then $L$ is ordered iff $L \neq (2h_1)$.
If $h_1 = 0$, then $L$ is ordered.
\endproclaim

\demo{Proof} First suppose that $h_1 \geq 1$. When   $L = (2h_1)$,   it
is easy to explicitly check that $L$ is not ordered.
Suppose
$L = (2h_n, v_{n-1}, \dots, v_1, 2h_1)$, $n \geq 2$.
By Proposition 7.1, to show $L$ is ordered,
it suffices to prove that $\Gamma^-(\lambda)[L]$ is not palindromic.  By
Theorem 6.1,
$$ \Gamma^-(\lambda)[L] =
 h_1 + h_2\lambda^{-v_1} +h_3\lambda^{-v_1-v_2}+ \dots +
h_n\lambda^{-v_1-v_2-\dots-v_{n-1}}, $$ and, since $n \geq 2$,
$\Gamma^-(\lambda)[L]$ is
not palindromic.

Next suppose $h_1 = 0$.
By Proposition 7.1, it suffices to verify that
$(1 + \lambda) + \lambda\cdot~\Gamma^-(\lambda)[L] \neq \lambda \cdot
\overline{\Gamma^-(\lambda)[L]}$.
By Theorem 6.1,
$$\lambda \cdot \overline{\Gamma^-(\lambda)[L]} = \sum_{m=-\infty}^\infty
a_m \lambda^m, \qquad
\text{with } a_1 = 0,
$$ while
$\lambda\cdot \Gamma^-(\lambda)[L] = \sum_{m = -\infty}^\infty b_m t^m$,
with  $ b_1 = 0,$
and thus
$$(1 + \lambda) + \lambda \cdot \Gamma^-(\lambda)[L] = \sum_{m =
-\infty}^\infty c_m \lambda^m, \qquad
\text{with } c_1 = 1.$$ Thus   $L$ must be ordered.
\qed
\enddemo

Distinct minimal legendrian
  versions
of the topological link
$$L_n = (2h_n, v_{n-1}, 2h_{n-1} \dots, 2h_2, v_1, 2h_1), \quad
h_n,\dots, v_1 \geq 1,\  h_1 \geq 1,$$
will now be enumerated.
For $h_1 \geq 1$, there are potentially $2\cdot 2^{n-1}$ legendrian
versions arising from the swap and the flype operations that can
be distinguished by the polynomials.  However,
the following proposition shows that the polynomials cannot
distinguish all these swaps and flypes:
the swap
operation always produces a link with the same polynomials as some
flype.

\proclaim{Proposition 7.3}  For $h_1 \geq 1$, $p_1 \in \{0, \dots,
2h_1-1\}$, and
 $p_i \in \{0, \dots, 2h_i\}$ when $i \geq 2$,
consider the minimal legendrian links
$$\align
 L_1 &= (2h_n, v_{n-1}, 2h_{n-1}^{p_{n-1}}, \dots, 2h_2^{p_2}, v_1,
2h_1^{p_1}), \\
   L_2 &= (2h_n, v_{n-1}, 2h_{n-1}^{p_{n-1}}, \dots, 2h_2^{p_2}, v_1,
2h_1^{p_1+1}),
\endalign
$$
Let $\overline{L_1}$ be the  swap of $L_1$.  Then
$\Gamma^-(\lambda)[\overline{L_1}] = \Gamma^-(\lambda)[L_2],$ and
$\Gamma^+(\lambda)[\overline{L_1}] = \Gamma^+(\lambda)[L_2].
$
\endproclaim

\demo{Proof}  It suffices to prove that
$\Gamma^-(\lambda)[\overline{L_1}] = \Gamma^-(\lambda)[L_2]$.  By
Proposition 7.1 and Theorem 6.2, for
$\sigma(j) = 1 + \sum_{i=1}^j p_i \mod 2$,
$$\aligned
\Gamma^-(\lambda)[\overline{L_1}] &= \overline{\Gamma^-(\lambda)[L_1]}
= \overline{h_1 + \sum_{i=2}^n
h_i \lambda^{(-1)^{\sigma(1)}v_1 + (-1)^{\sigma(2)} v_2+ \dots +
(-1)^{\sigma(i-1)}v_{i-1}}}\\
&= {h_1 + \sum_{i=2}^n
h_i \lambda^{(-1)^{\sigma(1)+1}v_1 + (-1)^{\sigma(2)+1} v_2+ \dots +
(-1)^{\sigma(i-1)+1}v_{i-1}}} \\
&= \Gamma^-(\lambda)[L_2].\text{\hglue3.05in} \qed
\endaligned
$$
\enddemo

Thus   the swap and flype operations give at most $2^{n-1}$
legendrian versions of $L_n$ that can be distinguished by the
polynomials.  In fact,
 for $n =3$,
there are often at least $4 = 2^{3-1}$ versions of $L_3$ distinguishable
by the polynomials.

\proclaim{Theorem 7.4}  Consider the topological link
$$ L_3 = (2h_3, v_2, 2h_2, v_1, 2h_1), \qquad  h_3, h_2, v_2, v_1 \geq
1,\quad h_1 \geq 0.$$
If $h_1 = 0$, then there are at least $4$ minimal legendrian versions of
$L_3$.
If $h_1 \geq 1$ and either
$v_2 \neq 2v_1$ or $h_2 \neq h_3$,
then there are at least $4$ minimal legendrian versions of   $L_3$.
\endproclaim

\demo{Proof}  For $h_1 = 0$, consider
$$ \gather L_3^0 = (2h_3, v_2, 2h_2^0, v_1, 0),\quad L_3^1 = (2h_3, v_2,
2h_2^1, v_1, 0), \\
L_3^2 = \overline{(2h_3, v_2, 2h_2^0, v_1, 0)}, \quad
L_3^3 = \overline{(2h_3, v_2, 2h_2^1, v_1, 0)}.
\endgather$$
By Theorem 6.2 and Proposition 7.1,
$$\gather
\Gamma^-(\lambda)[L_3^0] = h_2 \lambda^{-v_1} + h_3\lambda^{-v_1 -v_2}, \quad
\Gamma^-(\lambda)[L_3^1] = h_2\lambda^{-v_1} + h_3\lambda^{-v_1 +v_2}, \\
\Gamma^-(\lambda)[L_3^2] = 1 + \lambda + h_2\lambda^{+v_1} +
h_3\lambda^{+v_1 +v_2},\\
\Gamma^-(\lambda)[L_3^3] = 1 + \lambda + h_2\lambda^{+v_1} +
h_3\lambda^{+v_1 -v_2}.
\endgather
$$
It is easy to verify that for all choices of $h_i, v_i \geq 1$, these must be
distinct polynomials, and thus $L_3^i$, $i = 0, \dots 3$, are legendrian
distinct.

For $h_1 \geq 1$, consider
$$ \gather L_3^{(0,0)} = (2h_3, v_2, 2h_2^0, v_1, 2h_1^0),\quad
L_3^{(1,0)} =  (2h_3, v_2, 2h_2^1, v_1, 2h_1^0), \\
L_3^{(0,1)} =  (2h_3, v_2, 2h_2^0, v_1, 2h_1^1 ), \quad
L_3^{(1,1)} = (2h_3, v_2, 2h_2^1, v_1, 2h_1^1).
\endgather$$
By Theorem 6.2,
$$
\gather
\Gamma^-(\lambda)[L_3^{(0,0)}] = h_1 + h_2\lambda^{-v_1} + h_3\lambda^{-v_1
-v_2}, \\
\Gamma^-(\lambda)[L_3^{(1,0)}] = h_1 + h_2\lambda^{-v_1} + h_3\lambda^{-v_1
+v_2}, \\
\Gamma^-(\lambda)[L_3^{(0,1)}] = h_1 + h_2\lambda^{+v_1} + h_3\lambda^{+v_1
+v_2}, \\
\Gamma^-(\lambda)[L_3^{(1,1)}] = h_1 + h_2\lambda^{+v_1} + h_3\lambda^{+v_1
-v_2}.
\endgather
$$
The condition $v_2 \neq 2v_1$  or $h_2 \neq h_3$ implies that all  these
polynomials are
distinct. Thus there are at least $4$ distinct legendrian versions of
$L_3$.  \qed
\enddemo

The following condition on $h_i$ will guarantee that all
the flypes  have distinct $\Gamma^-(\lambda)$ polynomials. Such
sets arise in Additive Number Theory; see \cite{G}.

\definition{Definition 7.5}  A set $\{h_1, \dots, h_n\}$ of integers is
said to have {\sl distinct subset sums} if the sums of all its $2^n$
subsets are
distinct. Such a set will be abbreviated as a {\sl \dss set}.
\enddefinition

   It is easy to verify that $\{1,2\}$, $\{ 1,2,4 \}$ and
$\{ 2,3,4,8\}$ are \dss sets, while $\{ 1,2, 3\}$ is not. In general,
$\{ 2^i \: 0 \leq i \leq k\}$ is a \dss set of order $k+1$.
There is an Erd\"os prize associated to finding the largest order of
a \dss set with
entries positive and bounded above by $2^k$.

\proclaim{Theorem 7.6}  For $n \geq 4$,
consider the topological link
 $$L_n = (2h_n, v_{n-1}, 2h_{n-1}^0, \dots, 2h_2^0, v_1, 2h_1^0),
\quad h_n, v_{n-1}, \dots, v_1 \geq 1, h_1 \geq 0.$$
If $h_1 \geq 1$, assume $\{h_1, h_2, \dots, h_n\}$ form a  \dss set
of order $n$, while if $h_1 =
0$, assume $\{1, h_2, \dots, h_n\}$ form a \dss set of order $n$.    Then
there exist at least $2^{n-1}$ legendrian versions of $L_n$.
\endproclaim

\demo{Proof}  For the case where $h_1 \geq 1$, consider the
links
$$ (2h_n, v_{n-1}, 2h_{n-1}^{p_{n-1}}, \dots, 2h_2^{p_2}, v_1, 2h_1^{p_1}),
\qquad
p_{n-1}, \dots, p_1 \in \Bbb Z_2.$$
 It will be shown
that  distinct choices of $(p_{n-1}, \dots, p_1) \in \Bbb Z_2 \times \dots
\times \Bbb Z_2$
give rise to   legendrian links with distinct $\Gamma^-(\lambda)$ polynomials.
Let
$$\align
L_1 &= (2h_n, v_{n-1}, 2h_{n-1}^{p_{n-1}}, \dots, 2h_2^{p_2}, v_1,
2h_1^{p_1}), \\
L_2 &= (2h_n, v_{n-1}, 2h_{n-1}^{q_{n-1}}, \dots, 2h_2^{q_2}, v_1, 2h_1^{q_1})
\endalign$$
and suppose that $\Gamma^-(\lambda)[L_1] = \Gamma^-(\lambda)[L_2]$.
By Theorem 6.2,
	$$\Gamma^-(\lambda)[L_1] = h_1 + h_2\lambda^{s_2} +   \dots + h_n
\lambda^{s_n}, \quad
\Gamma^-(\lambda)[L_2] = h_1 + h_2\lambda^{m_2}   + \dots + h_n\lambda^{m_n}$$
for
$$
\align
s_k &= (-1)^{\sigma(1)} v_1 + (-1)^{\sigma(2)}v_2 + \dots +
(-1)^{\sigma(k-1)}v_{k-1},\\
m_k &= (-1)^{\mu(1)} v_1 + (-1)^{\mu(2)}v_2+ \dots +
(-1)^{\mu(n-1)}v_{n-1}, \qquad k = 2,\dots, n,
\endalign
$$
where $\sigma(s) =1 + \sum_{t=1}^s p_t \mod 2$,
$\mu(s) = 1 + \sum_{t=1}^s q_t \mod 2$.
If  it is shown that $s_k = m_k$, $\forall k$, then   $\sigma(k) =
\mu(k)$, $\forall k$, and thus
$p_k = q_k$, $\forall k$.

First   $\Gamma^-(\lambda)[L_1]$ and $\Gamma^-(\lambda) [L_2]$ will be
rewritten in terms of distinct powers of
$\lambda$. Choose $I_0, \dots, I_N$ with $I_0 < \dots < I_N$ to be the
distinct elements of
$\{ 0, s_2, \dots, s_n \} =  \{ 0, m_2, \dots, m_n \}$.  Then
$$ \Gamma^-(\lambda)[L_1] = \sum_{\alpha = 0}^N \sum_{k=1}^{K_\alpha}
h_{i_k^\alpha} \lambda^{I_\alpha}, \qquad
\Gamma^-(\lambda)[L_2] = \sum_{\alpha = 0}^N \sum_{k=1}^{L_\alpha}
h_{j_k^\alpha} \lambda^{I_\alpha},
$$
where, for all $\alpha$,
 $h_{i_1^\alpha} < h_{i_2^\alpha} < \dots < h_{{i^\alpha_{K_\alpha}}}$, and
$h_{j_1^\alpha} < h_{j_2^\alpha} < \dots < h_{{j^\alpha_{L_\alpha}}}$.
Since $\{h_1, \dots, h_n\}$ are a  \dss set, $\Gamma^-(\lambda)[L_1] =
\Gamma^-(\lambda)[L_2]$ implies $K_\alpha = L_\alpha$, for all $\alpha$,
and that ${i_k^\alpha} =
{j_k^\alpha}$, for all $\alpha$ and $k$. This implies $s_k = m_k$, as
desired.

For the case where $h_1 = 0$, consider the
links
$$ (2h_n, v_{n-1}, 2h_{n-1}^{p_{n-1}}, \dots, 2h_2^{p_2}, v_1, 0), \qquad
\overline{(2h_n, v_{n-1}, 2h_{n-1}^{p_{n-1}}, \dots, 2h_2^{p_2}, v_1, 0) },$$
where $p_{n-1}, \dots, p_2 \in \Bbb Z_2.$
 It will be shown
that  the $2^{n-2}$ choices of $(p_{n-1}, ..., p_2)$\break$ \in \Bbb Z_2 \times
\dots \times \Bbb Z_2$
give rise to $2 \cdot 2^{n-2}$  legendrian links with distinct
$\Gamma^-(\lambda)$ polynomials.
By Theorem 6.2 and Propositions 7.1 and 7.3,
$$\align
\Gamma^-(\lambda)&[(2h_n, v_{n-1}, 2h_{n-1}^{p_{n-1}}, \dots, 2h_2^{p_2},
v_1, 0)] \\
&=\Gamma^-(\lambda)[(2h_n, v_{n-1}, 2h_{n-1}^{p_{n-1}}, \dots, 2h_2^{p_2},
v_1, 2^0)] - 1,\\
\Gamma^-(\lambda)&[\overline{(2h_n, v_{n-1}, 2h_{n-1}^{p_{n-1}}, \dots,
2h_2^{p_2}, v_1, 0)}] \text{\hskip 1in\ }\\
&=\overline{\Gamma^-(\lambda)[ (2h_n, v_{n-1}, 2h_{n-1}^{p_{n-1}}, \dots,
2h_2^{p_2}, v_1, 0)  ]} + 1 +
\lambda \\ &= \overline{\Gamma^-(\lambda)[ (2h_n, v_{n-1},
2h_{n-1}^{p_{n-1}}, \dots, 2h_2^{p_2}, v_1, 2^0)  ]
- 1} + 1 + \lambda\\ &= \overline{\Gamma^-(\lambda)[ (2h_n, v_{n-1},
2h_{n-1}^{p_{n-1}}, \dots, 2h_2^{p_2},
v_1, 2^0)  ] } + \lambda \\ &= \Gamma^-(\lambda)[(2h_n, v_{n-1},
2h_{n-1}^{p_{n-1}}, \dots, 2h_2^{p_2}, v_1,
2^1)] +
\lambda.
\endalign$$
Thus to show that these $2^{n-1}$ polynomials are distinct, the following three
statements will be proven.
\roster
\item If $(p_{n-1}, \dots, p_2) \neq (q_{n-1}, \dots, q_2) $, then
$$\align
\Gamma^-(\lambda)[(2h_n, v_{n-1}, 2h_{n-1}^{p_{n-1}}, \dots, &2h_2^{p_2},
v_1, 2^0)] - 1 \neq\\
\Gamma^-(\lambda)[(&2h_n, v_{n-1}, 2h_{n-1}^{q_{n-1}}, \dots, 2h_2^{q_2},
v_1, 2^0)] - 1;
\endalign$$
\item If $(p_{n-1}, \dots, p_2) \neq (q_{n-1}, \dots, q_2) $, then
$$\align
\Gamma^-(\lambda)[(2h_n, v_{n-1}, 2h_{n-1}^{p_{n-1}}, \dots, &2h_2^{p_2},
v_1, 2^1)] + \lambda \neq\\
\Gamma^-(\lambda)[(&2h_n, v_{n-1}, 2h_{n-1}^{q_{n-1}}, \dots, 2h_2^{q_2},
v_1, 2^1)] + \lambda;
\endalign$$
\item For all $(p_{n-1}, \dots, p_2), (q_{n-1}, \dots, q_2)$,
$$\align \Gamma^-(\lambda)[(2h_n, v_{n-1}, 2h_{n-1}^{p_{n-1}}, \dots,
&2h_2^{p_2}, v_1, 2^0)] - 1 \neq\\
\Gamma^-(\lambda)[(&2h_n, v_{n-1}, 2h_{n-1}^{q_{n-1}}, \dots, 2h_2^{q_2},
v_1, 2^1)] + \lambda.
\endalign$$
\endroster
Statements (1) and (2) follow from Theorem 6.2.   To verify statement (3),
suppose that
there exist
$(p_{n-1}, \dots, p_2), (q_{n-1}, \dots, q_2)$ such that
$$\align \Gamma^-(\lambda)[(2h_n, v_{n-1}, 2h_{n-1}^{p_{n-1}}, \dots,
&2h_2^{p_2}, v_1, 2^0)] - 1 = \\
\Gamma^-(\lambda)[(&2h_n, v_{n-1}, 2h_{n-1}^{q_{n-1}}, \dots, 2h_2^{q_2},
v_1, 2^1)] + \lambda.
\endalign$$
By writing the polynomials in terms of distinct powers of $\lambda$,
 $$\align
\Gamma^-(\lambda)&[(2h_n, v_{n-1}, 2h_{n-1}^{p_{n-1}}, \dots, 2h_2^{p_2},
v_1, 2^0)] - 1 = \sum a_i t^i, \\
&\text{ where }a_0 = h_{i_1} + \dots + h_{i_{K}}, \qquad
h_{i_1}, \dots, h_{i_K} \in \{ h_2, \dots, h_n\},
\endalign$$ while
$$\align
\Gamma^-(\lambda)&[(2h_n, v_{n-1}, 2h_{n-1}^{q_{n-1}}, \dots, 2h_2^{q_2},
v_1, 2^1)] + \lambda
= \sum b_i \lambda^i, \\
&\text{ where } b_0 = 1 + h_{j_1} + \dots + h_{j_{L}}, \qquad
h_{j_1}, \dots, h_{j_L} \in \{ h_2, \dots, h_n\}.
\endalign$$
The assumption that these polynomials are equal
contradicts the hypothesis that $\{1, h_2, \dots, h_n\}$ are a  \dss set of
order $n$.
Thus statement (3) is true, and it follows that the $2^{n-1}$ polynomials
are distinct.
 \qed
\enddemo

\subheading{\bf Acknowledgements}  During the preparation of this paper, I had
the pleasure of visiting the American Institute of Mathematics in Palo Alto, CA
for their program on contact geometry, and the Institute for Advanced Study
in Princeton, NJ.
I thank both of these institutions for  their
wonderful hospitality.  I also thank my reviewers for numerous useful
comments and
suggestions.

This research was supported in part by NSF grant DMS 9971374 and, while
visiting IAS, NSF grant DMS 97-29992.\eject

\enddocument